%% file: main.tex
\documentclass[preprint,12pt]{elsarticle}
\usepackage[utf8]{inputenc}
\usepackage[left=1.2in,right=1.2in,top=1.2in,bottom=1.5in]{geometry}
\usepackage{amsmath}
\usepackage{libertine}
\usepackage{amsfonts}
\usepackage{color,soul} % for \hl{} highlighting
\newcommand{\com}[1]{}
\usepackage{natbib}
\usepackage{listings}
\usepackage{url}
\usepackage{float}
\usepackage{gensymb}
\usepackage{comment}
\usepackage{threeparttable}
\newcommand{\tablefontsize}{\footnotesize}
\usepackage{arydshln}
\usepackage{hhline}
\usepackage[compact]{titlesec}
\usepackage{adjustbox}
\usepackage{booktabs} 
\usepackage[skip=-1pt]{subcaption}
\usepackage[version=3]{mhchem} 
\usepackage{appendix}
\usepackage{indentfirst}
\usepackage{changepage}
\usepackage[table,xcdraw,svgnames]{xcolor}
\usepackage[colorlinks]{hyperref}
\abstracttitle{Summary}
\RequirePackage{doi}

\graphicspath{{Figures/}}
\journal{Energy \& Environmental Science}

\begin{document}
\begin{frontmatter}
\title{Assessing the Value of Coupling Thermal Energy Storage with Air-Source Heat Pumps for Residential Space Heating in U.S. Cities.}

\affiliation[inst1]{organization={School for Environment and Sustainability, University of Michigan},%Department and Organization
            addressline={440 Church St.},
            city={Ann Arbor},
            state={MI 48109},
            country={USA}}

\affiliation[inst2]{organization={Department of Industrial and Operations Engineering, University of Michigan},%Department and Organization
            addressline={1891 IOE Building 1205, Beal Ave}, 
            city={Ann Arbor},
            state={MI 48109},
            country={USA}}            

\affiliation[inst3]{organization={Department of Mechanical Engineering, University of Michigan},%Department and Organization
            addressline={2350 Hayward St.}, 
            city={Ann Arbor},
            state={MI 48109},
            country={USA}}
            
\author[inst1]{An T. Pham}
\author[inst3]{Bryan Kinzer}
\author[inst1]{Ritvik Jain}
\author[inst3]{Rohini Bala Chandran}
\author[inst1,inst2]{Michael T. Craig}

\begin{abstract}
%% Text of abstract
Widespread air source heat pump (ASHP) adoption faces several challenges that on-site thermal energy storage (TES), particularly thermochemical salt hydrate TES, can mitigate. No techno-economic analyses for salt-hydrate-based TES in residential applications exist. We quantify the residential space heating value of four salt hydrate TES materials - \ce{MgSO4}, \ce{MgCl2}, \ce{K2CO3}, and \ce{SrBr2} - coupled with ASHPs across 4,800 representative households in 12 U.S. cities by embedding salt-hydrate-specific Ragone plots into a techno-economic model of coupled ASHP-TES operations. In Detroit, salt hydrate TES is projected to reduce household annual electricity costs by up to \$241 (8\%). Cost savings from TES can differ by over an order of magnitude between households and salt hydrates. We identify the most promising salt in this study, \ce{SrBr2}, due to its high energy density and low humidification parasitic load. Break-even capital costs of \ce{SrBr2}-based TES range from \$13/kWh to \$17/kWh, making it the only salt hydrate studied to reach and exceed the U.S. Department of Energy's \$15/kWh TES cost target. Sensitivities highlight the importance of variable TES sizing and efficiency losses in the value of TES. 
\end{abstract}

\begin{comment}
    %%Research highlights
\begin{highlights}
\item We value salt hydrate thermal energy storage (TES) coupled with an air-source heat pump for residential space heating in 4,800 U.S. single family homes.
\item We integrate Ragone plots for four salt hydrates into a technoeconomic framework for sizing, operating, and valuing TES.
\item \ce{SrBr2} yields highest TES break-even cost of over \$15/kWh thermal with time-of-use retail rates in cold climate cities, potentially exceeding Department of Energy's \$15/kWh TES cost target.
\item Salt hydrate TES can serve as a powerful demand response tool, capable of reducing annual peak heating load by over 12\% in homes in a cold climate city.
\end{highlights}
\end{comment}

\begin{keyword}
%% keywords here, in the form: keyword \sep keyword
\sep technoeconomic assessment \sep ResStock \sep thermal energy storage \sep salt hydrates \sep Ragone plots \sep air source heat pump \sep residential space heating \sep building decarbonization 
%% PACS codes here, in the form: \PACS code \sep code
%\PACS 0000 \sep 1111
%% MSC codes here, in the form: \MSC code \sep code
%% or \MSC[2008] code \sep code (2000 is the default)
%\MSC 0000 \sep 1111
\end{keyword}
\end{frontmatter}

\newpage
\linespread{1.5} \selectfont
\section{Introduction}
\label{sec:intro}
Decarbonizing the residential sector hinges on decarbonizing residential space heating. In the United States, over 21\% of primary energy consumption in 2021 was in the residential sector \citep{IEA_2022b}, 32\% of which was consumed for space heating \citep{EIA_2022}. A key decarbonization strategy for residential heating is replacing fossil-based heating with air source heat pumps (ASHPs) \citep{Lun_and_Tung_2020, Langevin_2019, Fischer_2014}. But widespread residential ASHP adoption faces several challenges, including poor efficiencies (or coefficients of performance (COPs)) in extreme cold despite cold-weather design improvements \citep{Wu_and_Skye_2021}, high capital costs, and long payback periods \citep{IEA_2022c}. Widespread ASHP adoption might also increase peak electricity demand, especially in winter peaking power systems \citep{Vaishnav_et_al_2020, Love_et_al_2017}, which must be accounted for in power system planning and operation.

Various energy storage technologies, including lithium ion batteries, flow batteries, and thermal storage (including ice storage), can be coupled with ASHPs in residential heating. Here, we perform the first broad assessment of a technology that might mitigate each of these challenges: salt hydrate thermal energy storage (TES). On-site short-term thermal energy storage (TES) is particularly promising for coupling with ASHP \citep{Arteconi_and_Polonara_2013, Fischer_2014, Moreno_2014, Odukomaiya_2021, Renaldi_2017, Alvaro_and_Cabeza_2017, Li_2018, Osterman_and_Stritih_2021, Yu_2021}. \citet{Odukomaiya_2021} shows that in a wide range of scenarios, TES can be more cost-effective than battery storage when coupled with ASHP. TES can use diverse materials as heating or cooling sources. Water tank heat storage, for instance, is widely used for building heating \citep{Zhao_and_Wang_2019}. But thermochemical heat TES using salt hydrates are a promising emerging technology due to their much higher energy densities (two to ten times higher than sensible heat capacity of water \citep{Cammarata2018, Donkers_2017, Michel2012}), lower losses than other residential TES materials \citep{Donkers_2017,Kenisarin_and_Mahkamov_2016, Wang_2019b, Kumar_2020, Wang_2020, Pardinas_et_al_2017}, reasonable stability, and low toxicity. In salt hydrate TES, heat is generated during salt hydration (i.e., when it uptakes water) and consumed during salt dehydration (i.e., when it releases water). During hydration, water (vapor) is impregnated into the salt lattice structure \citep{Kinzer_2023}. Despite their promise, salt hydrate TES has yet to be commercialized; rather, most existing work examines prototype salt hydrate TES devices \citep{NTsoukpoe_Kuznik_2021}. Furthermore, many potential salt hydrates might not be cost effective due to high material and/or device costs \citep{Kinzer_2023, NTsoukpoe_Kuznik_2021}. 

The cost effectiveness of TES in a coupled TES-ASHP system depends on multiple factors, including TES materials and designs, ASHP designs, the local climate, and buildings' heating load profiles. In the TES-heat pump system literature, techno-economic studies quantify the cost effectiveness of different TES-heat pump systems in realistically scaled building heating and cooling applications in diverse global geographical locations \citep{Wang_2019, Fischer_2014, Teamah_and_Lightstone_2019, Farah_2019, Wu_et_al_2020, Le_et_al_2019, Masip_2019, Zhu_et_al_2015, Arteconi_and_Polonara_2013, DEttorre_et_al_2019}. However, these analyses are limited to small samples of buildings (typically one or two buildings in a single city or location) and narrow time windows in several cases (e.g., one or two representative day(s) in a season). Conclusions from these studies are mixed, as some find significant potential in cost saving and peak load reduction from coupling TES with heat pumps \citep{Zhu_et_al_2015, Arteconi_and_Polonara_2013, DEttorre_et_al_2019}, while another finds TES is not cost-effective \citep{Le_et_al_2019}. More studies that can robustly quantify the value of TES when coupled with ASHPs are needed. Most importantly, no existing analyses consider high energy density salt hydrate TES, resulting in no techno-economic guidance for salt hydrate TES development in residential applications.

In this study, we conduct a techno-economic assessment of salt hydrate TES in residential space heating application. To model realistic operation of salt hydrate TES, we use Ragone plots, which indicate the relationship between energy and power density given material properties, and expected influences of device designs and operating conditions \citep{Kinzer_2023}. Ragone plots have been used in a wide range of applications, most notably for batteries to guide storage device designs \citep{Christen_Carlen_2000}. While the application of Ragone plots for TES devices is an emerging research topic, \citep{Woods_et_al_2021, James2022, Yazawa2019}, it has been used in techno-economic analyses of stand-alone battery storage \citep{Teichert_et_al_2023, Vykhodtsev_2022}. Thus far technoeconomic analyses have not been specifically applied to salts-based TES and preclude the incorporation of Ragone plots. To overcome these knowledge gaps, this study applies the Ragone framework, first developed by Kinzer et al. \citep{Kinzer_2023} for select salt hydrates to project cost-savings for salt-based TES systems. 

We quantify the value of coupling salt hydrate TES with ASHP across 4,800 representative U.S. households spread across 12 U.S. cities. We specifically quantify the value of four salt hydrates - \ce{MgSO4}, \ce{MgCl2}, \ce{K2CO3}, and \ce{SrBr2} - with varying power and energy densities and of six TES designs with varying charging and discharging limits and storage sizes (salt masses). We quantify TES value in terms of residential heating electricity cost savings and peak electricity consumption reduction. In so doing, we make two key contributions. First, we quantify the value of diverse TES materials and designs when coupled with ASHPs across a broad set of U.S. households, providing unique guidance for salt hydrate TES material and device development. Next, we incorporate salt-hydrate-specific Ragone plots, developed first by \citep{Kinzer_2023}, into a techno-economic analysis framework, advancing modeling capabilities to assess cost effectiveness of TES in residential applications. The purpose of our analysis is to generate a first-order estimate of economic viability of salt hydrate TES for space heating by analyzing its value across a large number of U.S. households. To achieve this desired scale, our analytical approach explicitly models key elements of TES-ASHP design and operation, and implicitly models other such elements (eliding some system design and operational details) using simple parameterizations. Due to uncertainty in capital costs of the emerging TES technologies considered in this study, we choose to minimize total operating costs, then use estimated operational savings to estimate break-even capital costs. Our break-even capital costs reflect the maximum cost a household could spend on a TES device while still achieving net savings from a TES device, thereby providing cost benchmarks for TES developers to reach. 

We develop an optimization-based techno-economic model of meeting residential space heating demand with an ASHP coupled with a salt hydrate TES device (hereafter referred to as the residential space heating, or RSH, model). In the coupled TES-ASHP system, the TES and ASHP can serve heating demand, and the ASHP can charge the TES. For each studied residential home, the RSH model minimizes space heating costs or peak demand by deciding ASHP capacity and/or ASHP and TES operations. We embed TES material specific Ragone plots into the RSH model using a piece-wise linear approximation. Using the RSH model, we quantify two key values of salt hydrate TES: (1) cost savings from avoided electricity purchases for space heating and (2) reductions in peak electricity demand for space heating. We run the RSH model for 4,800 representative single detached homes in 12 climatically-diverse U.S. cities, or 400 representative homes in each U.S. city, on an hourly basis for 2018 weather. This home sample size represents variability in housing characteristics within the total building stock in each city \citep{Deetjen_2021}. For each studied residential home, we run the RSH model in two steps for computational tractability. In the first step, the model optimizes the home's ASHP's capacity based on the home's annual peak heating demand. In the second step, using the optimized ASHP capacity from the first step as a parameter input, the model optimizes hourly operations of the TES and ASHP, which determines electricity purchases. In the second step, to quantify the value of TES, we run the RSH model for each household twice, once with only an ASHP and once with a coupled ASHP-TES system. Key RSH model inputs include hourly heating demands, hourly ASHP COPs, residential electricity retail prices, and TES Ragone plots and efficiency. Given uncertainties in future TES designs, we quantify the sensitivity of our results to TES sizes and designs. 

\section{Results}
\label{results}

\subsection{ASHP and TES operations across TES materials in a cold climate city}
\label{tes_hp_operation}
We first present detailed results for Detroit, a cold climate city, for our four salt hydrates -\ce{MgSO4}, \ce{MgCl2}, \ce{K2CO3}, and \ce{SrBr2} - to build intuition for their operation, sizing, and economics, then scale our results to 12 U.S. cities. Across the four salts considered, there are clear patterns of TES operations (Figures \ref{fig:tes_discharge_opt}). In general, TES discharges when ASHP's COP declines and/or is at a local minimum, and charges when ASHP's COP rises and/or is at a local maximum. While this operational relationship applies to all studied TES materials, TES discharging and charging differ between salts due to their different power and energy densities and the relationships between their power and energy densities (Figures \ref{fig:tes_discharge2}). Salts that have relatively higher energy densities compared to their power densities such as \ce{MgCl2} can discharge and charge for larger magnitude and for longer periods of time and can engage in more cycles per day than relatively lower energy density salts such as \ce{SrBr2}.

\begin{figure}[H]
     \includegraphics[scale=0.42]{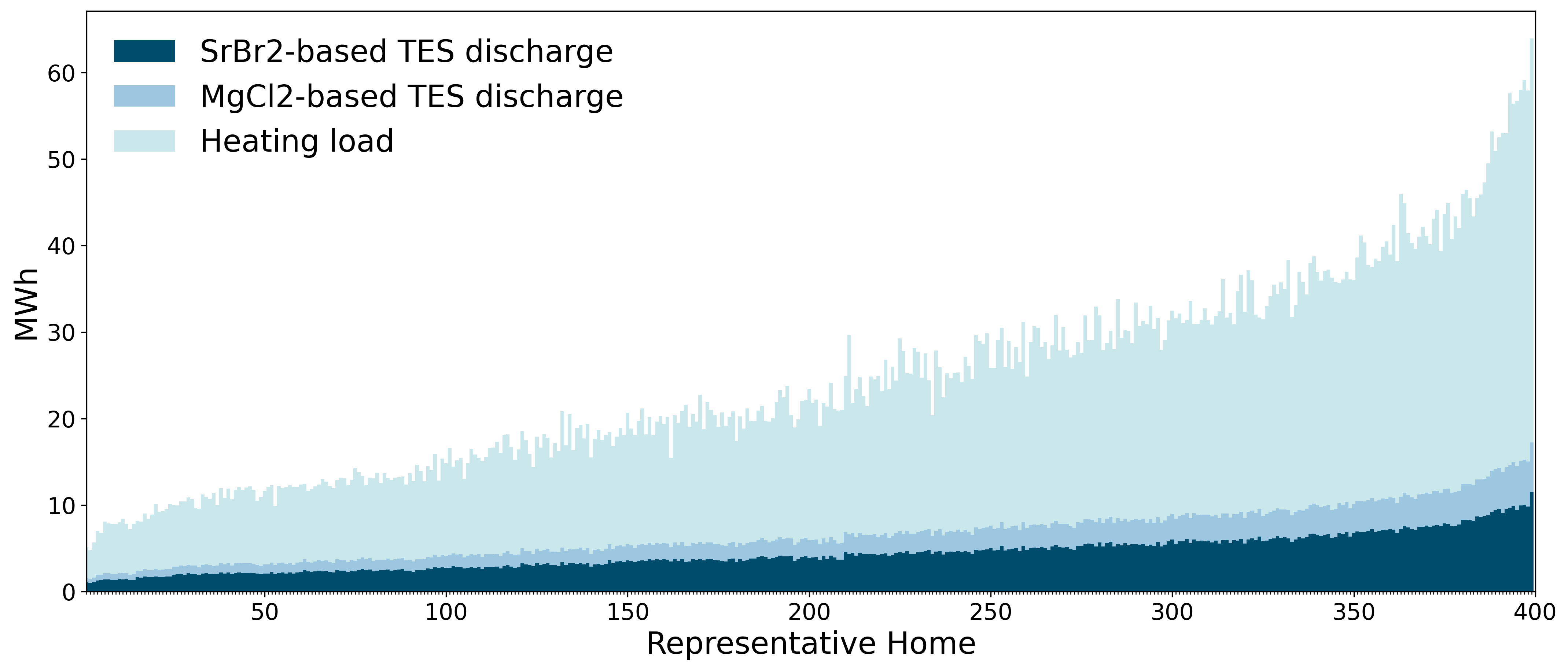} 
    \caption{\ce{MgCl2} and \ce{SrBr2}-based TES's energy outputs across 400 representative Detroit homes in 2018. Homes are sorted from lowest to highest annual peak load. The heights of the azure, light blue, and dark blue colors bars indicate heating loads, \ce{MgCl2}-based TES discharges, and \ce{SrBr2}-based TES discharges across homes, respectively.}
    \label{fig:tes_output_combined_opt}
\end{figure}

Figure \ref{fig:tes_output_combined_opt} compares annual TES discharge as a fraction of annual heating load across our 400 Detroit residential single family detached homes (hereafter "residential home" or just "home") when a high energy density salt (\ce{MgCl2}) and a lower energy density salt (\ce{SrBr2}) are used as TES materials. Across individual homes, loads shifted by TES differ widely. \ce{MgCl2} shifts between 1.2-17.2 MWh (26-27\%) of annual heating load, while \ce{SrBr2} shifts between to 0.8-11.5 MWh (17-18\%) of annual heating load. Aggregating across 400 representative homes, TES can shift up to 2.6 GWh (27\%) of total annual space heating load when coupled with ASHP. Specifically, \ce{MgSO4}, \ce{MgCl2}, \ce{K2CO3}, and \ce{SrBr2}-based TES, respectively shifts 0.15 GWh (1.5\%), 2.6 GWh (27\%), 1.8 GWh (19\%), and 1.7 GWh (18\%) of total annual space heating load (Figure \ref{fig:tes_output_opt}). Each of these 400 Detroit homes represents 664 actual homes (Table \ref{tab:scale_up_factor}), therefore, city-wide, \ce{MgSO4}, \ce{MgCl2}, \ce{K2CO3}, and \ce{SrBr2}-based TES can respectively shift 100 GWh, 1,727 GWh, 1,196 GWh, and 1,129 GWh of the city's total annual space heating load.

Salt hydrates that have relatively higher energy density (and thus higher energy capacity) (\ce{MgCl2}, see Table \ref{tab:tes_salts}) have the potential to shift more load within and across homes (Figures \ref{fig:tes_output1_opt} and \ref{fig:tes_output2_opt}) than salts with relatively higher power density but lower energy density (and thus limited energy capacity) (\ce{K2CO3} and \ce{SrBr2}) (Figures \ref{fig:tes_output3_opt} and \ref{fig:tes_output4_opt}). While \ce{MgSO4} has a relatively higher energy density, it also has a large parasitic load due to humidification requirements, reducing its load shifting utilization. 

\subsection{Economic value of TES for residential space heating in a cold climate city}
\label{economic_value}

To quantify the economic value of coupling TES with ASHP in meeting residential space heating demand, we calculate total annual costs of purchasing electricity for each home for two heating systems: (1) ASHP only and (2) a coupled TES-ASHP system. In Detroit, a cold climate city, the total annual cost of meeting space heating demands for 400 representative homes with only ASHP is \$602,000. Since each representative home represents 664 actual Detroit homes in ResStock (Table \ref{tab:scale_up_factor}) \cite{Deetjen_2021}, the total annual cost of meeting space heating demands with ASHP for all single detached family homes in Detroit is \$400 million. 

\begin{figure}[H]
     \includegraphics[scale=0.41]{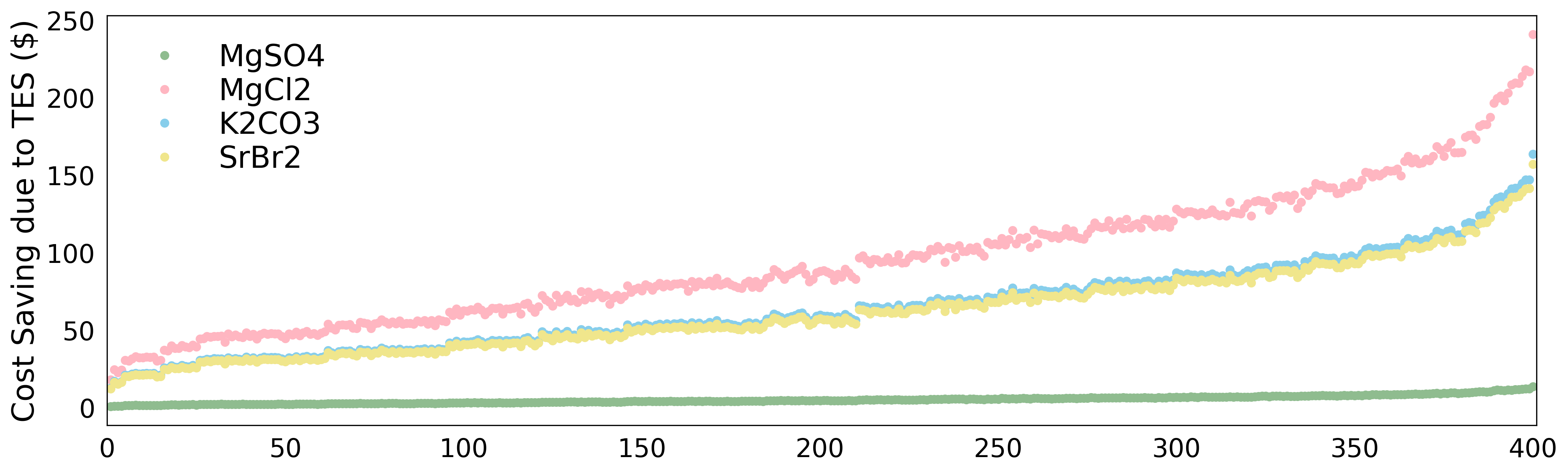}
     \caption{Annual cost savings in 400 representative Detroit homes from coupling TES of different salt hydrates with ASHP. Homes are sorted in order of lowest peak load to highest peak load.}
     \label{fig:cost_savings_opt_dollars}
\end{figure}

Relative to the only ASHP system, coupling TES with ASHP yields annual heating cost savings across 400 representative Detroit homes of \$2,012 (0.3\%), \$38,009 (6.3\%), \$25,764 (4.3\%), and  \$24,746 (4.1\%) for \ce{MgSO4}-, \ce{MgCl2}-, \ce{K2CO3}-, and \ce{SrBr2}-based TES, respectively. These are equivalent to annual cost savings across all actual Detroit homes of \$1.3 million, \$25.2 million, \$17.1 million, and \$16.4 million for  \ce{MgSO4}, \ce{MgCl2}, \ce{K2CO3}, and \ce{SrBr2}-based TES, respectively. Between homes, annual cost savings from TES differ widely, from \$0.7 to \$14 (0.2\%- 0.4\%), \$18 to \$241 (5.4\% - 7.6\%), \$13 to \$164 (3.6\% - 5.2\%), and \$12 to \$157 (3.5\% - 5.0\%) if \ce{MgSO4}, \ce{MgCl2}, \ce{K2CO3}, and \ce{SrBr2} are TES materials, respectively (Figure \ref{fig:cost_savings_opt_dollars}). These wide differences in cost saving across homes reflect the homes' wide range of heating loads, with higher heating load homes having higher TES discharging and cost savings. 

The cost-effectiveness of each TES salt depends the size of salt mass needed for TES operation. We model the size of each of the salt hydrates used as TES materials in kg. Here, we size TES using variable sizing, in which the TES salt mass is sized based on annual peak load of each building. Details on how we size TES are in \ref{appendix:tes_sizing} . We also consider incremental sizing and fixed sizing methods, which we include in sensitivity analysis.

To compare the cost effectiveness among TES salts, we also quantify the annual cost savings from TES per kg of salt. \ce{SrBr2} is the most cost-effective salt, saving \$1.8/kg to \$2.4/kg per year across our 400 Detroit homes, or \$2.2/kg on average per year. \ce{SrBr2} is followed by \ce{K2CO3}, which saves between \$1/kg and \$1.3/kg (\$1.2/kg on average) per year; then \ce{MgCl2}, which saves between \$0.7/kg and \$0.9/kg (\$0.8/kg on average); and finally \ce{MgSO4}, yielding cost saving between \$0.08/kg and \$0.15/kg (\$0.13/kg on average). Cost-effectiveness (or annual cost saving per kg) of each salt depends on the operation of each TES salt, which in turn depends on their energy density and parasitic load. \ce{MgSO4} has low cost-effectiveness due to its significant humidification parasitic load, while \ce{K2CO3} and \ce{MgCl2} have lower energy density (\ce{K2CO3}) or lower energy and power density (\ce{MgCl2}) than \ce{SrBr2}, limiting the amount and duration of their discharge.

\begin{figure}[H]
        \includegraphics[scale=0.41]{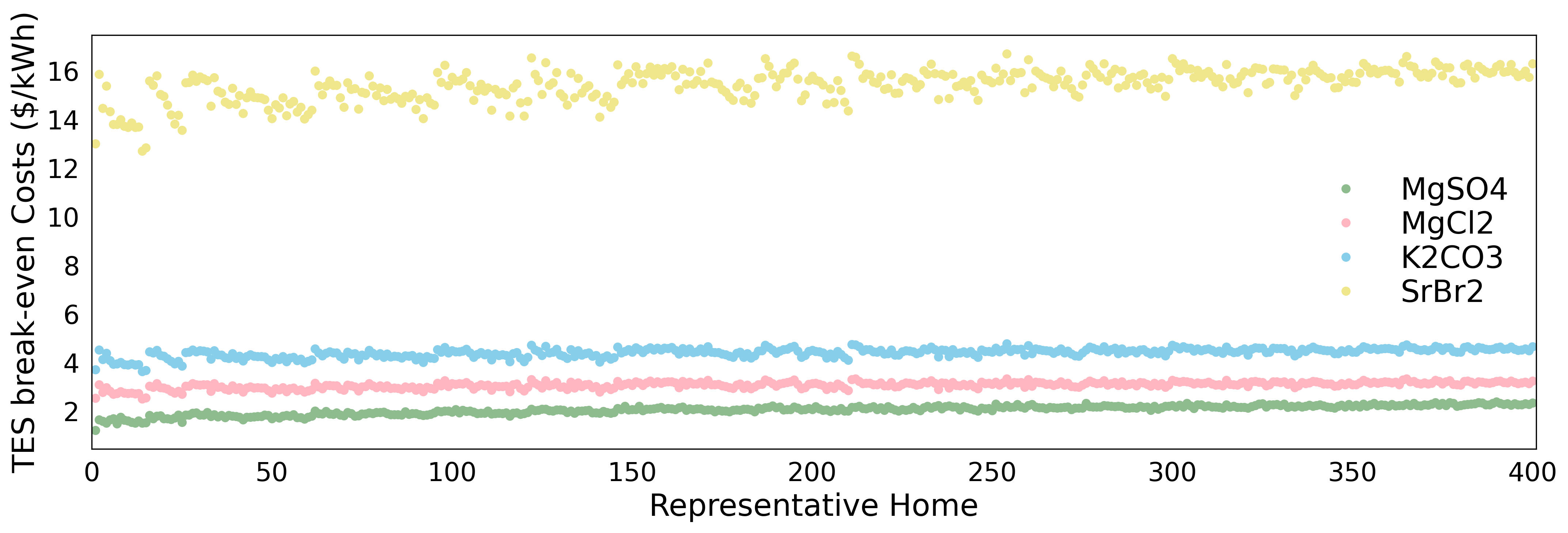}
        \caption{Break-even costs in 400 representative Detroit homes from coupling TES of different salt hydrates with ASHP. Homes are sorted in order of lowest peak load to highest peak load.}
        \label{fig:break_even_kWh}
\end{figure}

Future salt hydrate TES capital and installation costs are highly uncertain. Using the cost savings from TES calculated above, we calculate a break-even cost per MWh for each TES at each household. These break-even costs reflect the maximum cost a household could spend on a TES device while still achieving net savings from a TES device. Assuming a targeted TES device lifetime of 20 years \citep{DOE_2022}, the break-even costs of \ce{MgSO4}-, \ce{MgCl2}-, \ce{K2CO3}-, and \ce{SrBr2}-based TES across 400 representative Detroit homes are on average \$2.7, \$16, \$24, and \$43 per kg, respectively (Figure \ref{fig:break_even_kg_var}). We calculate these break-even costs as the average across households of total 20-year TES cost savings (or 20 times annual cost savings quantified above) divided by the TES salt mass. Due to different enthalpies among salts, these \$/kg break-even costs are, on average, equivalent to break-even costs of \$2.0/kWh, \$3.1/kWh, \$4.4/kWh, and \$15.4/kWh thermal for \ce{MgSO4}-, \ce{MgCl2}-, \ce{K2CO3}-, and \ce{SrBr2}-based TES, respectively (Figure \ref{fig:break_even_kWh}). Between homes, the break-even costs of \ce{MgSO4}-, \ce{MgCl2}-, \ce{K2CO3}-, and \ce{SrBr2}-based TES vary from \$1.2 to \$2.4, \$2.5 to \$3.3, \$3.6 to \$4.8, and \$12.7 to \$16.7 per kWh thermal, respectively. Differences in break-even costs between salts exceed differences in break-even costs between households.

\subsection{Load shifting and economic value of TES across U.S. cities}
\begin{figure}[H]
\centering
     \begin{subfigure}[c]{1\textwidth}
     \hspace*{-0.00001in}
         \includegraphics[scale=0.51]{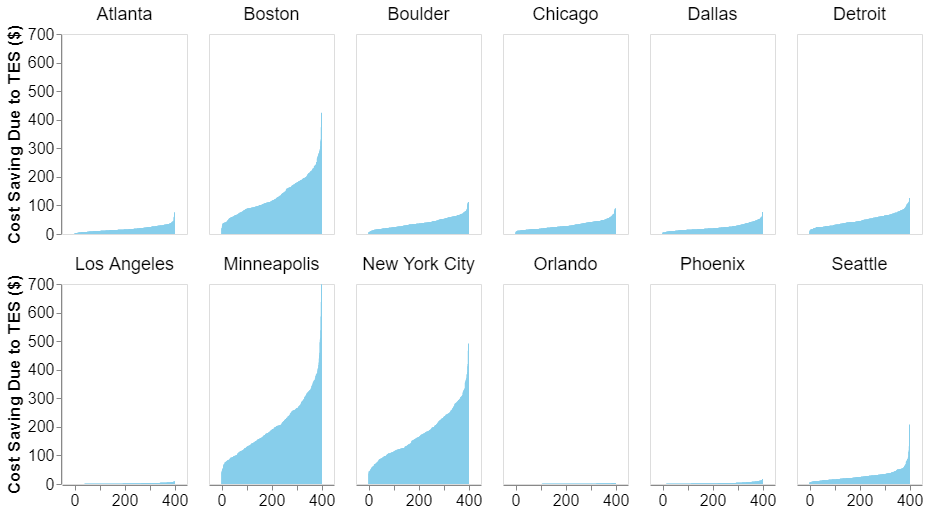}
         \caption{Cost saving due to TES in 400 representative homes in 12 U.S. cities (\$)}
         \label{fig:cost_savings_all1}
     \end{subfigure}    
     \begin{subfigure}[c]{1\textwidth}
     \centering
         \includegraphics[scale=0.51]{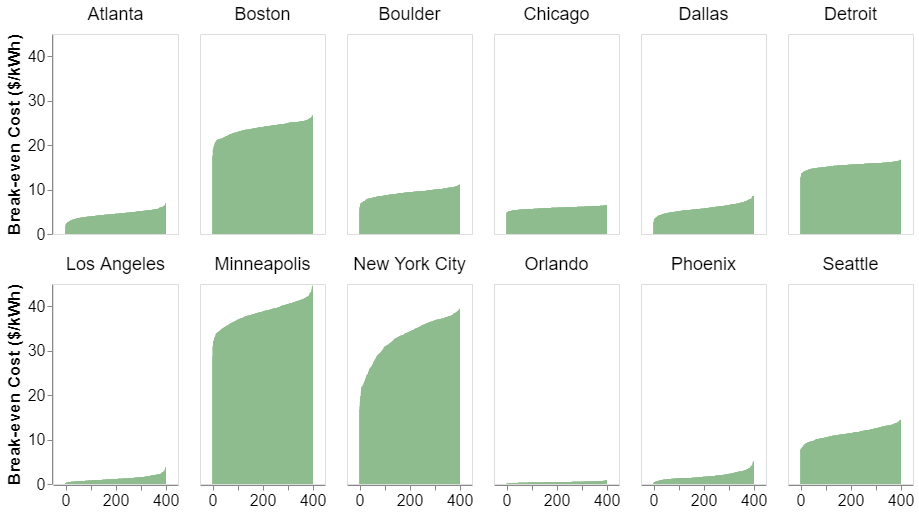}
         \caption{TES break-even costs across 400 representative homes in 12 U.S. cities (\$/kWh thermal)}
         \label{fig:cost_savings_all2}
     \end{subfigure}
        \caption{Cost savings due to \ce{SrBr2}-based TES in \$ (a) and break-even cost in \$/kWh thermal (b) in 400 representative homes in 12 U.S. cities. Cost savings and break-even costs in each city are sorted in order from lowest to highest. TES sizes here are determined using variable sizing method which depends on each home's peak load.}
        \label{fig:cost_savings_all}
\end{figure}

\begin{figure}[H]
\hspace*{-0.26in}
     \begin{subfigure}[c]{0.51\textwidth}
        \caption{Annual cost saving (million \$) due to TES}
         \includegraphics[scale=0.45]{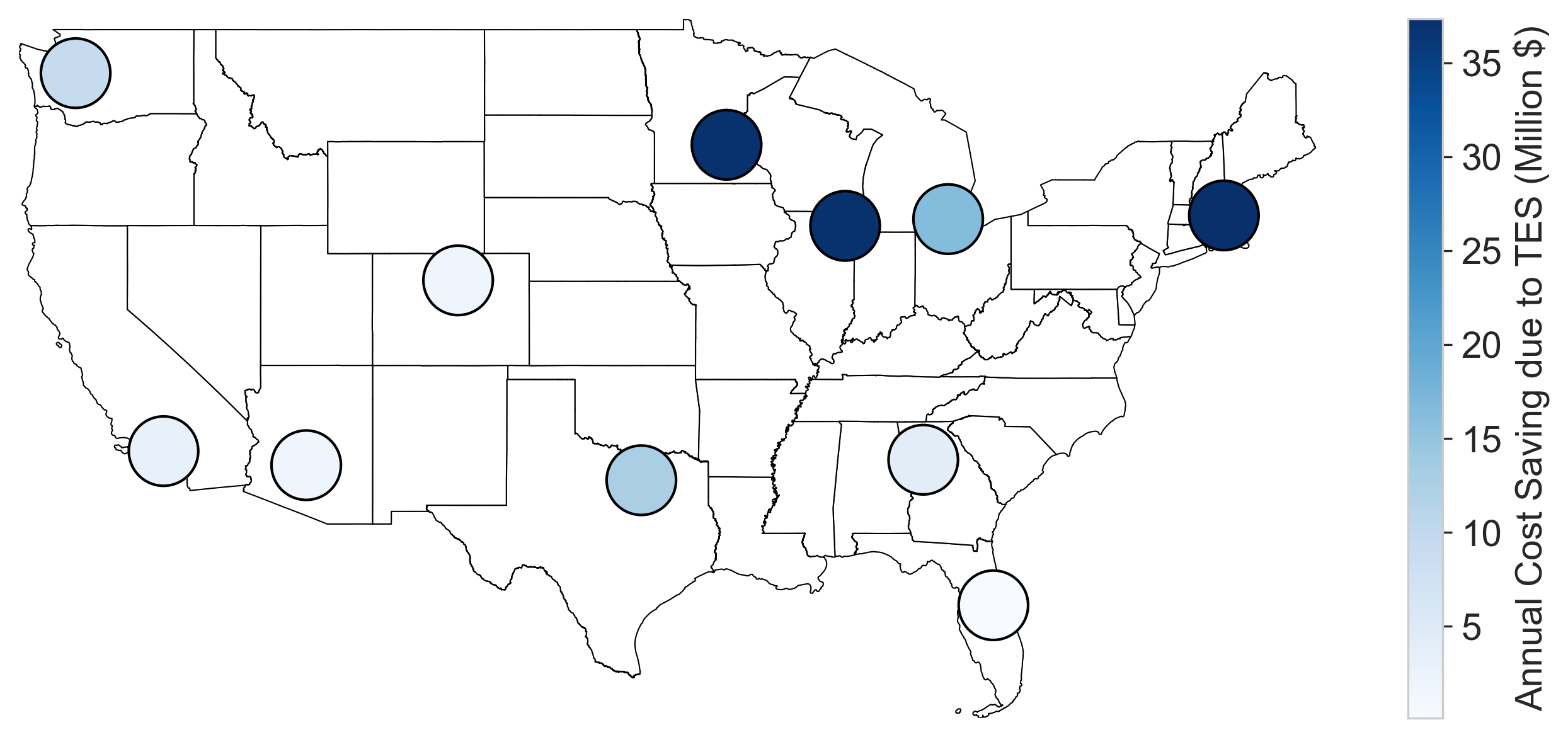}
         \label{fig:cost_savings_map1}
     \end{subfigure}
     \begin{subfigure}[c]{0.51\textwidth}
      \caption{TES break-even cost (\$/kWh thermal)}
         \includegraphics[scale=0.452]{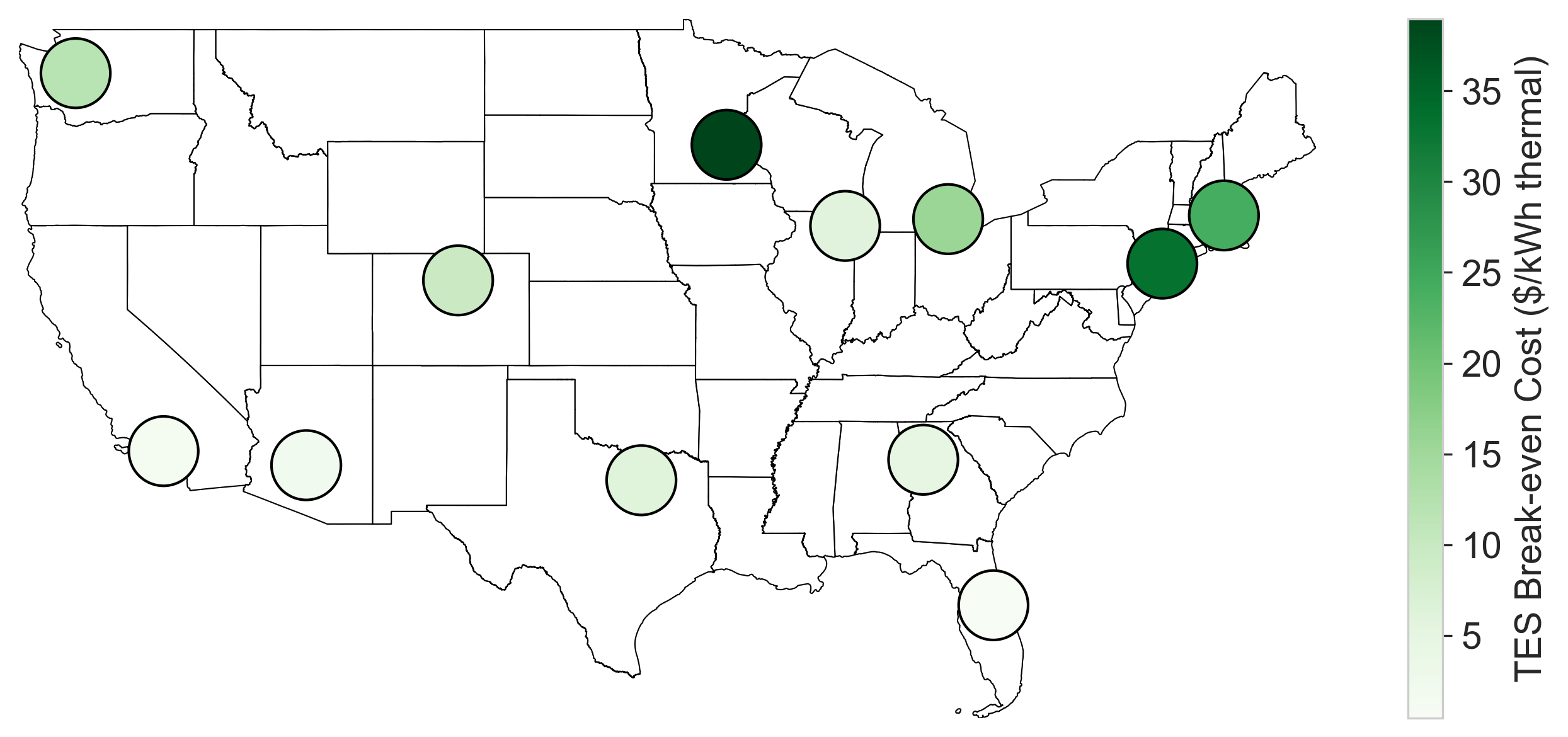}
         \label{fig:cost_savings_map2}
     \end{subfigure}
\hspace*{-0.23in}
     \begin{subfigure}[c]{0.51\textwidth}
     \caption{Peak Load Reduction (\%) due to TES}
         \includegraphics[scale=0.45]{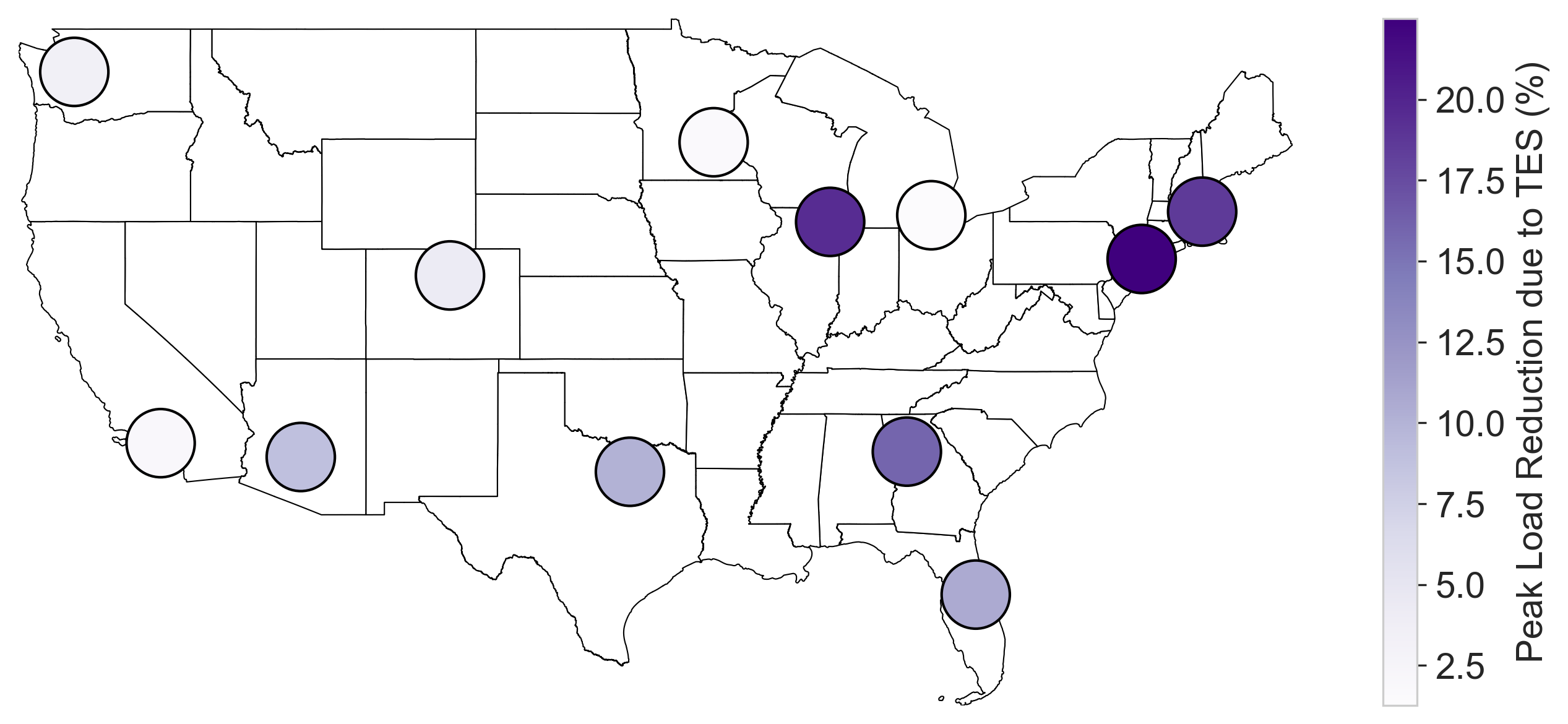}
         \label{fig:load_shifting_map1}
     \end{subfigure}
     \begin{subfigure}[c]{0.51\textwidth}
     \caption{Annual Load Shifting (\%) due to TES}
         \includegraphics[scale=0.45]{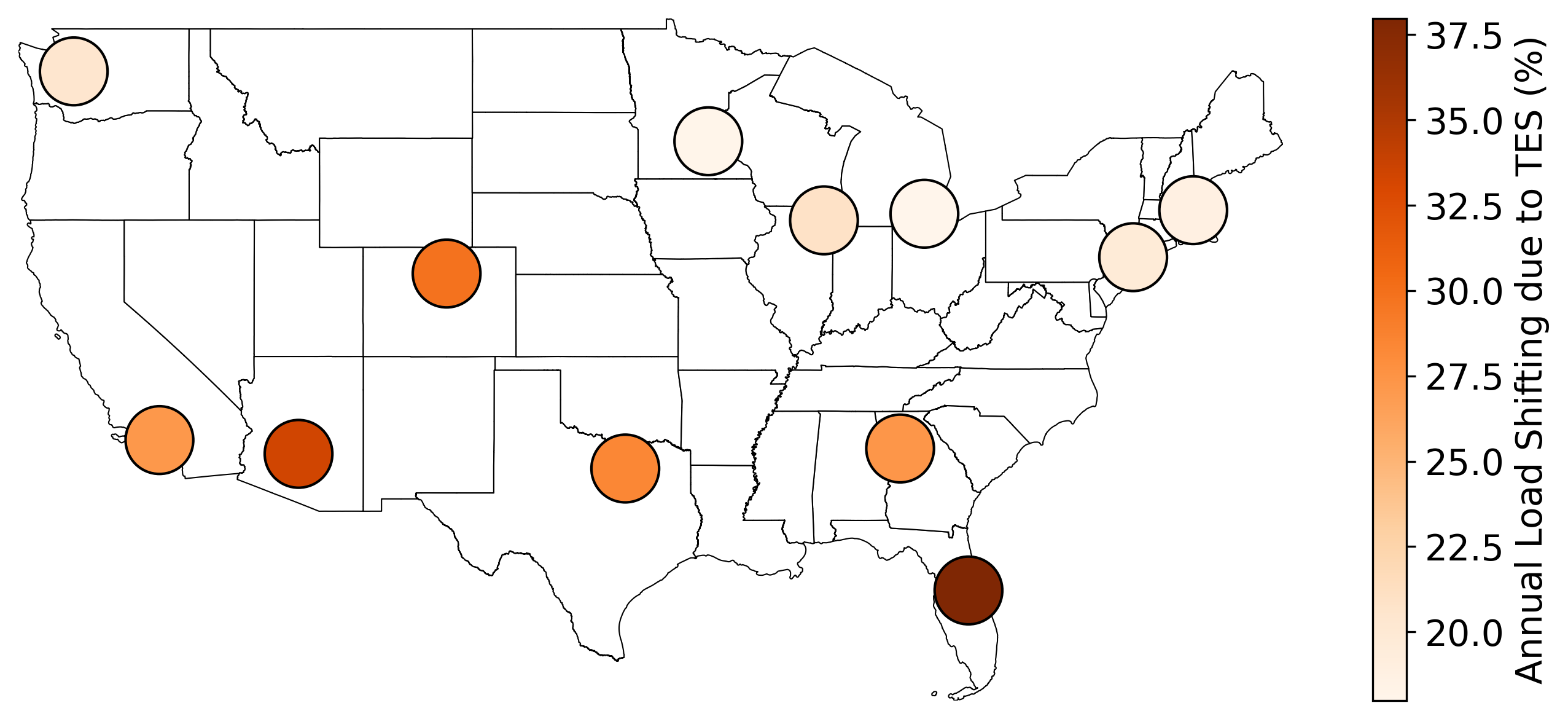}
         \label{fig:load_shifting_map2}
     \end{subfigure}
        \caption{Annual cost savings (\$) (a) and average break-even costs (\$/kWh thermal) (b) in 12 U.S. cities due to \ce{SrBr2}-based TES. Total annual amount of peak load reduction (\%) (c) and total annual amount of load being shifted (\%) (d) due to \ce{SrBr2}-based TES. New York City is omitted from panel (a) due to being an outlier. TES sizes here are determined using variable sizing method which depends on each building's peak load.}
        \label{fig:cost_savings_map}
\end{figure}

We run our RSH model for 400 representative homes and the most promising salt identified above (\ce{SrBr2}) in each of 12 cities that span three climates - cold (Boston, Boulder, Chicago, Detroit, Minneapolis, and New York), mild (Atlanta and Seattle), and hot (Dallas, Los Angeles, Orlando and Phoenix) (\citep{DOE_2017}). Summing across these representative homes by city, \ce{SrBr2}-based TES operations shifts 0.04 to 11.6 TWh (20\% to 38\%) of total annual heating loads (Figure \ref{fig:load_shifting_map2}), and incidentally reduces annual peak heating loads between 10.7\% to 22.5\% (16 MWh - 4.1 GWh) (Figure \ref{fig:load_shifting_map1}). Using TES to maximize peak load shifting rather than minimizing household costs could increase peak load reduction from 2\% to almost 12\% in a cold climate city (Section \ref{modified_model} and Figure \ref{fig:load_reduction}). Load shifting by TES saves up to \$699 (19\%) of total space heating costs annually in an individual home across cold climate cities, \$208 (12\%) annually in an individual home in mild climate cities, and \$77 (15\%) annually in an individual home in hot climate cities (Figure \ref{fig:cost_savings_all1}).

Cities in colder climates in general can save more with TES because they have higher heating load and more hours in which their ASHPs become inefficient. As a result, break-even costs of \ce{SrBr2}-based TES in cold climate cities are much higher than those in mild and hot climate cities (Figure \ref{fig:cost_savings_all2}). Break-even costs are highest in Minneapolis at \$29-\$45/kWh thermal across households, followed by in New York City at \$17-\$39/kWh thermal and Boston at \$18-\$27/kWh thermal. In mild climate cities, Seattle has the highest break-even costs for TES at \$8-\$14/kWh thermal across households. In hot climate cities, Dallas has the highest break-even costs for TES at \$3-\$9/kWh thermal across households.

\newgeometry{left=0.9in,right=0.8in,top=0.8in}
\restoregeometry

\newpage
\linespread{1.5} \selectfont
Scaling our results to city-wide level, Figure \ref{fig:cost_savings_map} illustrates annual cost savings, break-even costs, and annual peak load reduction and total load shifting due to \ce{SrBr2}-based TES across 12 cities. Annual city-wide cost savings reach up to \$602 million in New York City (omitted from Figure \ref{fig:cost_savings_map1} as it is an outlier) but other 11 cities have total cost savings between to \$0.11 million to \$37 million (Figure \ref{fig:cost_savings_map1}). Cost savings differ widely across cities. In a cold climate city except New York City, annual cost saving due to TES ranges from \$1.8 million (7\%) (Boulder) to \$37 million (10\%) (Boston). Average break-even costs for \ce{SrBr2}-based TES in a cold climate city can reach \$39/kWh thermal but they are at most \$12/kWh thermal in mild and \$6/kWh thermal in hot climate cities (Figure \ref{fig:cost_savings_map2}). These results emphasize the potential cost-effectiveness of TES in cold climate homes. 

\newgeometry{left=0.8in,right=0.5in}
\restoregeometry

\linespread{1.5} \selectfont

\subsection{Sensitivity analysis}
We test the robustness of our results for either all salt hydrates or just \ce{SrBr2} to alternative TES and ASHP design and operational schemes in a cold climate city (Detroit) (Table \ref{tab:sa_scenario}). One sensitivity modestly increases the value and cost-effectiveness of TES; all others modestly or significantly reduce its value and cost-effectiveness. For all four salt hydrates, increasing non-parasitic losses from 2\% to 5\% and from 2\% to 10\% would respectively reduce annual cost savings by between 3\% to 11\%, and between 13\% to 25\% across households and salt hydrates (Table \ref{tab:sa_results} and Section \ref{tes_materials}). Alternative TES operational schemes could require greater humidification needs, resulting in parasitic losses of 44\% to 70\% across salt hydrates (compared to our base case where only \ce{MgSO4} has humidification-driven losses of 39\%) (Table \ref{tab:sa_results} and Section \ref{tes_materials}). At these parasitic losses in Detroit, TES cost savings decrease 39\% to 93\% across households and salt hydrates, and break-even costs decrease by 9\% to 82\%. \ce{SrBr2} remains the salt hydrate with the greatest break-even costs for households, with break-even costs of above \$10 even under 44\% parasitic loss and additional non-parasitic losses. Finite rather than variable TES sizing, which might occur if TES is commercialized, would reduce TES break even capital costs by 9\% to 42\% due to over-sizing across salts despite increasing annual household cost savings by 23\% to 71\%. Using a gas furnace rather than electric baseboard as the ASHP backup would reduce \ce{SrBr2} cost savings by roughly 19\% to 24\% across homes and lead to up to 20\% lower break-even costs. Finally, running TES at a constant specific power output \citep{Zhang2020} (rather than variable output based on a Ragone plot) would increase break-even costs of \ce{SrBr2} by 3\%. For more details, see \ref{appendix:sa}. 

\section{Discussion}
\label{discussion}
To guide salt hydrate TES research and development, we provided the first broad techno-economic analysis of salt hydrate TES in residential space heating applications using Ragone plots for four TES materials \citep{Kinzer_2023} and heating profiles for 4,800 households across 12 U.S. cities. We found that salt hydrate TES can effectively reduce residential space heating costs when coupled with ASHP, particularly in cold climate cities. In Boston, Minneapolis and New York City, for instance, we found salt hydrate TES could save 10\%-44\% of household heating costs. Given these savings, we found households in Detroit could pay up to \$2.7, \$16, \$24, and \$43 per kg (or \$1.0, \$5.7, \$8.5, and \$15.3 per kWh thermal output) for \ce{MgSO4}, \ce{MgCl2}, \ce{K2CO3}, and \ce{SrBr2} salt, respectively, in a salt hydrate TES device and still break-even, i.e. save money. \ce{SrBr2} yields the highest cost saving per kg among the four salts due to its relatively high reaction enthalpy and power rating, and low parasitic humidification needs at typical residential indoor temperature, driving high energy output from a relatively small volume of salt. Our break-even costs provide upper bounds, or R\&D targets, for total TES device costs for economic viability in residential space heating. The other lower energy density salts' break-even costs fall far below that of \ce{SrBr2}, which means much higher cost reductions would be needed for them to be cost-effective in residential space heating applications. Our household-level results also provide an estimate of the size the market potential for TES in residential space heating. In cold climates, for instance, \ce{SrBr2}-based TES break-even costs exceed \$10/kWh in most homes in Boston, Boulder, Detroit, Minneapolis, and New York City, but in no homes in Chicago. In mild and hot climates, salt hydrate TES is less cost-effective due to lower heating needs and more efficient ASHP operations. In these climates, except for Seattle, \ce{SrBr2}-based TES could only reach break-even costs of less than \$6 /kWh thermal. 

Through a wide range of sensitivity analysis, we emphasized the importance of TES sizing methods and designs that developers of future TES devices should consider. First, strategically sizing TES based on household energy needs significantly increases the value of TES rather than using predetermined TES sizes, which can lead to under- or over-sizing of TES. Second, homes with electric baseboard backup, which is the most common ASHP backup system, benefit more from TES than homes with gas furnaces backup. Finally, consistent with \citep{NTsoukpoe_Kuznik_2021}, our results across all salts show that choosing TES materials with low humidification parasitic load is key in raising the cost-effectiveness of TES-ASHP systems. Government support for research and development should therefore target development of salt hydrates with low humidification parasitic loads.

In Detroit, \ce{SrBr2}-based TES can surpass the Department of Energy (DOE)'s TES subprogram's target of \$15/kWh thermal break-even capital costs in cold climates in most homes (Figure \ref{fig:break_even_kWh}. This DOE target can also be achieved with \ce{SrBr2}-based TES in all other cold climate cities except Chicago. TES also has significant potential to be an effective peak demand reduction tool, particularly when operated with the main objective of maximizing peak demand shifting (i.e., as a demand response tool). When operating TES to maximize peak demand shifting, \ce{SrBr2}-based TES reduced between 12\% and 27\% peak load across all homes, and reduced 395 MWh (11.6\%) total peak load city-wide without raising total system costs. Therefore, mass adoption of TES can help utilities and balancing authorities reduce energy generating resources capacity deployment and contribute to improving systems flexibility and resource adequacy.Energy policy could incentivize TES adoption by households, e.g. through tax credits that are commonly used to encourage adoption of other technologies. TES could be compensated for peak reduction either directly from a utility or from a third-party aggregator, which is a supplemental revenue stream we do not capture in our analysis. However, given the rarity at which demand response events are called, demand response payments would likely be a small share of the overall TES cost savings from daily cycling. 

Our work has several limitations that future work can address. First, to maintain computational tractability, we model heating demand and energy outputs from ASHP and TES energetically. Future research can combine detailed household designs with more detailed building integration modeling approaches, e.g. exergetic approaches \citep{DEttorre_et_al_2019, Arteconi_and_Polonara_2013, Farah_2019}, to obtain more finely tuned valuations for high-value TES designs identified in our study. To what extent an exergetic approach would change the value of TES is unclear. Our ASHP COP estimates, for instance, are in general higher than those estimated in ResStock; higher COP estimates reduce the potential value of TES, so our results might undervalue the cost-effectiveness of TES. Second, we focus on four salt hydrates as representative TES materials because they are widely studied for building heating and have experimentally reported operational data. Future studies could examine the many other promising salt hydrates for residential space heating using our analytical framework and datasets \citep{Clark_et_al_2022}. Finally, our study does not consider any TES material-specific challenges, for example, stability to thermal cycling. Additionally, the Ragone plots for salts are more representative of ideal material performance dictated by intrinsic kinetics. Future studies could apply our techno-economic framework to Ragone plots that capture more real-world TES performance. However, given that our studied salt hydrate TES face challenging economics in the form of low break-even capital costs (especially non-\ce{SrBr2} hydrates), better capturing real-world TES design and operations (e.g., through inclusion of piping costs or Ragone plots for specific reactor designs) would not qualitatively affect our insights.

Despite these limitations and opportunities for future work, our research yields several important insights for salt hydrate TES research and development. Our findings emphasize the potential market of commercially produced salt hydrate TES to provide energy input for residential space heating demand and the potentially large role salt hydrate TES could play in decarbonizing the residential sector. We highlighted the attributes of salt hydrates that would make them more economically effective (high energy density, moderately high power density, and low humidification parasitic load), and the importance of TES sizing in driving optimal TES break-even costs. Finally, we identified \ce{SrBr2} as the most promising salt of the four we studied to be used as material in TES devices that can potentially bridge and exceed DOE's TES capital cost target in many cold climate cities. 

\section{Experimental Procedures}
\label{methods}
Figure \ref{fig:model_overview} illustrates how we model TES integration within residential space heating systems. Our building-level model allows for ASHP to serve hourly heating load as well as charge TES using electricity purchased from utilities. When ASHP runs, its heating fan operates to circulate air. The salt hydrate TES system is operated with a fan and humidifier. When TES charges, it uses air from ASHP. When salt-hydrate TES discharges, a fan and humidifier run to maintain the home's humidification levels, and the TES uses recirculating air in the home. Our objective in this study is to provide a broad valuation of coupling salt hydrate TES systems with ASHPs for space heating within the U.S. residential sector. We capture physical constraints on TES operations by integrating salt-hydrate-specific Ragone plots (see \ref{appendix:tes_materials}) within a techno-economic model of a residential space heating system (see \ref{model_overview}). How a specific TES system would be integrated within each household's heating system would require detailed engineering analyses, which would require heating system schematics not available for a broad swath of U.S. households. We therefore do not conduct specific integration analyses for each household. Rather, in our techno-economic model, we capture key factors that would affect the integration and performance of an open TES system, specifically balance-of-system losses and energy required for hydration during discharge of the TES device (see Section \ref{tes_materials} for more details). 

\subsection{Model for space heating with ASHP and TES}
Our techno-economic model of residential space heating (``RSH model") is an optimization-based linear program that minimizes the total annual cost of purchasing electricity from a utility to meet residential space heating demand at the individual residential single family detached home (``residential home" or just ``home") level. Decision variables in the RSH model are hourly operation of TES, hourly operation of ASHP, and hourly electricity purchases to power the ASHP.  In the absence of TES, residential space heating is served by an ASHP with electric baseboard heating, electric furnace, or gas furnace as back-up systems at low temperatures. When TES is coupled with the ASHP/electric baseboard heater system, energy output from ASHP can be used to charge (or store heat in) TES and/or to serve heating demand in any given hour. TES can use stored heat at a later time to provide heat to serve demand. Due to the characteristics of TES materials used in this study, TES cannot charge and discharge at the same time. Therefore, during the hours when TES is discharging to serve heating load, ASHP can assist in serving heating load only, and during the hours when TES is charging, ASHP runs to serve heating load and charge TES. Constraints require residential heating loads to be met; enforce salt-hydrate-specific relationships between state of charge and maximum charge and discharge rates via piecewise linear approximations of the materials-specific Ragone plots; and vary the ASHP's COP with ambient air temperatures. By running this model for residential space heating when TES is and is not coupled with an ASHP, we quantify cost savings and ASHP operational changes due to TES, as well as TES operations. For model details, please see \ref{model_overview} in the SI.

To maintain computational tractability and a linear optimization, we model heating demand and energy outputs from ASHP and TES energetically (in units of kWh). The ASHP and TES devices used in this study are assumed to operate within the ranges of temperatures required for thermal comfort demands of our studied residential homes. Through our simplification of an energetic modeling approach, we can estimate the cost savings of diverse TES materials across a large number of households and climate zones, providing broad insights into the value of TES. Future research can adopt more household-specific modeling approaches, e.g. exergetic approaches \citep{DEttorre_et_al_2019, Arteconi_and_Polonara_2013, Farah_2019}, to obtain more finely tuned valuations for high-value TES designs identified in our study in smaller subsets of households.

To study the peak demand reduction, or demand response, potential of TES, we run a modified version of our RSH model that maximizes household peak load reduction instead of minimizing total cost (\ref{modified_model} in the SI). 

\subsection{TES materials and designs}
\label{tes_materials}
For thermal energy storage materials, we consider four types of salt hydrates: magnesium sulfate (\ce{MgSO4}), magnesium chloride (\ce{MgCl2}), potassium carbonate (\ce{K2CO3}), and strontium bromide (\ce{SrBr2}) (Table \ref{tab:tes_salts}). These salts are chosen as representative materials because they are widely studied for heating; have experimentally reported data for the reaction rates of hydration/dehydration of the salts; and have phase change dynamics vis-a-vis temperature and relative humidity amenable to residential space heating \citep{Kinzer_2023}. \ce{MgSO4} represents a salt with high reaction enthalpy but relatively low specific power \citep{Linnow2014}. This is contrasted by \ce{K2CO3}, which has low reaction enthalpy but high specific power \citep{Gaeini2019}. \ce{SrBr2} and \ce{MgCl2} offer characteristics between these two extremes \citep{Cammarata2018,Fisher2021,Hawwash2020}. Finally, all four salts have temperature inputs and outputs during charging and discharging aligned with ASHP outlet temperatures and residential space heating needs, respectively \citep{Kinzer_2023}. As such, we ignore temperature in our RSH model and use instead an energy-based model, allowing us to quantify value across a large number of residential households.  

\begin{table}[H]
\tablefontsize
\centering
\begin{tabular}{lccc} 
\hline 
\hline 
 Salt & Reaction Enthalpy (kWh/kg) & Specific Power at Peak Load (kW/kg) & Data Source\\
\hline
\ce{MgSO4} & 0.750 & 0.281 & \citep{Linnow2014}\\
\ce{MgCl2} & 0.193 & 0.085 & \citep{Hawwash2020}\\
\ce{K2CO3} & 0.186 & 1.646 & \citep{Gaeini2019}\\
\ce{SrBr2} & 0.356 & 0.811 & \citep{Cammarata2018}\\
\bottomrule 
\end{tabular} 
\caption{Power and Energy Requirements for the four salts used as TES materials in our study.}
\label{tab:tes_salts}
\end{table}

Reaction rates influence power densities for discharging (or heat release) and charging (or heat storage) during hydration and dehydration reactions, respectively. Thermogravimetric (TGA) measurements are typically operative under conditions where intrinsic, materials-specific kinetics is the main limiting factor \citep{Mahmoudi2021}, and have been previously applied to salt hydrate TES \citep{Kinzer_2023, Mahmoudi2021, Ye_et_al_2022}. Therefore, TGA measurements are agnostic to device design and represent the best-case power densities that can be accessed for a specific salt.While TGA discharge rates can be challenging to achieve in larger reactors, kinetic limitations become the primary constraining factor at high flow rates (around 0.25-0.5 m/s depending on salt and reactor design) \citep{Kinzer_2023, de_Jong_thesis}, with additional optimization required to minimize additional parasitic losses including pumping costs to overcome air pressure drop. For instance, Kinzer et al. and Aarts et al. (\citep{Kinzer_2023, Aarts2022}) show the benefits of using larger millimeter sized and porous salt hydrate pellets to minimize pressure drop in a packed bed reactor. We use data from TGA measurements to obtain best-fit functions for the discharging and charging steps based on operating reaction temperatures and humidity \citep{Kinzer_2023}. In general, as the SOC approaches 0, reaction rates non-linearly decrease as there is less salt left to be hydrated \citep{Gaeini2019}. During charging, the reactants and products flip, so reaction rates non-linearly decrease as the SOC approaches 1 (Figure \ref{fig:Ragone}). By controlling air temperature and relative humidity it is possible to achieve similar reaction times for both recharging and discharging \citep{Gaeini2019}. This allows us to simplify our model by assuming operating conditions that lead to a recharge curve that is a mirror image of the discharge curve with $sP_{discharging}(SOC=X) = sP_{charging}(SOC = 1-X)$. Further explanation for the Ragone plots and specific traits of each salt are included in \ref{appendix:tes_materials} and Figure \ref{fig:Ragone}. 

For each salt hydrate, we assume integration would occur via an open TES system. When discharging, recirculating indoor air is humidified if needed and passed through the TES device. During discharging, we account for fan power consumption for recirculating air and account for parasitic losses due to humidification \citep{NTsoukpoe_Kuznik_2021}. We estimate fan power consumption by generalizing results from NREL's ResStock analysis tool (see \ref{appendix:hp_fan}), and incorporate this fan power consumption in our RSH model (Section \ref{appendix:heating_load}). TES humidification needs will vary between households based on the household-specific space heating design, which is out of our analytical scope. Instead, we estimate parasitic losses from humidification for several alternative operational schemes. In each scheme, we estimate parasitic losses based on the heat of vaporization and the reaction enthalpy and vapor pressure requirements of each salt. Given the use of recirculating indoor air in our base case, we assume air flowing into the humidifier prior to the TES device is at 22 degrees Celsius and 20\% relative humidity. In this base case, only \ce{MgSO4} requires humidification of incoming air during discharging, resulting in a parasitic load of 39\%. We test the sensitivity of our results to larger humidification needs that result in parasitic loads ranging from 44\% to 70\% (Section \ref{sensitivity_analysis}). Parasitic losses are incorporated in our model as a discharging efficiency penalty equal to one minus the parasitic load from humidification (parameter $F_d$ in equation \eqref{eq:soc_TES}). When charging, warm air from the ASHP is passed through the TES device. Power consumption during charging is accounted for in the ASHP operations (equation \eqref{eq:purchased_elec}). Finally, we account for other system losses through a round-trip efficiency penalty of 2\% \citep{DeBoer_et_al_2014}, but vary this parameter as well via sensitivity analysis. 

We size each salt hydrate TES in three ways, reflecting different ways the mass of the salts for TES is determined when commercialized in the future: 1) Variable TES sizing; 2) Incremental (or step-wise) TES sizing; and 3) Fixed TES sizing. Details of each sizing approach are in \ref{appendix:tes_sizing}. Our main scenarios (discussed below) assume variable TES sizing, in which specific TES salt mass is sized to meet each home's annual peak load (kWh thermal), while other sizing methods are discussed in sensitivity analysis.

\subsection{Residential data}
We run our RSH model for one year (2018, to be consistent with the residential data we use taken from \citep{Deetjen_2021}) for each of 400 representative residential homes in one cold climate U.S. city (Detroit, MI) as our main city of study, and for each of 400 representative single detached family homes in additional 11 U.S. cities across different climates to capture the values of TES across different climate areas (see \ref{scenarios}). This study design requires three key data inputs detailed in this section: hourly building-specific heating loads, hourly city-specific ASHP COPs, and hourly city-specific retail electricity rates. 

We simulate the space heating load of these residential homes using NREL's ResStock analysis tool (version 2.2.5) \citep{Wilson_2017}. The ResStock analysis tool is a database of household characteristics for millions of U.S. residential homes. Characteristics for each household are contained in an EnergyPlus file, allowing for bottom-up simulation of each household using EnergyPlus. From the ResStock database, we select 400 representative homes in each city to simulate in ResStock and then to run in our RSH model. A sample size of 400 homes captures most variability between single family detached homes in a given city, including 88-96\% of a city's annual heating demand \citep{Deetjen_2021}, allowing us to balance computational tractability with generalizability. Furthermore, ResStock was designed to run representative rather than all homes in a given region. To simulate hourly residential space heating load for each of our studied residential homes, we modify the building's heating source to electric baseboard heater (COP=1), then simulate the building's operations in EnergyPlus for 2018. By using electric baseboard heat as the only heating device in a home, we capture the actual heating demand of that home given a COP of one, i.e. hourly energy input to the heater equals its hourly energy output to serve heating demand. Further details of how we use ResStock to simulate residential space heating loads by household are described in \ref{appendix:heating_load}.

In the RSH model, we assume that each home already has an ASHP with a back up electric baseboard heater to serve heating demand. The electric baseboard heater provides heat at a COP of 1.0 during extreme cold (lower than -17\degree C) when the ASHP's COP falls below 1.0 (see \ref{appendix:hp_COP}). 

We estimate hourly ASHP COPs for each of our households using hourly city-specific weather data from 2018 and temperature-COP curves for a median cold climate ASHP \citep{Waite_and_Modi_2022} (Figure \ref{fig:COP}). Hourly outdoor air temperature data for all homes in each city are obtained from Automated Surface/Weather Observing Systems (ASOS/AWOS) at each city's major airport location \citep{IEM}, which is the same weather data source as that used in ResStock. 

To estimate cost savings from TES, we assume electricity to power the ASHP in our RSH model is purchased from a local distribution utility at a time-of-use retail rate that is the default/most popular rate of each state's major utility. We use time-of-use retail rate in our analysis to capture the ongoing trend of many utilities (such as DTE Energy, San Diego Gas and Electric, and Sacramento Municipal Utility District) increasingly implementing time-of-use rate as default retail rate for their customers, either as an alternative option to fixed retail rate or mandatory pricing rule to replace fixed retail rate. We obtain time-of-use rates from the most recent electric rate books of the utilities that are in charge of providing electricity to the buildings in the cities of our study. For Detroit, MI, which is within the service territory of DTE Energy, the fixed retail rate is 18.3\textcent/kWh and the time-of-use rate varies depending on the month and hour. Other cities' retail rates are shown in Table \ref{tab:retail_rates}.

\subsection{Scenarios}
\label{scenarios}
To quantify the value of different TES designs coupled with ASHPs across climate zones, retail electricity rates, and households, we run five main scenarios for all representative homes in each city (Table \ref{tab:main_scenario}). Scenario 1 meets space heating demand with only ASHP, which is powered by purchased electricity from grid at retail rate. Scenario 2, which is the Reference scenario, couples \ce{SrBr2}-based TES using variable sizing with ASHP. Scenarios 3 through 5 facilitate similar comparisons to Scenario 2 but using different salt hydrates as TES materials. Comparing these scenarios to Scenario 1 captures the economic value of each TES design.

To quantify the value of TES in different climate areas, we run the above scenarios in 12 U.S. cities: Atlanta, GA; Boston, MA; Detroit, MI; Chicago, IL; Dallas, TX; Denver, CO; Los Angeles, CA; Minneapolis, MN; New York, NY; Orlando, FL; Phoenix, AZ; and Seattle, WA. 

\begin{table}[H]
\tablefontsize
\centering
\begin{tabular}{lccc} 
\hline 
\hline 
 Scenario & Coupled TES & TES Material & TES Sizing Method  \\
\hline
\multicolumn{4}{l}{\textit{No TES (Baseline):}}\\
Scenario 1 & No & -- & -- \\
\hline \\ \hline
\multicolumn{4}{l}{\textit{TES Sizes Vary Depending on Buildings' Peak Loads:}}\\
Scenario 2 \textit{(Reference)} & Yes & \ce{SrBr2} & Variable \\
Scenario 3  & Yes & \ce{MgCl2} & Variable \\
Scenario 4  & Yes & \ce{K2CO3} & Variable \\
Scenario 5  & Yes & \ce{MgSO4} & Variable \\
\bottomrule 
\end{tabular} 
\caption{Main scenarios for each of the 400 representative homes in each city of our study}
\label{tab:main_scenario}
\end{table}

\subsection{Sensitivity analysis}
\label{sensitivity_analysis}
We test the robustness of our results to alternative TES and ASHP configurations and operations (Table \ref{tab:sa_scenario}). Scenarios 6 through 13 facilitate similar comparisons to Scenario 1 discussed above in \ref{scenarios} using different salt hydrates but with incremental and fixed TES sizing. Scenario 16 assumes constant \ce{SrBr2}-based TES specific power output of 100 W/kg, which is near the upper end of the range shown to be achievable with similar experimental prototypes\citep{Yan2021,Zhang2020,Michel2012}. Scenario 15 assumes the TES-ASHP system is backed up by a gas furnace instead of an electric baseboard. Scenarios 16 through 19 consider alternative operational schemes with consequences on the parasitic load from humidification of air during TES discharging. Previous scenarios use humidification to fully provide the required humidity to the TES device when discharging, thereby not affecting inside humidity during TES discharging. Scenarios 16 and 17 increase the humidity of outdoor air assuming 0 degrees Celsius air at 10\% RH to that required for TES discharging. In these scenarios, parasitic loads for our salts range from 44-70\%. In Scenarios 18 and 19, we assume non-parasitic TES losses are 5\% and 10\%. All of these scenarios are run for Detroit, MI.

\begin{table}[H]
\tablefontsize
\centering
\begin{tabular}{lccccc} 
\hline 
\hline 
 Scenario & Coupled & TES & TES Sizing  & City & TES Design \\
  & TES & Material & Method &  & \\
\hline
\multicolumn{6}{l}{\textit{Examining value of TES with incremental sizes:}}\\
Scenario 6 & Yes & \ce{MgSO4} & Incremental  & Detroit & Reference \\
Scenario 7  & Yes & \ce{MgCl2} & Incremental & Detroit & Reference\\
Scenario 8 & Yes & \ce{K2CO3} & Incremental  & Detroit & Reference\\
Scenario 9  & Yes & \ce{SrBr2} & Incremental   & Detroit & Reference\\
\hline \\ \hline
\multicolumn{6}{l}{\textit{Examining value of TES with fixed sizes:}}\\
Scenario 10 & Yes & \ce{MgSO4} & Fixed  & Detroit & Reference \\
Scenario 11  & Yes & \ce{MgCl2} & Fixed & Detroit & Reference\\
Scenario 12 & Yes & \ce{K2CO3} & Fixed & Detroit & Reference\\
Scenario 13  & Yes & \ce{SrBr2} & Fixed & Detroit & Reference\\
\hline \\ \hline
\multicolumn{6}{l}{\textit{Examining value of TES with different designs:}}\\
Scenario 14 & Yes & \ce{SrBr2} & Variable & Detroit & Constant rating (100W/kg)\\
Scenario 15 & Yes & \ce{SrBr2} & Variable & Detroit & Gas furnace backup \\
\hline \\ \hline
\multicolumn{6}{l}{\textit{Examining the consequences of high parasitic loads and other non-parasitic losses:}}\\
Scenarios 16-17 & Yes & All four salts & Variable & Detroit & High Parasitic Load$^*$ \\
Scenarios 18-19 & Yes & \ce{SrBr2} & Variable & Detroit & Low round-trip efficiency$^*$\\
\bottomrule
\multicolumn{6}{l}{\textit{\footnotesize{$^*$ See \ref{appendix:sa_scenario_2} for details.}}}\\
\end{tabular} 
\caption{Selected additional scenarios and sensitivity analysis of the study for each of the 400 representative homes in Detroit, MI}
\label{tab:sa_scenario}
\end{table}
\section{Data and Code Availability}
Data and model code are available at \url{https://github.com/atpham88/TES}.

\section{Acknowledgements}
We thank Parth Vaishnav from the University of Michigan's School for Environment and Sustainability (SEAS) for invaluable housing stock data guidance and suggestions, Claire McKenna from SEAS for expert insights into heat pump settings and operations, and Pamela Wildstein from SEAS for assistance with utilities' electric rate books. We also thank Steven Kiyabu from the University of Michigan's Department of Mechanical Engineering for initial assistance in choosing promising salts for this study. Funding for this research was provided by the Graham Sustainability Institute’s Carbon Neutrality Acceleration Program (CNAP).

\section{Author Contributions}
\textbf{An T. Pham:} Conceptualization, techno-economic model development, model input preparation, results analysis \& interpretation, writing – original draft, writing – review \& editing, visualization. \textbf{Bryan Kinzer:} Conceptualization, performing analyses to generate the Ragone plots as model input, identifying and selecting the salts investigated in this study, formulating and developing the Ragone framework for salt hydrates, writing – review \& editing. \textbf{Ritvik Jain:} Performing ResStock simulation to obtain heating load profiles as model input. \textbf{Rohini Bala Chandran:} Conceptualization, identifying and selecting the salts investigated in this study, formulating and developing the Ragone framework for salt hydrates, results analysis \& interpretation, writing – review \& editing, supervision, securing project funding. \textbf{Michael T. Craig:} Conceptualization, model development, results analysis \& interpretation, writing – review \& editing, supervision, securing project funding. 

\section{Declaration of Interests}
The authors declare no competing interests.

\bibliographystyle{chicago}
\bibliography{tes_bib}

\appendix
\include{tes_model_appendix}

\end{document}

%% file: tes_model_appendix.tex
\renewcommand{\thetable}{A.\arabic{table}}
\setcounter{table}{0}
\renewcommand{\thefigure}{A.\arabic{figure}}
\renewcommand{\theequation}{\Alph{section}.\arabic{equation}}
\setcounter{figure}{0}
\setcounter{equation}{0}
\setcounter{page}{1}
\begin{center}
\textbf{Appendix}  \\
\vspace{0.1cm}
\textsc{Assessing the Value of Coupling Thermal Energy Storage with Heat Pumps for Residential Space Heating in U.S. Cities.} \\
\vspace{0.2cm}
\textit{For Online Publication.} \\

\end{center}

\section{The RSH Model to Minimize Heating Costs}
\label{model_overview}
The Residential Space Heating Model (RSH Model) is an optimization-based techno-economic model to quantify the cost saving from avoided purchased electricity from utilities for residential space heating when a thermal energy storage unit is coupled with existing heat pump compared to when residential space heating demand is supplied by solely heat pump. For each studied residential building, we run the RSH Model in two steps. In the first step, the model optimizes heat pump's capacity and purchased electricity to power said heat pump. In the second step, using the optimal heat pump capacity from the first step as parameter input, the model optimizes hourly operation of TES (whether TES is charging or discharging, TES's state of charge, energy discharge or charge from/to TES to shift load, and TES's power rating based on its state of charge), hourly operation of heat pump (whether heat pump is charging TES or not, heat pump's energy output to TES, and heat pump's energy output to serve load), and purchased electricity to power heat pump/resistance heater. The RSH Model runs on an hourly basis for a calendar year using several inputs -- hourly residential space heating load, hourly heat pump's imputed coefficient of performance (COP) from external temperature, residential electricity prices, and different TES characteristics (charging efficiency, power specifics-energy specifics trade offs). With model outputs, we quantify operations of heat pump and TES, and variable costs to meet residential space heating demand with and without TES.

\begin{figure}[H]
    \centering
    \includegraphics[scale=0.5]{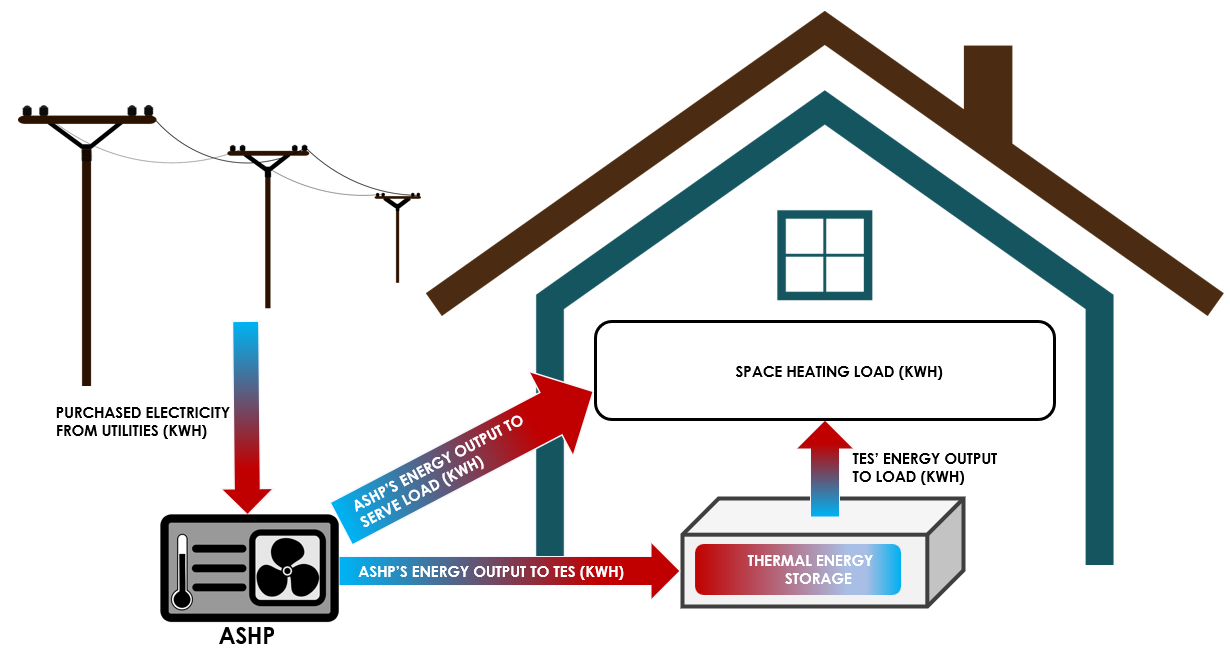}
    \caption{Overview of the RSH Model in Each Residential Home.}
    \label{fig:model_overview}
\end{figure}

The RSH model minimizes total annual costs of residential space heating for each household, where heating costs equal the sum of hourly purchased electricity times the hourly electricity rate:
\begin{align}
\label{eq:obj}
    TC_j = \sum_t P^R_{jt} d^{HP}_{jt}
\end{align}
where $j$ and $t$ index home and time (hour), respectively; $P^R_{jt}$ is the hourly retail electricity rate (\$/kWh); and $d^{HP}_{jt}$ is the hourly electricity purchased from utilities (kWh) to power ASHP when its COP $\geq$ 1 or backup electric baseboard when ASHP's COP $<$ 1. 

For each home, the model ensures hourly energy supply meets residential space heating demand:
\begin{align}
\label{eq:mc_electricity}
    g^{HP-L}_{jt} + g^{TES}_{jt} & \geq D_{jt}, \, \forall j \in \mathbb{J}, \, \forall t \in \mathbb{T}
\end{align}
where $g^{HP-L}_{jt} \geq 0$ is the hourly energy output from ASHP to serve heating demand (kWh); $g^{TES}_{jt} \geq 0 $ is the hourly energy output (or discharge) from TES to serve heating demand (kWh); and $D_{jt}$ is the hourly residential space heating demand (kWh).

Energy outputs from ASHP are constrained by ASHP's capacity:
\begin{align}
    g_{jt}^{HP-L} + g_{jt}^{HP-TES} & \leq K_j^{HP}, \, \forall j \in \mathbb{J}, \, \forall t \in \mathbb{T} \label{eq:hp_capmax}
\end{align}
where $g^{HP-TES}_{jt} \geq 0$ is hourly energy output from ASHP to charge TES in home $j$ at time $t$ (kWh thermal) and $K_j^{HP}$ is the maximum generating capacity of ASHP in home $j$ at any given time $t$ (kWh) calculated based on the home's peak heating load.

The efficiency, or coefficient of performance ($COP_{jt}$), of an ASHP varies with the ambient air temperature and becomes smaller as the air temperature decreases \citep{Waite_and_Modi_2022} (See \ref{appendix:hp_COP} for details). As COP declines for a constant heating requirement, more electricity will be consumed by the ASHP, as captured in \eqref{eq:purchased_elec}:
\begin{align}
    \label{eq:purchased_elec}
    COP_{jt} d^{HP}_{jt} \geq g^{HP-L}_{jt} + g^{HP-TES}_{jt}, \, \forall j \in \mathbb{J}, \, \forall t \in \mathbb{T}
\end{align}
where $COP_{jt}$ is calculated as a function of home $j$'s external temperature at time $t$ following a methodology from \citep{Waite_and_Modi_2022}, as discussed in \ref{appendix:hp_COP}. 

TES's hourly energy output, input (i.e., ASHP's output to charge TES, $g^{HP-TES}_{jt}$), hourly power rating, and state of charge are bounded by zero and maximum ratings, as captured in equations \eqref{eq:TES_power_rating}, \eqref{eq:TES_soc}, \eqref{eq:TES_power_rating_charge}, and \eqref{eq:soc_TES_bounds}:
\begin{subequations}
    \begin{align}
        g^{TES}_{jt} & \leq k^{TES,D}_{jt}, \, \forall j \in \mathbb{J}, \, \forall t \in \mathbb{T} \label{eq:TES_power_rating} \\
        g^{TES}_{jt} & \leq x^{TES}_{jt}, \, \forall j \in \mathbb{J}, \, \forall t \in \mathbb{T} \label{eq:TES_soc}\\
        g^{HP-TES}_{jt} & \leq k^{TES,C}_{jt}, \, \forall j \in \mathbb{J}, \, \forall t \in \mathbb{T}
        \label{eq:TES_power_rating_charge}\\
        0 \leq x^{TES}_{jt} & \leq E^{TES}_j, \, \forall j \in \mathbb{J}, \, \forall t \in \mathbb{T}, \label{eq:soc_TES_bounds}
    \end{align}
\end{subequations}
where $k^{TES,D}_{jt}$ and $k^{TES,C}_{jt}$ are TES' hourly discharging and charging (kWh) in home $j$ at time $t$, both of which are dependent on TES's state of charge which equals $\dfrac{x^{TES}_{jt}}{E^{TES}_j}$. The relationship between charging and discharging power ratings, $k^{TES,C}_{jt}$ or $k^{TES,D}_{jt}$, and TES's state of charge, $\dfrac{x^{TES}_{jt}}{E^{TES}_j}$ is captured in the Ragone plot of each of our salt-specific TES systems illustrated in Figure \ref{fig:Ragone}. We assume two-segment piece-wise linear approximations of these Ragone plots in the model, which allow us to increase computational tractability yet still maintain good approximations of salt based TES's state of charge given its charging/discharging capabilities at any given time. Details on formulation of these approximations of Ragone plots are in \ref{add_formulation}.

TES state of charge varies with charging and discharging decisions (\ref{eq:soc_TES}), and is bounded between zero and a maximum state of charge (\ref{eq:soc_TES_bounds}):
\begin{subequations}
    \begin{align}
        x^{TES}_{jt} & = x^{TES}_{j(t-1)} + F^c g^{HP-TES}_{jt} - F^d g^{TES}_{jt}, \, \forall j \in \mathbb{J}, \, \forall t \in \mathbb{T}, \label{eq:soc_TES} \\
        x^{TES}_{j(t=1)} & = 0, \, \forall j \in \mathbb{J}, \, \forall t=0 \in \mathbb{T} \label{eq:soc_TES_initial}
    \end{align}
\end{subequations}
where $x^{TES}_{j(t-1)}$ is TES's state of charge in home $j$ at previous hour $t-1$ (kWh), and $F^c$ and $F^d$ are TES charging and discharging efficiencies, respectively (see Section \ref{tes_materials} for efficiency values). The initial state of charge of TES $x^{TES}_{j(t=1)}$ is assumed to be 0 (\ref{eq:soc_TES_initial}). 

\section{Modified RSH Model to Maximize Peak Demand Reductions}
\label{modified_model}
To study the peak demand reduction, or demand response, potential of TES, we run a modified version of our RSH model that maximizes household peak load shifting instead of minimizing total cost. Our modified objective function, in this case, is:
\begin{align}
    LR_j = \max d^s_j
\end{align}
where $d^s_j$ is total peak load being shifted in representative home $j$.

Our modified model limits total system's operational cost to be at most the total operational cost from the non-modified RSH model (ref prior section) when all heating load is met by ASHP (i.e., without TES):
\begin{align}
    \sum_t P^R_{jt} d^{HP}_{jt} \leq \bar{TC}_j
\end{align}
where $\bar{TC}_j$ is the least cost solution from the RSH model when TES is not coupled with ASHP.

Finally, the energy outputs from ASHP must be less than the reduced peak load:
\begin{align}
    g^{HP-L}_{j(t=\text{peak hour})} + g^{HP-TES}_{j(t=\text{peak hour})} \leq D^{p}_{j(t=\text{peak hour})} - d^s_j
\end{align}

The rest of the modified model's constraints remain the same as in the original model.

We programmed the RSH model in Python's Pyomo and solved it using CPLEX Version 20.1.0.1 (see \ref{code}).

\section{Piece-wise Linear Functions of Ragone Plots}
\label{add_formulation}
To reduce computational time, we linearize the Ragone plots for the salts of our study upon both charging (heating) and discharging (cooling) into two-segment piece-wise linear functions. These functions are defined over domain of $[0,1]$ which is the state of charge of the salt hydrate specific TES devices. The two-segment piece-wise linear Ragone plot for discharging is defined as $f_d: [0,1] \to \mathbb{R} = m_d(x) + b$ being a continuous function such that $f_d$ over $[0,a]$ is a linear interpolation from $f_d(0)$ to $f_d(a)$ and $f_d$ over $[a,1]$ is a linear interpolation from $f_d(a)$ to $f_d(1)$. Similarly, the two-segment piece-wise linear Ragone plot for charging is defined as $f_c: [0,1] \to \mathbb{R} = -m_d(x) + b$ being a continuous function such that $f_c$ over $[0,a]$ is a linear interpolation from $f_c(0)$ to $f_c(a)$ and $f_c$ over $[a,1]$ is a linear interpolation from $f_c(a)$ to $f_c(1)$. $m_d$ and $b$ are fitted coefficients of the Ragone plots based on empirical data (see \ref{appendix:tes_materials}). The different formulas for Ragone plots for discharging and charging reflect the different power ratings at a given SOC when TES discharges or charges. Specifically, TES' maximum discharging and charging (kWh) in building $j$ at time $t$ are dependent on TES's SOC at that same hour, $x^{TES}_{jt}$ (kWh), and TES's energy capacity, $e^{TES}_j$ (kWh):
\begin{subequations}
    \begin{align}
        k^{TES,D}_{jt} = \dfrac{y_B - y_A}{x_B-x_A} \dfrac{x^{TES,D}_{jt}}{e^{TES}_{j}} + y_A - \dfrac{y_B - y_A}{x_B-x_A} x_A, \, \forall j \in \mathbb{J},  \, \forall t \in \mathbb{T} \label{appendix:TES_ragone_discharging} \\
        k^{TES,C}_{jt} = \dfrac{y_B - y_A}{x_B-x_A} \left(1-\dfrac{x^{TES,D}_{jt}}{e^{TES}_{j}}\right) + y_A - \dfrac{y_B - y_A}{x_B-x_A} x_A, \, \forall j \in \mathbb{J},  \, \forall t \in \mathbb{T} \label{appendix:TES_ragone_charging} \\
        \quad \text{for} \quad x_A \leq x^{TES}_{jt}  \leq x_B \nonumber
    \end{align}
\end{subequations}

where $(x_A, y_A)$ and $(x_B, y_B)$ are the coordinates of line segment $AB$ of piece-wise linear function of TES's state of charge and discharging and charging power rating. The different formulas for TES' discharging and charging power ratings reflect the different power ratings at a given SOC when TES discharges or charges.

Approximated coordinates of each segment $AB$ for each salt are shown in the figure below.
\begin{figure}[H]
    \centering
    \includegraphics[scale=0.35]{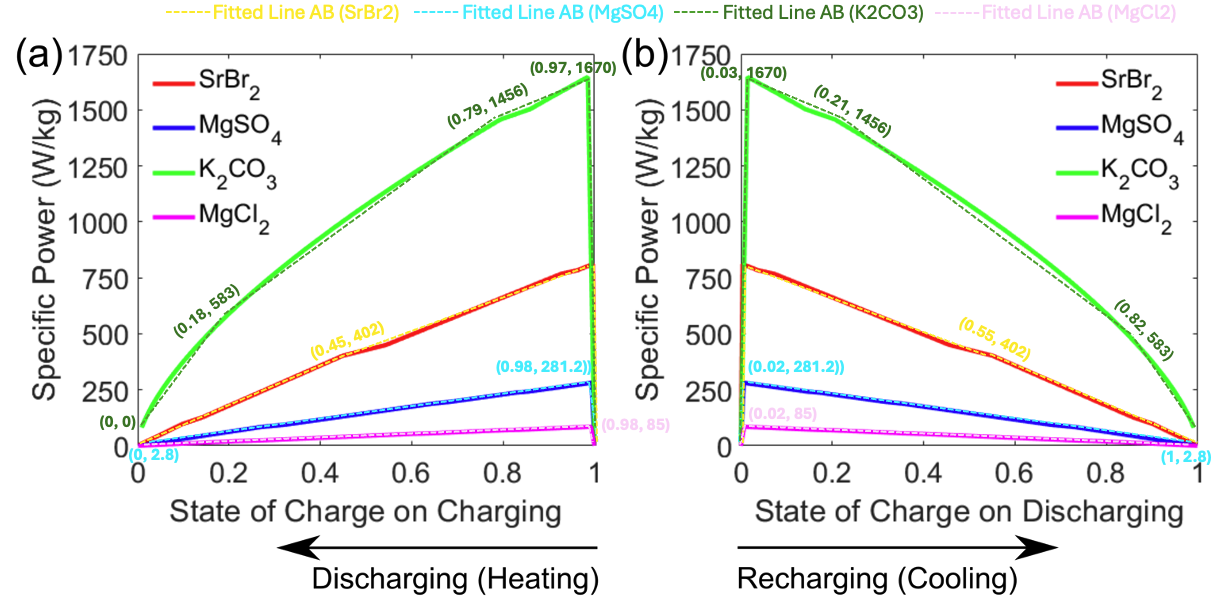}
\end{figure}

\section{Humidification Penalties}
\label{humid_penalties}
We estimate parasitic loads from humidification during TES discharging using Clasius Clapeyron phase diagrams for each of our salt hydrates. The phase diagrams indicate the relative humidity required for phase changes of each TES salt. In particular, hydration of salt hydrates (during discharging of heat from TES) requires humid air input. Humidification of air can pose a large parasitic load on the TES system \citep{NTsoukpoe_Kuznik_2021}. We assume parasitic load of dehydration (or charging) is negligible. To estimate parasitic load of hydration, we assume outlet temperatures from our TES of 80\degree F (or roughly 27\degree C), which is consistent with temperatures of hydration reactions from each of our salts and from an ASHP. Phase diagrams from each salt yield the relative humidity required at 27\degree C for hydration reaction \citep{Fisher2021, Linnow2014, Lele_and_Tamba_2017}. Given vapor pressures during hydration, we calculate the parasitic losses of humidification as:

\begin{align}
    \delta & = \dfrac{P_{target} - P_{outside}}{P_{target}} \\
    \beta & = \dfrac{\delta \left(\Delta H_{vap}\right)}{\Delta H_{rx,salt}} \\
    \text{Efficiency} & = 1 - \beta
\end{align}

where $\Delta H_{\text{vap}}$ vap is the heat of vaporization (2.26 MJ/kg \ce{H2O} ); $\Delta H_{\text{rx,salt}}$ is the reaction enthalpy of each salt hydrate; and $P_{\text{target}}$ and $P_{\text{outside}}$ are the vapor pressures required during TES discharge and of the ambient air passed through the TES device during discharging, respectively. We generate our results for several values of $\delta$ (i.e., for several levels of required humidification) to understand the effect of varying operational modes on the value of TES in residential space heating applications. TES humidification needs will vary between households based on the household-specific space heating design, which is out of our analytical scope. Our base case assumes that our open TES system uses recirculating indoor air at 20\% relative humidity and 22\degree C when discharging, so humidification is only necessary when indoor humidity is less than that required by the TES device during discharging. Under this base case, only \ce{MgSO4} requires energy for humidification, resulting in a parasitic load of 39\%. We also test the sensitivity of our results to fully providing the required humidity to the TES device when discharging, thereby not affecting the humidity levels inside the home during TES discharging, and to increasing the humidity of outdoor air assuming 0\degree C air at 10\% relative humidity. In these cases, parasitic loads for our salts range from 54 to 70\% and from 44-52\%, respectively (see Table below).

\begin{table}[H]
\tablefontsize
\centering
\begin{tabular}{lcccc}
\hline 
\hline
Salt & Efficiency &  Efficiency & Efficiency & Efficiency \\
& at 0\% RH & at 10\% RH & at 30\% RH & at 20\% RH \\
\hline
\ce{K2CO3} & 0.34 & 0.483 & 1 & 1\\
\hline
\ce{MgCl2} & 0.30 & 0.55 & 1 & 1\\
\hline
\ce{SrBr2} & 0.40 & 0.56 & 1 & 1\\
\hline
\ce{MgSO4} & 0.46 & 0.48 & 0.55 & 0.61\\
\bottomrule 
\end{tabular} 
\caption{Efficiencies of salt hydrates given humidification penalties for varying assumptions regarding relative humidity of input air during discharging. We test the sensitivity of our results to each penalty via sensitivity analysis. Our base case assumes that our open TES system uses recirculating indoor air at 20\% RH when discharging.}
\label{appendix:humidity_salts}
\end{table}

\section{Scenarios}
\subsection{Main Scenarios}

\begin{table}[H]
\tablefontsize
\centering
\begin{tabular}{lccc} 
\hline 
\hline 
 Scenario & Coupled TES & TES Material & TES Sizing Method   \\
\hline
\multicolumn{4}{l}{\textit{No TES (Baseline):}}\\
Scenario 1 & No & -- & -- \\
\hline \\ \hline
\multicolumn{4}{l}{\textit{TES Sizes Vary Depending on Buildings' Peak Loads:}}\\
Scenario 2 (Reference) & Yes & \ce{SrBr2} & Variable \\
Scenario 3  & Yes & \ce{MgCl2} & Variable \\
Scenario 4  & Yes & \ce{K2CO3} & Variable \\
Scenario 5  & Yes & \ce{MgSO4} & Variable \\
\bottomrule 
\end{tabular} 
\caption{Main scenarios of the study for each of the 400 representative homes in each city of our study}
\label{appendix:main_scenario}
\end{table}

\subsection{Sensitivity Analysis}

\begin{table}[H]
\tablefontsize
\centering
\begin{tabular}{lccccc} 
\hline 
\hline 
 Scenario & Coupled & TES & TES Sizing & City & TES Design \\
  & TES & Material & Method &  & \\
\hline
\multicolumn{6}{l}{\textit{Examining value of TES with incremental sizes:}}\\
Scenario 6 & Yes & \ce{MgSO4} & Incremental  & Detroit & Reference \\
Scenario 7  & Yes & \ce{MgCl2} & Incremental & Detroit & Reference\\
Scenario 8 & Yes & \ce{K2CO3} & Incremental  &  Detroit & Reference\\
Scenario 9  & Yes & \ce{SrBr2} & Incremental  & Detroit & Reference\\
\hline \\ \hline
\multicolumn{6}{l}{\textit{Examining value of TES with fixed sizes:}}\\
Scenario 10 & Yes & \ce{MgSO4} & Fixed  & Detroit & Reference \\
Scenario 11  & Yes & \ce{MgCl2} & Fixed & Detroit & Reference\\
Scenario 12 & Yes & \ce{K2CO3} & Fixed  & Detroit & Reference\\
Scenario 13  & Yes & \ce{SrBr2} & Fixed  & Detroit & Reference\\
Scenario 14 & Yes & \ce{SrBr2} & Variable & Detroit & Constant power rating (100W/kg)\\
Scenario 15 & Yes & \ce{SrBr2} & Variable & Detroit & Gas furnace backup \\
\bottomrule 
\end{tabular} 
\caption{Selected additional scenarios and sensitivity analysis of the study for each of the 400 representative homes in Detroit, MI}
\label{appendix:sa_scenario_1}
\end{table}

\begin{table}[H]
\tablefontsize
\centering
\begin{tabular}{lccccccc} 
\hline 
\hline 
 Scenario & City & \ce{MgSO4}-TES & \ce{MgCl2}-TES & \ce{K2CO3}-TES & \ce{SrBr2}-TES & Other \\
 &  & Parasitic Load & Parasitic Load &  Parasitic Load & Parasitic Load & TES Loss\\
\hline
\multicolumn{6}{l}{\textit{High Parasitic Load:}}\\
Scenario 16 & Detroit & 52\%  & 45\% & 52\% & 44\% & 2\%\\
Scenario 17 & Detroit & 54\%  & 70\% & 66\% & 70\% & 2\%\\
\hline \\ \hline
\multicolumn{6}{l}{\textit{High Non-parasitic TES Loss:}}\\
Scenario 18 & Detroit & 39\%  & 0\% & 0\% & 0\% & 5\%\\
Scenario 19  & Detroit & 39\%  & 0\% & 0\% & 0\% & 10\%\\
\bottomrule 
\end{tabular} 
\caption{Sensitivity analysis to test the consequences of high parasitic loads and other non-parasitic losses for each of the 400 representative homes in Detroit, MI}
\label{appendix:sa_scenario_2}
\end{table}

\section{Data}
\label{appendix:data}
In this section, we discuss the data and intermediate steps to calculate the parameters that are used in the model.

\subsection{Residential space heating load}
\label{appendix:heating_load}
The characteristics and space heating load of 400 representative residential homes are simulated using NREL's ResStock Model \citep{Wilson_2017}. The ResStock Model contains a housing database that captures characteristics of millions of residential homes in the U.S. By capturing the housing characteristics using probability distributions based on the home's attributes, the model allows us to generate 400 virtual residential homes that accurately represent the attributes of the actual homes. The selected 400 representative homes in each city used in our model are expected to adequately capture the variability between actual homes in each city, which in total can capture between 88\% to 96\% of the city's annual heating demand \citep{Deetjen_2021}. To capture hourly residential space heating load from ResStock for 400 homes in each city, we run ResStock for these homes in each city, assuming no heat pump installation in these homes. This means we modify ResStock parameters to make sure these homes' space heating loads are satisfied by electronic furnace (COP=1) which is powered by purchased electricity from utilities.

Using a variety of data sources (census data, builders surveys, energy consumption surveys, etc.), ResStock constructed a series of probability distributions for different housing parameters in U.S. cities. For instance, for housing age, ResStock assumes there are a number of houses built before 1940, between 1940 and 1950 and so on. Similar distributions were built for heating fuel used, type of windows, insulation thickness, and many other housing attributes that affect energy use. These distributions were constructed based on correlations between the various parameters - the distribution of housing ages looks very different in the Northeast than the South, for instance, and the amount of insulation will vary depending on how old your house is. ResStock takes all this into account.

To obtain residential space heating load for the 400 single family detached homes of our study in each city, we follow this following framework to run ResStock:

\textbf{Step 1:} Locate the housing characteristics folder in ResStock 2.2.5. This folder has all the housing characteristics in tsv format with the probability distributions: \\\textit{C:\textbackslash resstock-2.2.5\textbackslash project\textunderscore singlefamilydetached\textbackslash housing\textunderscore characteristics}

\textbf{Step 2:} Open the ``Location Region" file, in Excel it would look like this:
\begin{figure}[H]
    \centering
    \includegraphics[scale=1]{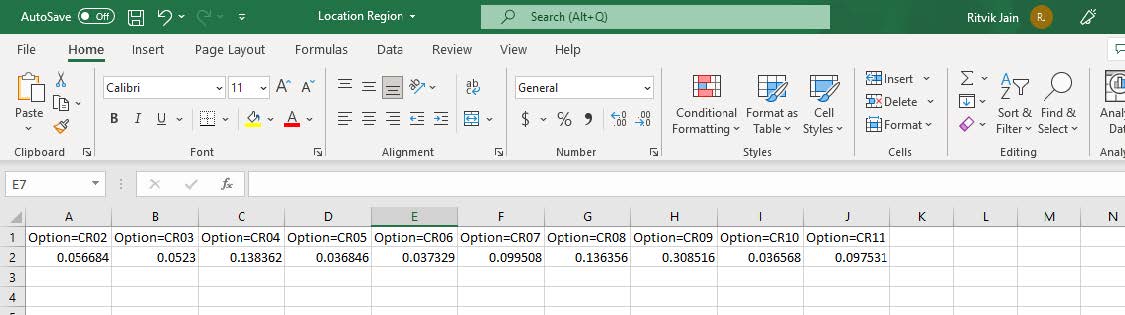}
\end{figure}
The values in row 2 present probabilities that we discussed earlier. The sum of this should always equal 1. CR stands for custom region and this is how CR’s are represented in the US.
\begin{figure}[H]
    \centering
    \includegraphics[scale=0.9]{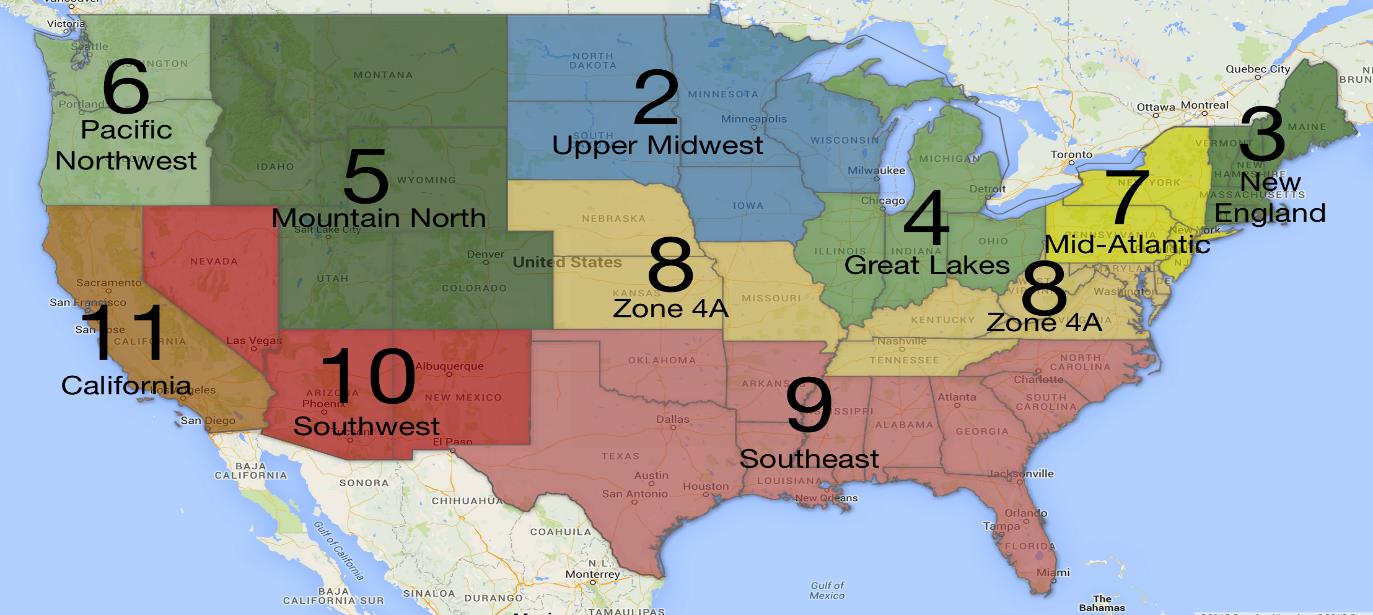}
\end{figure}
For example if we want to generate samples only for California you set Option=CR11 as 1 and all other values as 0. Now, if we want to go further and generate samples from a particular city (for example San Francisco), we can do the following: After setting CR11 as 1 save the .csv file and open the file ``Location” in the same folder. It would look like this in excel:
\begin{figure}[H]
    \centering
    \includegraphics[scale=0.9]{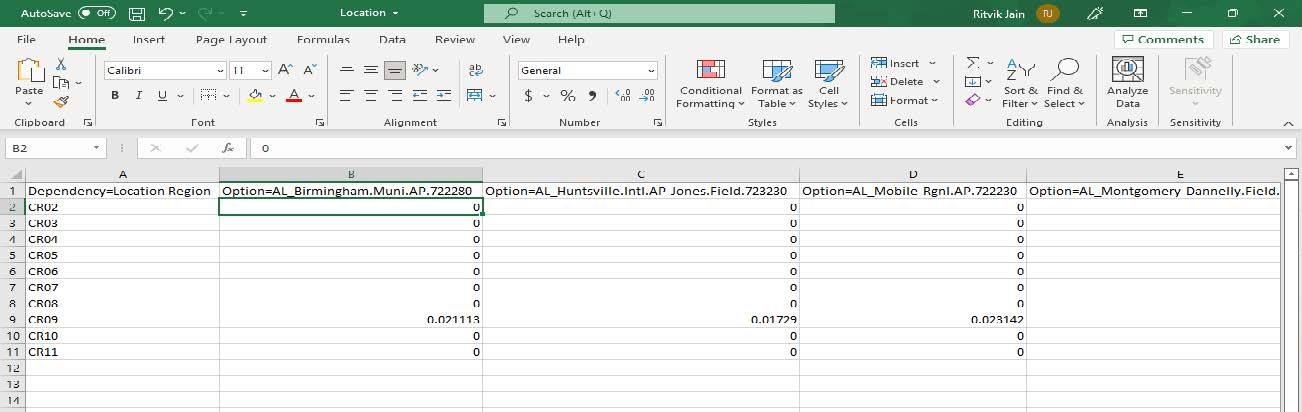}
\end{figure}
First Row and column indicate dependencies, which are very important. In this context it means that ``Location" is dependent on ``Location Region". Dependencies are common in housing characteristics. A detailed characterization of dependencies can be found in \citep{RCW}. Next, search for San Francisco and set its value as 1 and 0 for all other cities in the same custom Region. After doing this if we want to generate a file that would represent 400 homes in San Francisco we would have to do so via Ruby (have to set up and install Ruby, in this case we installed Ruby 2.6.8, and make sure installation of ruby is carried out in the same drive where Resstock 2.2.5 is present (preferably C in Windows)). After setting up Ruby, open command prompt and set the current directory to be where Resstock is located on your drive.
    
\textbf{Step 3:} Run Ruby: \\
Once the location is set, enter the following code in the command prompt Ruby:\\
\textit{resources\slash run\textunderscore sampling.rb -p ``project name” -n “sample size” -o buildstock.csv}

Here is the code for single family detached with 400 samples:\\
\textit{ruby resources\slash run\textunderscore sampling.rb -p project\textunderscore singlefamilydetached -n 400 -o buildstock.csv}

After running this to open the file that has the sample head to: \textit{C:\textbackslash resstock-2.2.5\textbackslash resources} and open the buildstock.csv file to check if the all the samples are from the same city.

\subsection{Scaling-up factors}
\label{appendix:scaling_factor}
Since each home in ResStock represents 242 homes \citep{Wilson_2017}, we estimate relative weight of each city in the national sample by counting the number of homes in that city that ResStock simulates, as below.
\begin{table}[H]
\tablefontsize
\centering
\begin{tabular}{lc} 
\hline 
\hline 
City & Number of actual homes represented by each virtual home\\
\hline
Atlanta & 576.7 \\
Boston & 701.8 \\
Boulder & 113.1 \\
Chicago & 2,964.5 \\
Detroit & 664.3 \\
Dallas & 1,354.6\\
Los Angeles & 3,610.0\\
Minneapolis & 448.9 \\
New York & 8,509.9\\
Orlando & 311.0\\
Phoenix & 1,507.7\\
Seattle & 803.4\\
\bottomrule 
\end{tabular} 
\caption{Number of actual homes represented by each virtual home simulated in ResStock in 12 U.S. Cities of our study.}
\label{tab:scale_up_factor}
\end{table}

\subsection{Residential retail rates}
\label{appendix:elec_prices}
The time of use residential utility rates of electricity are taken from the most updated electric rate books of the utilities that are in charge of providing electricity to the buildings in the cities of our study. Each utilities have a wide variety of rates that the homes can opt-in. In this study, we use the most popular, usually default time-of-use rates in each utility that vary depending on the month and hours. 

All the time-of-use retail rates of the 12 U.S. cities of our studies are as following:
\begin{table}[H]
\tablefontsize
\centering
\begin{tabular}{lccccc} 
\hline 
\hline 
City & On-peak rate & On-peak rate  & Off-peak rate & On-peak hours & Utility \\
 & -Summer (\$/kWh) & -Winter (\$/kWh) & (\$/kWh) &  (weekday) &  \\
  & (mid-peak rate) & (mid-peak rate) &  & &  \\
\hline
Atlanta & 0.117993 & 0.117993 & 0.012614 & 2pm-7pm & Georgia Power \cite{GW_2023} \\
\hline
Boston & 0.28783 & 0.28783 & 0.13477 & 7am-8pm & Eversource \cite{Eversource_2023} \\
 \hline
Boulder & 0.19 & 0.19 & 0.12 & 3pm-7pm (high) &  Xcel \cite{Xcel_CO_2023} \\
 & (0.15 mid) & (0.15 mid) & & 1pm-3pm (mid) &  \\
\hline
Chicago & 0.06117 & 0.06117 & 0.02447 &2pm-7pm (high) & ComED \cite{ComED_2023} \\
 & (0.03970 mid) & (0.03970 mid) &  & 6am-2pm \& &  \\
  &  & & &7pm–10pm (mid) &  \\
\hline
Detroit & 0.1809 & 0.2240 & 0.1673 & 3pm-7pm & DTE \cite{DTE_2023}\\
 & (June-September) & (October-May) &  &   & \\
\hline
Dallas & 0.245241 & 0.245241 & 0.077926 & 1pm-7pm & Xcel \cite{Xcel_TX_2023}\\
 & (June-September) & (October-May) &  &  & \\
\hline
Los Angeles & 0.21659 & 0.21659 & 0.16826 & 10am-8pm & LADWP \cite{LADWP_2023} \\
\hline
Minneapolis & 0.25879 & 0.21408 & 0.05171 & 9am-9pm & Xcel \cite{Xcel_MN_2023} \\
 & (June-September) & (October-May) &  & & \\
\hline
New York & 0.3305 & 0.1223 & 0.0233 & 8pm-12am & Con Edison  \cite{ConEd_2023}\\
 & (June-October) & (November-May) &  &  & \\
 \hline
Orlando & 0.08828 & 0.08828 & 0.06520 & 7am-10am \& & OUC \cite{OUC_2023}\\
 & (April-October) & (November-March) &  & 6pm-9pm (high) & \\
  &  &  &  &  10am-6pm (mid) & \\
  \hline
Phoenix & 0.34396 & 0.32543 & 0.12345 & 4pm-7pm & APS \cite{APS_2023} \\
 & (May-October) & (November-April) &  &  &  \\
 \hline
Seattle & 0.15 & 0.15 & 0.08 & 5pm-9pm (high) & Seattle City Light \cite{SCL_2023} \\
 & (0.13 mid) & (0.13 mid) &  & 6am-5pm \&  & \\
  &  &  &  & 9pm-12am (mid) & \\
\bottomrule 
\end{tabular} 
\caption{Time-of-Use retail rates in 12 U.S. cities of our study.}
\label{tab:retail_rates}
\end{table}

\subsection{ASHP's coefficient of performance}
\label{appendix:hp_COP}
\begin{figure}[H]
    \centering
    \includegraphics[scale=0.9]{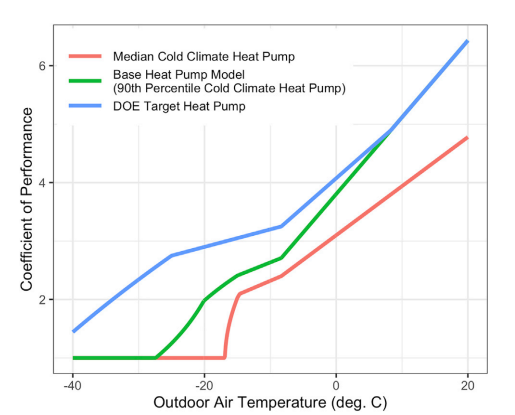}
    \caption{Heat Pump Model Temperature-Dependent Coefficient of Performance  \citep{Waite_and_Modi_2022}}
    \label{fig:COP}
\end{figure}

\begin{figure}[H]
    \centering
    \includegraphics[scale=0.6]{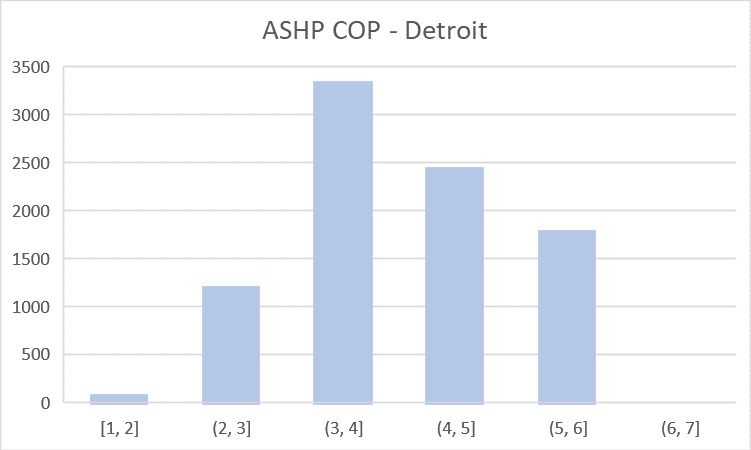}
    \caption{Hourly Heat Pump COP in Detroit}
    \label{fig:COP_detroit}
\end{figure}

Heat pump's coefficient of performance for a given time (hour) measures how efficient the heat pump is during that particular hour. It is calculated based on hourly outdoor temperature at the building's location, using a model of heat pump's temperature-dependent COP taken from \citep{Waite_and_Modi_2022} (Figure \ref{fig:COP}). For this study, we assume median cold climate heat pump (red line). In Detroit, most hours have COP between 3-4, with 6.1 being the highest COP in an hour (Figure \ref{fig:COP_detroit}).  
%However, \citet{Waite_and_Modi_2022} also models two additional types of heat pumps using a heat pump database with more than 1,000 available heat pumps: the 90th percentile performance heat pumps (green line), and the Department of Energy Target heat pump (blue line).
%

\subsection{Fan energy consumption}
\label{appendix:hp_fan}
\begin{figure}[H]
     \includegraphics[scale=0.25]{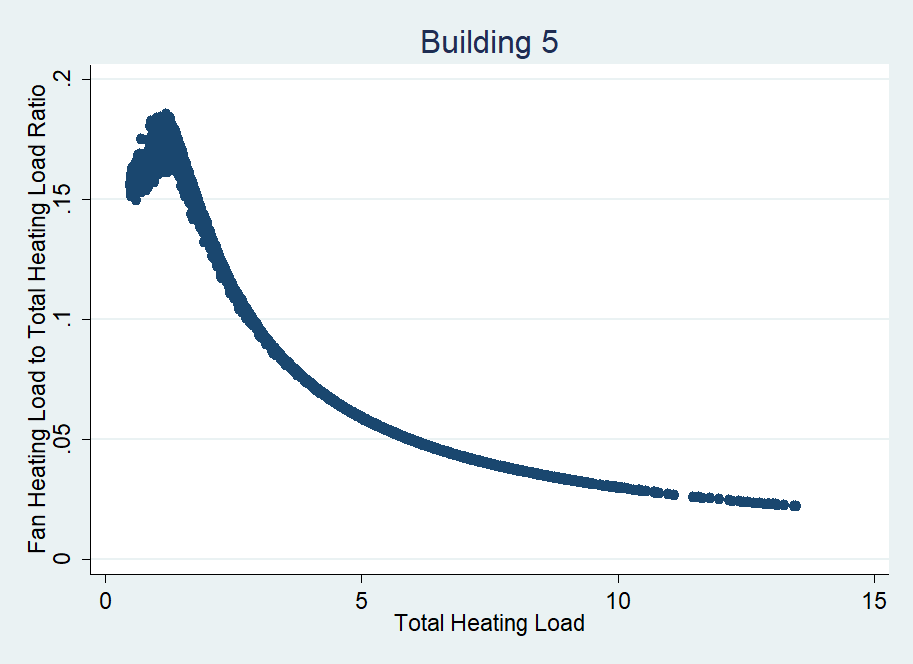} \includegraphics[scale=0.25]{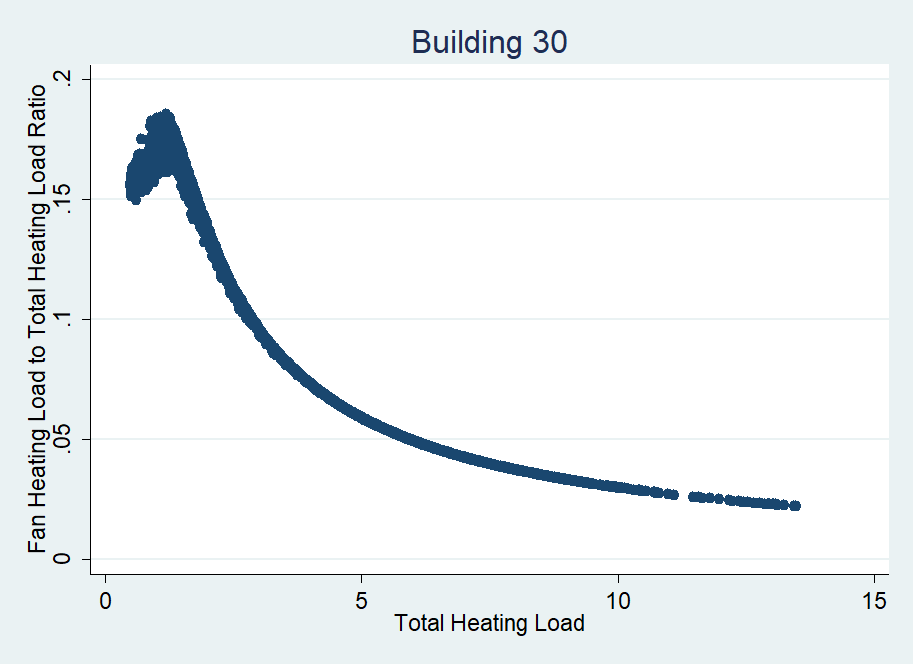}
     \includegraphics[scale=0.25]{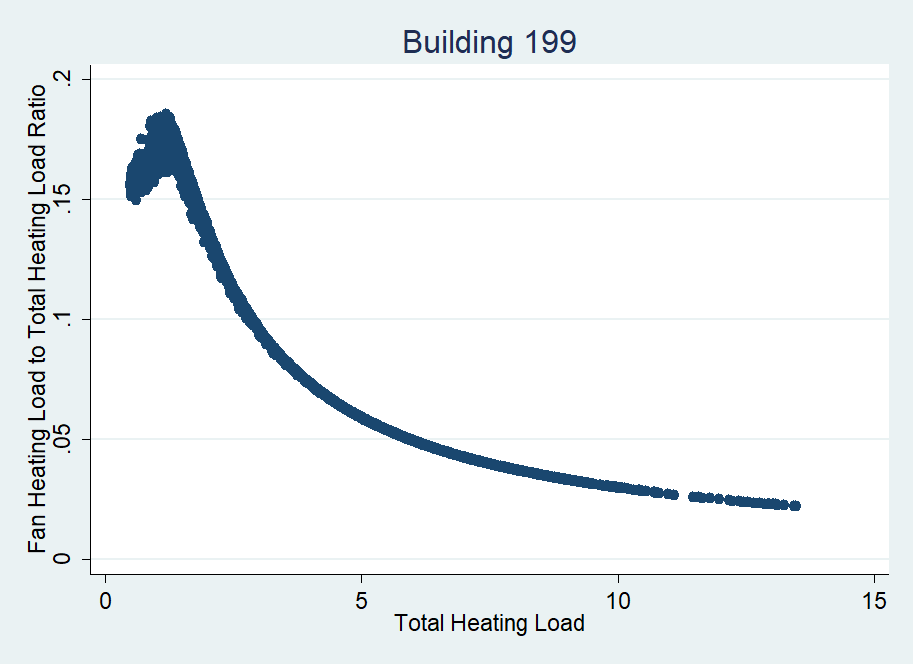}
     \includegraphics[scale=0.25]{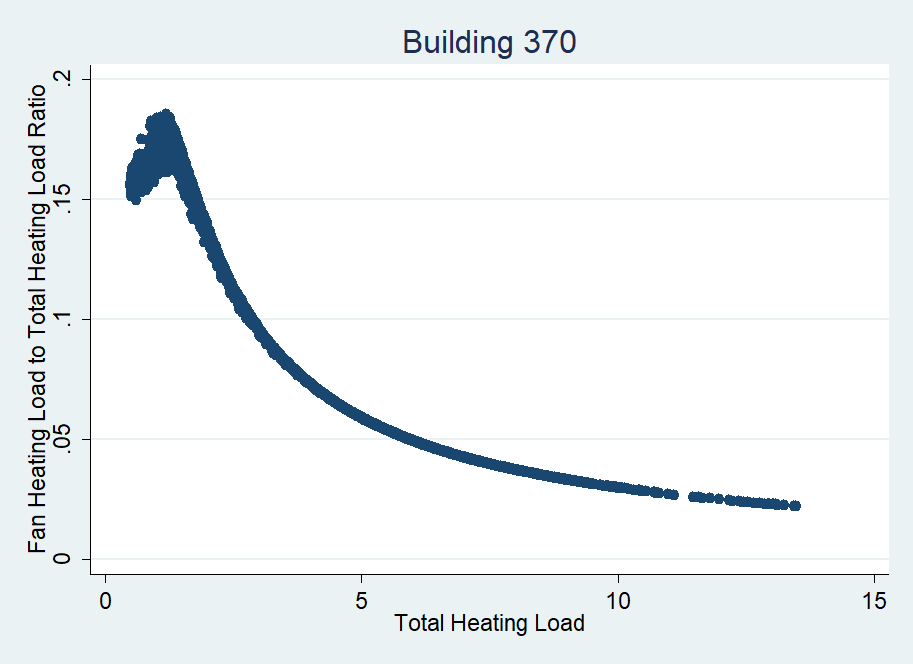}
    \caption{Scattered Plots of Fan Heating Load-Total Heating Load Ratio and Total Heating Load in four Detroit homes of different load profiles.}
    \label{fig:example_fanload}
\end{figure}

\begin{figure}[H]
    \centering
    \includegraphics[scale=0.4]{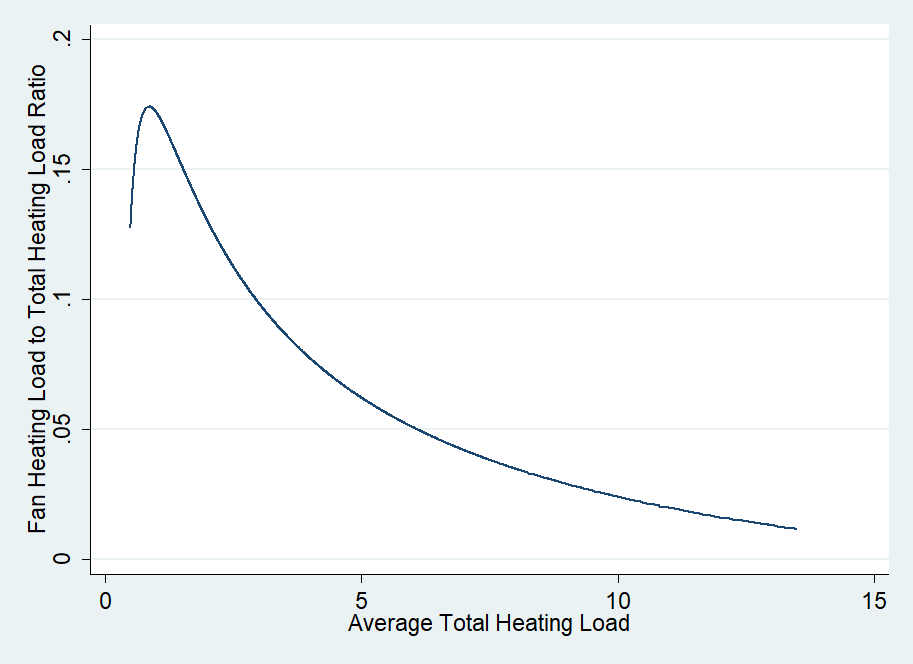}
    \caption{Fitted Polynomial Function of Fan Heating Load-Total Heating Load Ratio and Total Heating Load}
    \label{fig:fan_heating_load_ratio}
\end{figure}

To capture fan heating load for both ASHP and TES, we run ResStock again, following the steps in \ref{appendix:heating_load}, for the same 400 homes in Detroit, assuming each home has installed ASHP. To do this, we modify ResStock parameters to make sure these homes' heating loads are supplied by ASHP with electric baseboard as backup. We then calculate the ratio of total fan heating load to total heating load for each home and observe that in all homes, the relationship between this ratio and the home's total heating load follows a polynomial function (Figure \ref{fig:example_fanload}). Therefore, we then calculate the average hourly total heating loads and fan heating loads across these homes and perform polynomial regression analysis to fit fan heating ratio against total heating load ($R^2 = 0.92$) (Figure \ref{fig:fan_heating_load_ratio}). We assume this relationship between total heating output and fan heating load is the same for TES discharging process.

\subsection{TES materials}
\label{appendix:tes_materials}
The reaction rate equations used to generate the Ragone plots are curve fits based on thermo-gravimetric analysis (TGA) experiments as listed in Table \ref{tab:Raw_TGA}. The TGA experiments are of the hydration step all beginning from a fully dehydrated state and finishing with a fully hydrated/reacted salt \cite{Gaeini2019}. The rate equations vary based on set salt temperatures and relative humidity, so these are listed for reference. The instantaneous specific power, $sP$, 
 is calculated according to equation \ref{eq:sp_SOC} 

\begin{equation}
    sP = \frac{\partial SOC(t)}{\partial t} enthalpy_s
    \label{eq:sp_SOC}
\end{equation}

\begin{figure}[H]
    \centering
    \includegraphics[scale=0.8] {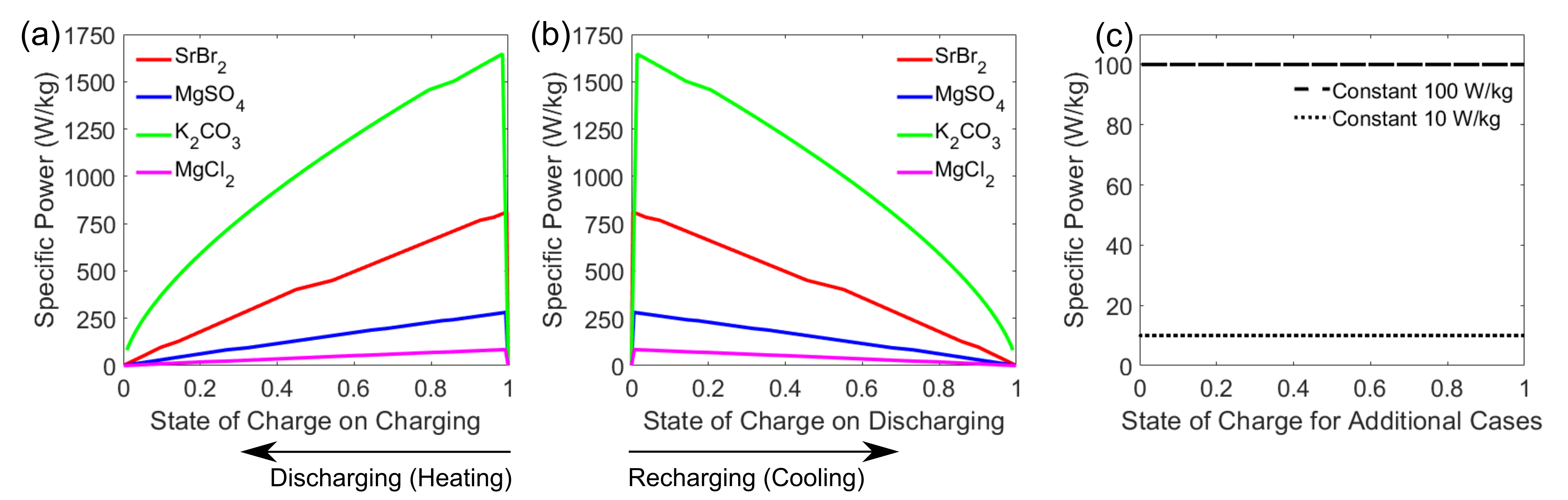}
    \caption{Ragone Plots of Salt Hydrates upon (a) discharging/heating and (b) recharging/cooling. (b) is a mirror image of (a). Solid line curves come from thermo-gravimetric analysis (TGA) data. In (c) the dotted lines represent constant power outputs more typically seen in reactors.}
    \label{fig:Ragone}
\end{figure}

When equation \ref{eq:sp_SOC} is combined with the data from Table \ref{tab:Raw_TGA} and Table \ref{tab:tes_salts} this results in the Ragone plots of Figure \ref{fig:Ragone} (a,b). The constant power outputs of 10 and 100 W/kg are shown in Figure \ref{fig:Ragone}(c).

\begin{table}[H]
\tablefontsize
\centering
\begin{tabular}{clc}
\hline 
\hline 
 Salt Hydrate, ref & Rate Equation (t in min) & Conditions\\
\hline
\ce{SrBr2}, \citep{Cammarata2018} & $SOC = exp(-3.81 \times 10^{-2} t)$ & $25^{\circ} C,\ 50\% $ \ RH \\
\ce{MgSO4}, \citep{Linnow2014}  & $SOC = exp(-6.28 \times 10^{-3} t)$ & $25^{\circ} C,\ 80\% $ \ RH \\
\ce{K2CO3}, \citep{Gaeini2019} & $SOC = (3-0.182t)^3/27$ & $26^{\circ} C,\ 39\% $ \ RH \\
\ce{MgCl2}, \citep{Fisher2021} & $SOC = exp(-7.35 \times 10^{-3} t)$ & $25^{\circ} C,\ 24\% $ \ RH \\
\bottomrule 
\end{tabular} 
\caption{Salt Hydrate Rate Data for Discharging from Thermo-gravimetric Analysis Experiments}
\label{tab:Raw_TGA}
\end{table}

\section{TES Sizing}
\label{appendix:tes_sizing}
The size of TES in each building in our model is the mass (in kg) of each of the salt hydrates used as TES materials. There are multiple ways to size TES. In this paper, we use two ways: 1) Variable sizing, in which TES for each salt is sized based on annual peak load of each building, 2) Incremental sizing, in which TES for each salt is sized based on annual peak load of each building then rounded up the the nearest 25kg, and 3) Fixed sizing, in which TES size is fixed to one mass for all salts and all buildings.
\subsection{Variable TES sizing}
\label{appendix:opt_sizing}
\begin{figure}[H]
    \centering
    \includegraphics[scale=0.6]{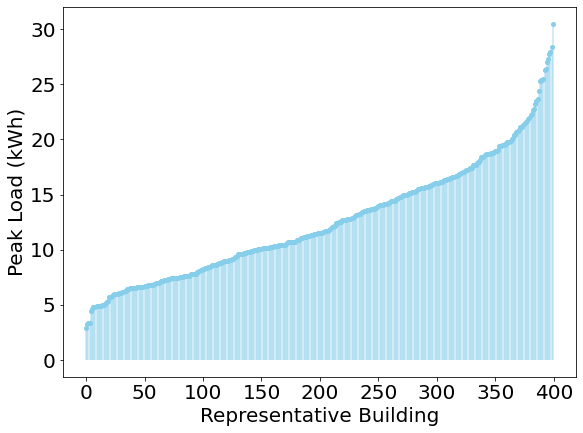}
    \caption{2018 Peak Loads of 400 Representative Residential Homes in Detroit, MI}
    \label{fig:peakload}
\end{figure}

For variable sizing for a city, we first look at all the buildings annual peak load (Figure \ref{fig:peakload} shows annual peak loads for all representative buildings in Detroit, MI, for example). 

Based on these peak loads, TES sizes (in kg of salt mass) are determined based on the power requirement and energy requirement (reaction enthalpy) to reach these peak loads, that are unique to each salt hydrate used as TES material Table \ref{tab:tes_salts}. The salt must be the size of the larger of the two requirements.

\begin{figure}[H]
    \centering
    \includegraphics[scale=0.35]{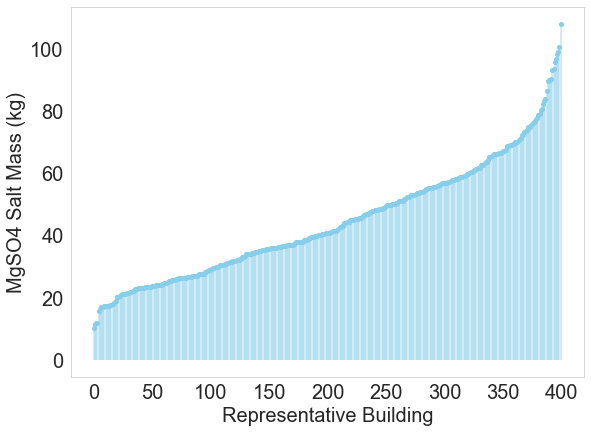}
    \includegraphics[scale=0.35]{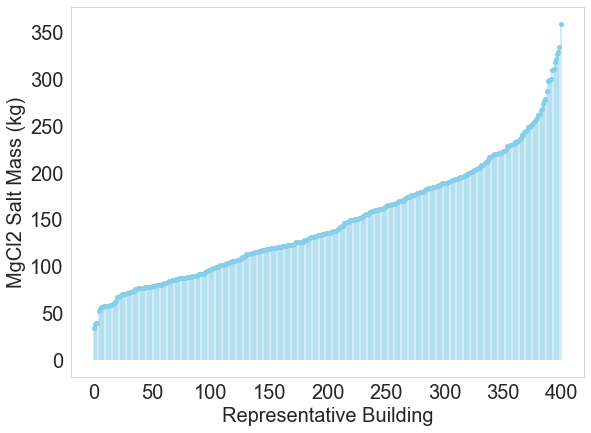}\\
    \includegraphics[scale=0.35]{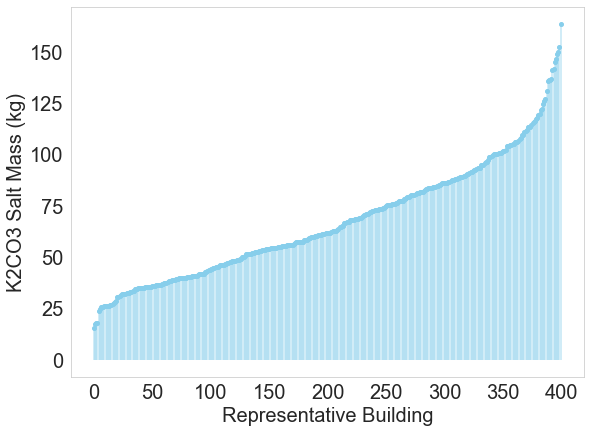}
    \includegraphics[scale=0.35]{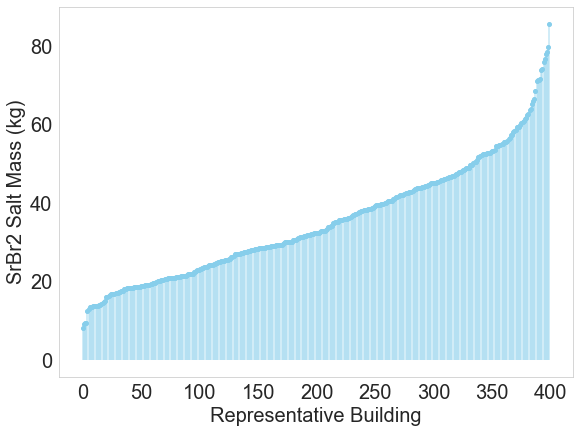}
    \caption{TES sizes for four salts for 400 representative buildings in Detroit, MI, using optimal sizing method.}
    \label{fig:salt_mass_opt}
\end{figure}

Therefore, the TES size for each type of salt for each building is:
\begin{align}
    mass^{TES}_{js} = \max \left\{\dfrac{D^{Peak}_{j}}{enthalphy_{s}}, \dfrac{D^{Peak}_{j}}{powerSpec_{s}}\right\}
\end{align}
where $mass^{TES}_{js}$ the the size of TES using salt $s$ in building $j$, $D^{Peak}_{j}$ is annual peak demand of building $j$, $enthalphy_{s}$ is reaction enthalphy of salt $s$, and $powerSpec_{s}$ is power specific of salt $s$ to reach the peak demand. Using this TES sizing method, the TES sizes for each salt for each building in Detroit, MI, are shown in Figure \ref{fig:salt_mass_opt}.

\subsection{Incremental TES sizing}
\label{appendix:var_sizing}
\begin{figure}[H]
    \centering
    \includegraphics[scale=0.35]{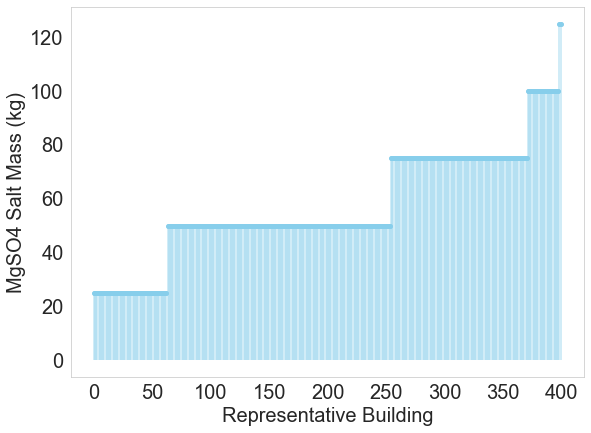}
    \includegraphics[scale=0.35]{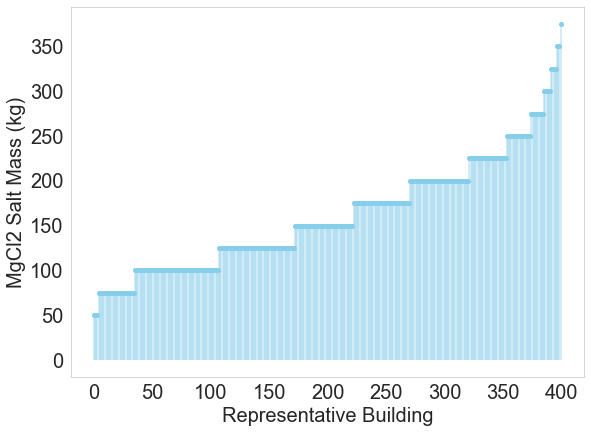}\\
    \includegraphics[scale=0.35]{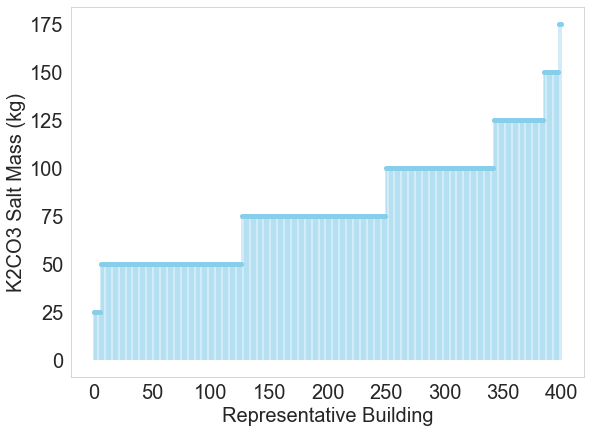}
    \includegraphics[scale=0.35]{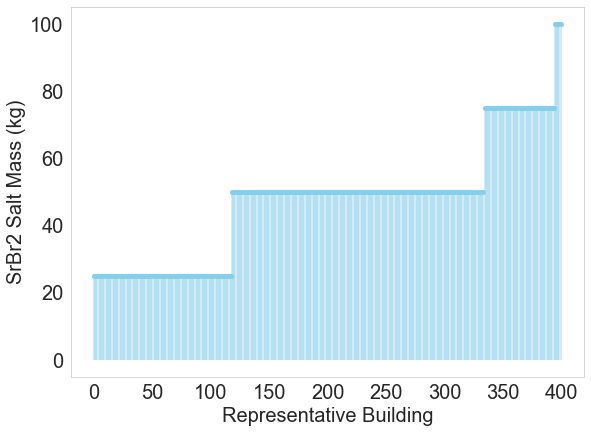}
    \caption{TES sizes for four salts for 400 representative buildings in Detroit, MI, using variable sizing method.}
    \label{fig:salt_mass}
\end{figure}

For incremental sizing, we round up the salt mass resulted using optimal sizing to the nearest 100 kg. Using this TES sizing method, the TES sizes for each salt for each building in Detroit, MI, are shown in Figure \ref{fig:salt_mass}.

\subsection{Fixed TES sizing}
\label{appendix:fixed_sizing}
Another way to size TES is to fix the TES size for all buildings in a city for all salts. This sizing method allows us to compare the effectiveness of these salts as TES materials better. Based on variable sizing method, the average TES sizes for all buildings in Detroit are 100.5 kg, 195.25 kg, 114.5 kg, and 100 kg for \ce{MgSO4}, \ce{MgCl2}, \ce{K2CO3}, and \ce{SrBr2} respectively, we assume the fixed TES size is 150 kg for all salts for all buildings.

\section{Model Code and Data Availability}
\label{code}
Model code and data are available at \url{https://github.com/atpham88/TES}. 

\newpage
\section{Lists of Model Variables, Parameters, and Full Model Formulation}
\begin{table}[H]
\caption{List of Variables}
\label{var_list}
\tablefontsize
\linespread{1.5} \selectfont
\centering
\begin{tabular}{clc} 
\hline 
\hline 
Variable & Definition & Unit \\
\hline
$k^{TES,D}_{jt}$ & Discharging Power rating for TES in building $j$ at time $t$  & KW \\
$k^{TES,C}_{jt}$ & Charging Power rating for TES in building $j$ at time $t$  & KW \\
$g^{TES}_{jt}$ & TES's energy output in building to shift load in building $j$ at time $t$ & kWh \\
$g^{HP-TES}_{jt}$ & Heat pump's energy output to TES in building $j$ at time $t$ & kWh \\
$g^{HP-L}_{jt}$ & Heat pump's energy output to serve load in building $j$ at time $t$ & kWh \\
$d^{HP}_{jt}$ & Purchased electricity to power heat pump/resistance heater in building $j$ at time $t$ & MWh \\
$x^{TES}_{jt}$ & State of charge for TES in building $j$ at time $t$ & kWh \\
%$v^{TES}_{{j}t}$ & Whether TES is discharging (1) or not (0) in building $j$ at time $t$ & Binary \\
%$v^{HP-TES}_{{j}t}$ & Whether heat pump is charging TES (1) or not (0) in building $j$ at time $t$ & Binary \\
\bottomrule 
\end{tabular} 
\end{table}

\begin{table}[H]
\caption{List of Parameters and Sets}
\label{param_list}
\tablefontsize
\linespread{1.5} \selectfont
\centering
\begin{tabular}{clc} 
\hline 
\hline 
Parameter/Set & Definition & Unit/Value \\
\hline
\textit{Parameters:} \\
$F^d$ & Discharging efficiency & \% \\
$F^c$ & Charging efficiency & \%\\
$F^s$ & Humidification parasitic TES load using salt $s$ & \%\\
$(x_A,y_A), (x_B, y_B)$ & Cut-off points for piece-wise function of TES power rating and SOC  & (\%, kW)\\
$E^{TES}_{j}$ & Energy capacity for TES in building $j$ & KWh \\
$P^{R}_{jt}$ & Residential electricity price at location of building $j$ at time $t$ & \$/kWh \\
$T^{O}_{t}$ & Outdoor temperature at location of building $j$ at time $t$ & \degree C \\
$K_j^{HP}$ & Heat pump capacity in building $j$  & KWh \\
$COP_{jt}$ & Heat pump Coefficient of performance in building $j$ at time $t$ & \degree C \\
$D_{jt}$ & Heating load in building $j$ at time $t$ & KWh \\
\hline \\ \hline
\textit{Sets:} \\
$\mathbb{J}$ & Set of buildings, index $j$ & --\\
$\mathbb{T}$ & Set of hours in a typical year, index $t = \{1,2,3,...,8760\}$ & --\\
$\mathbb{Z}_{2}$ & Set of binary numbers, $\mathbb{Z}_{2} = \{0,1\}$ & --\\
$\mathbb{Z}_{1}$ & Set of real numbers $\geq 1$, $\mathbb{Z}_{1} = \{x \in \mathbb{Z}_{1}: x\geq 1\}$ & --\\
\bottomrule 
\end{tabular} 
\end{table}

\section{Additional Results}
\subsection{TES operation across different TES materials (Reference case)}
\begin{figure}[H]
     \centering
     \begin{subfigure}[c]{1\textwidth}
         \caption{\ce{MgSO4}-based TES}
         \includegraphics[scale=0.42]{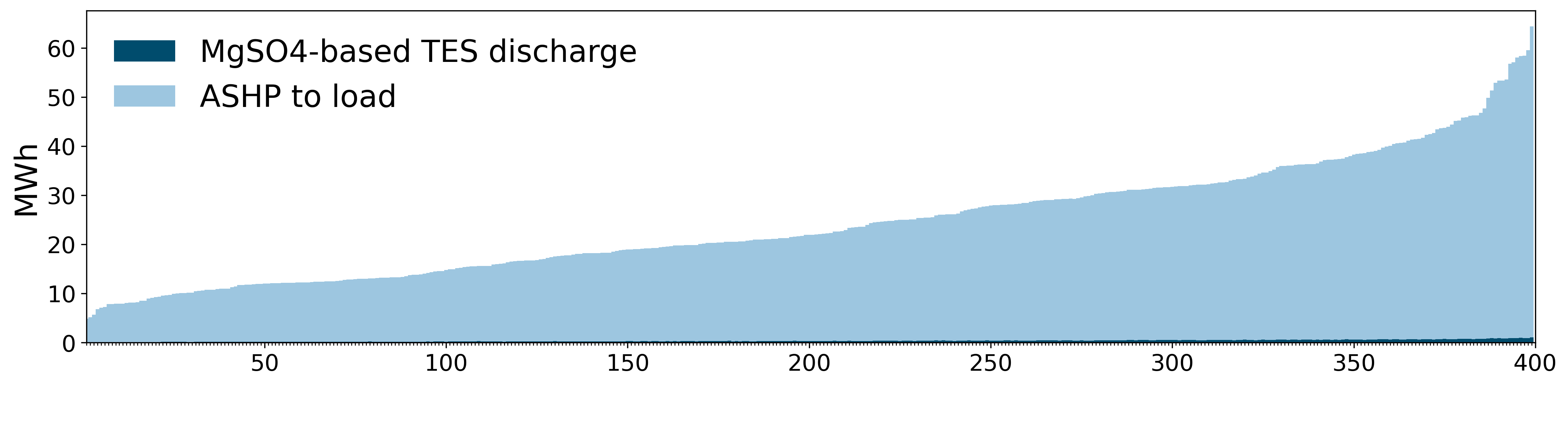} 
         \label{fig:tes_output1_opt}
     \end{subfigure}
     \begin{subfigure}[c]{1\textwidth}
         \caption{\ce{MgCl2}-based TES}
         \includegraphics[scale=0.42]{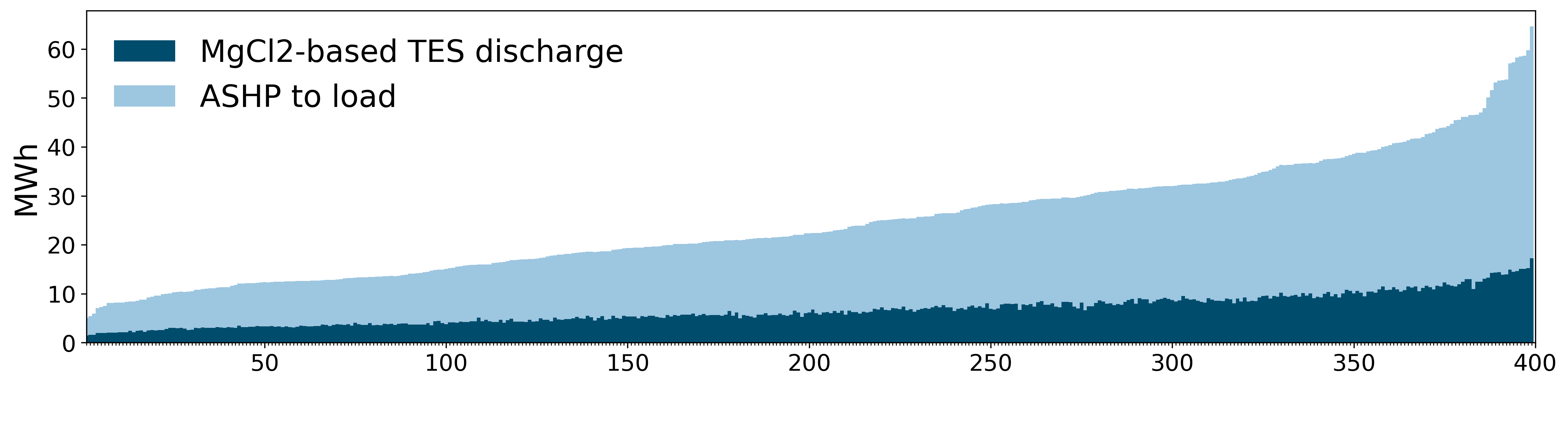}
         \label{fig:tes_output2_opt}
     \end{subfigure}
     \begin{subfigure}[c]{1\textwidth}
         \caption{\ce{K2CO3}-based TES}
         \includegraphics[scale=0.42]{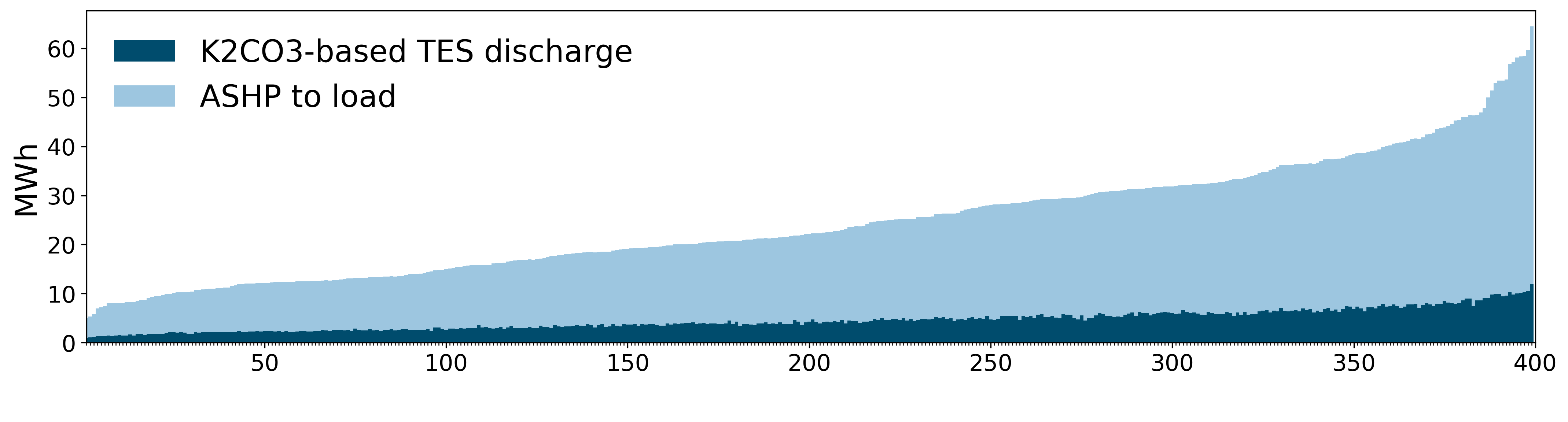}
         \label{fig:tes_output3_opt}
     \end{subfigure}
     \begin{subfigure}[c]{1\textwidth}
         \caption{\ce{SrBr2}-based TES}
         \includegraphics[scale=0.42]{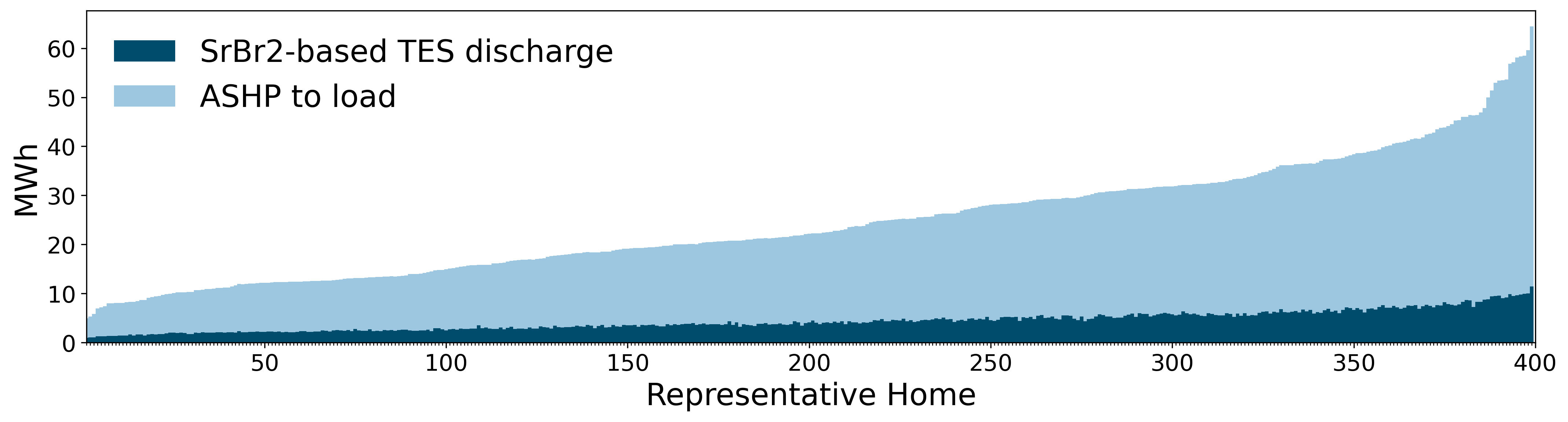}
         \label{fig:tes_output4_opt}
     \end{subfigure}
        \caption{TES's energy outputs as fraction of total annual loads across 400 representative Detroit homes.}
        \label{fig:tes_output_opt}
\end{figure}

\begin{figure}[H]
     \centering
     \begin{subfigure}[c]{1\textwidth}
         \centering
         \includegraphics[scale=0.41]{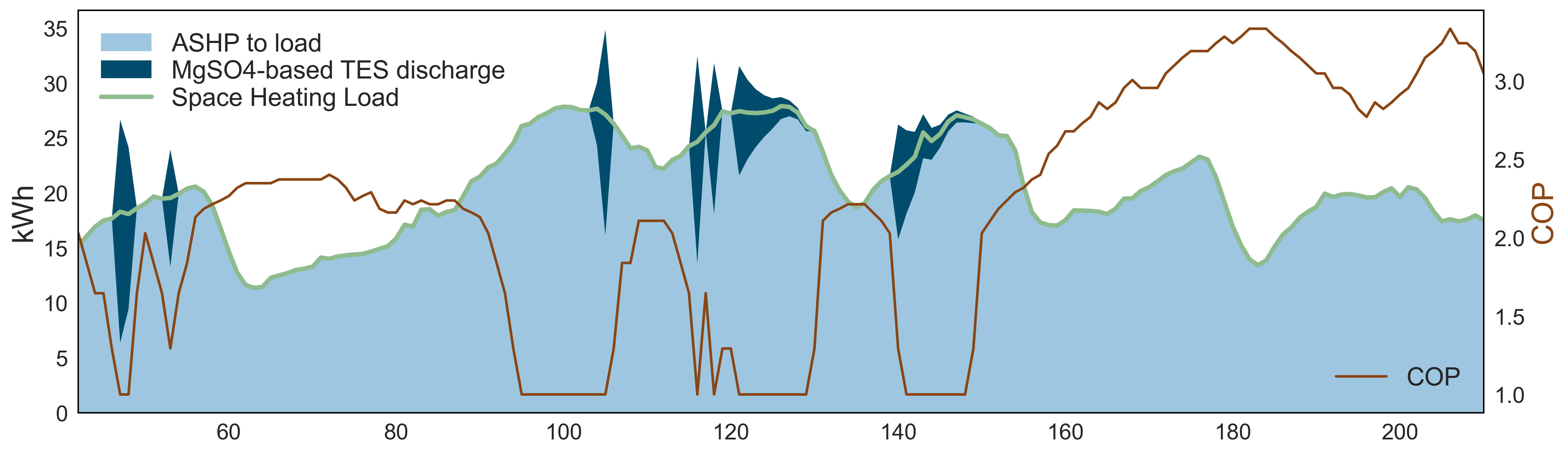}
         \label{fig:tes_discharge_opt1}
     \end{subfigure}
     \begin{subfigure}[c]{1\textwidth}
         \centering
         \includegraphics[scale=0.41]{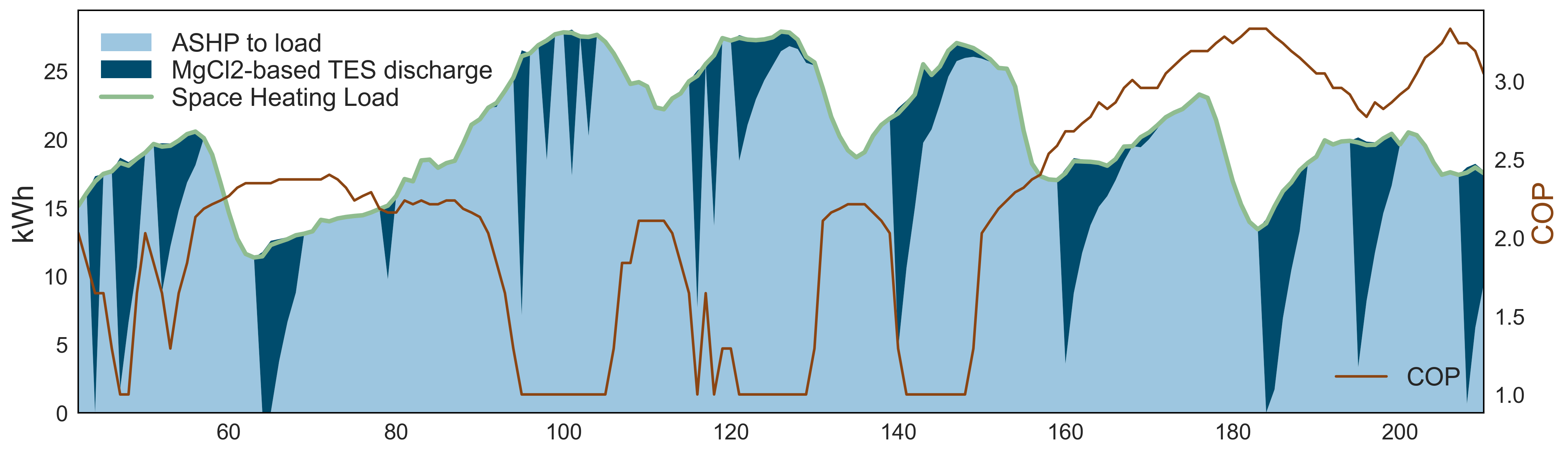}
         \label{fig:tes_discharge_opt2}
     \end{subfigure}
     \begin{subfigure}[c]{1\textwidth}
         \centering
         \includegraphics[scale=0.41]{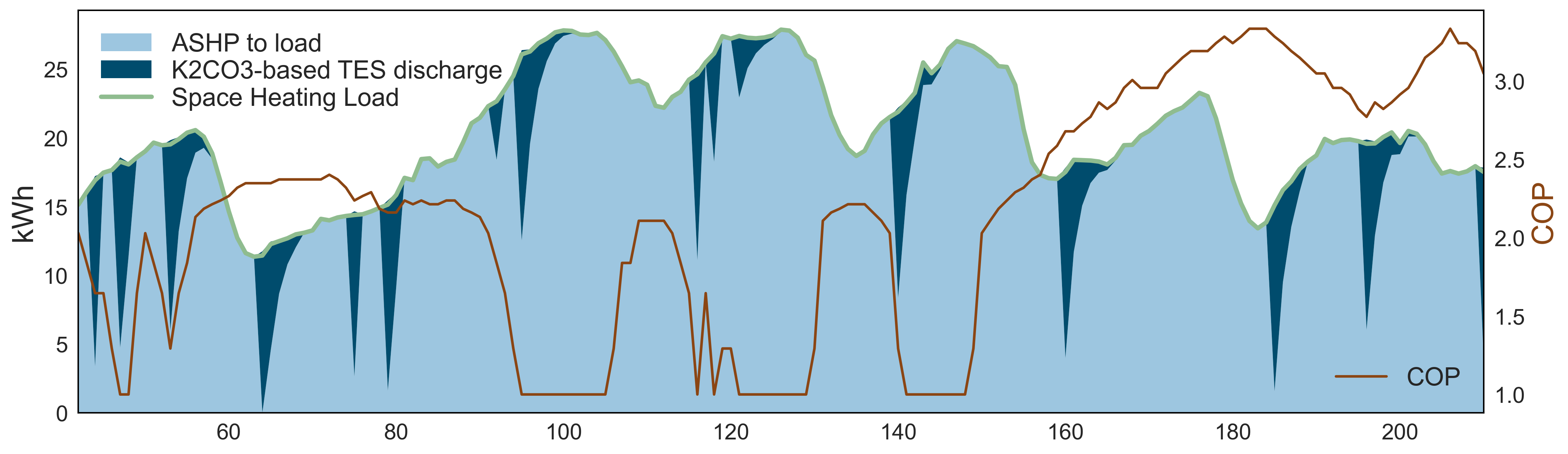}
         \label{fig:tes_discharge_opt3}
     \end{subfigure}
     \begin{subfigure}[c]{1\textwidth}
         \centering
         \includegraphics[scale=0.41]{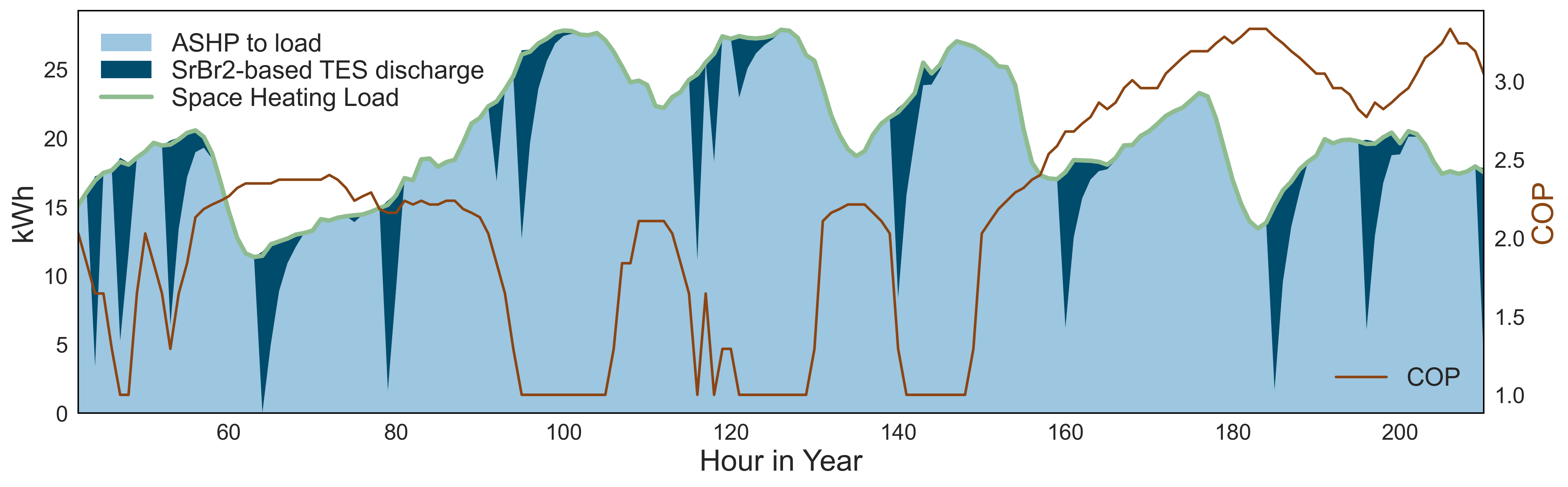}
         \label{fig:tes_discharge_opt4}
     \end{subfigure}
        \caption{TES's energy outputs as fraction of total annual loads across 400 representative homes in Detroit with variable TES sizing. The amount of \ce{MgSO4}-based TES discharge above load is the parasitic load due to \ce{MgSO4}'s humidification requirement during discharge.}
        \label{fig:tes_discharge_opt}
\end{figure}

\begin{figure}[H]
     \begin{subfigure}[c]{1\textwidth}
     \centering
         \caption{Operation of ASHP and \ce{SrBr2}-based TES}
         \includegraphics[scale=0.41]{Figures/tes_discharge_SrBr2.png}
         \label{fig:tes_discharge1}
     \end{subfigure}
     \hspace{-0.5in}
     \begin{subfigure}[c]{1\textwidth}
         \caption{Comparing charging and discharging of \ce{MgCl2} and \ce{SrBr2}-based TES}
         \includegraphics[scale=0.41]{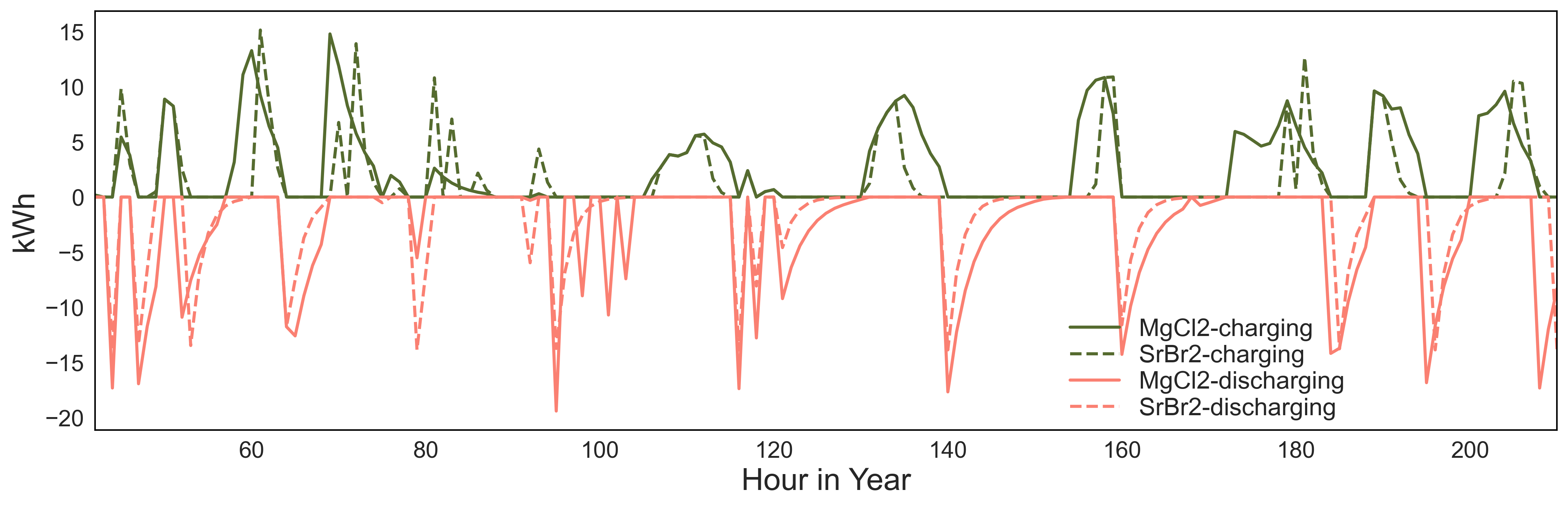}
         \label{fig:tes_discharge2}
     \end{subfigure}
        \caption{TES and ASHP's operations to serve space heating load in a high load Detroit home during the highest peak heating load week in 2018.}
        \label{fig:tes_discharge}
\end{figure}

\subsection{ASHP and TES operation under peak load shifting maximization}

\begin{figure}[H]
\centering
\includegraphics[scale=0.41]{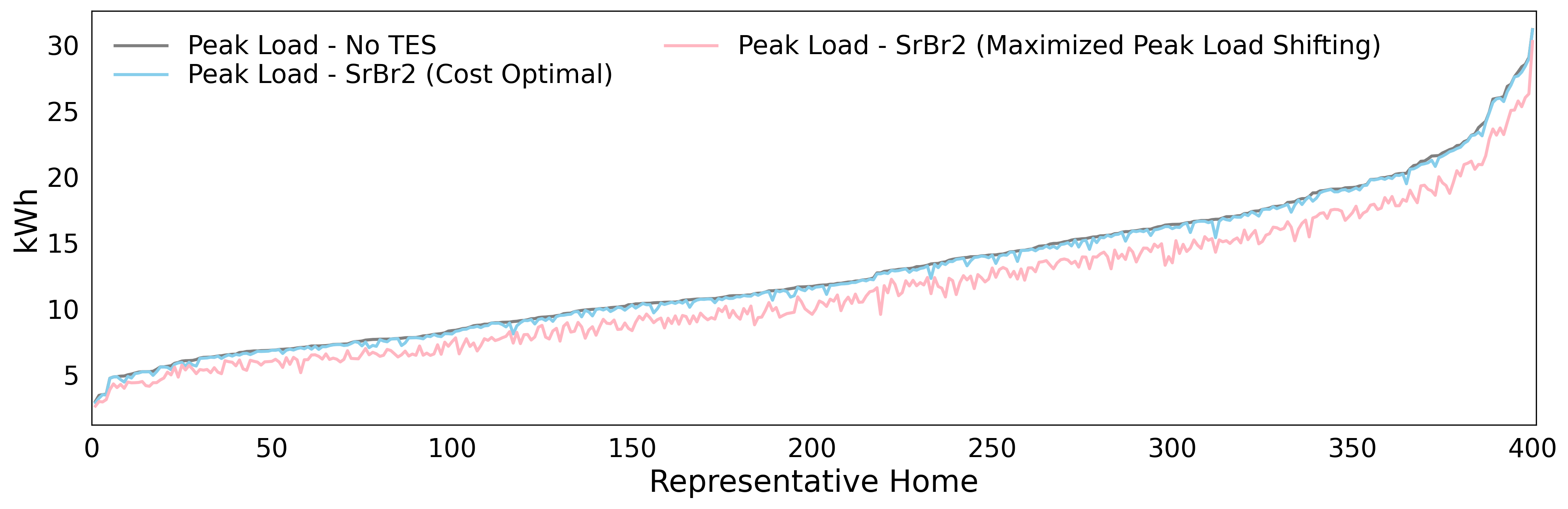}
\caption{Peak load reduction in 400 representative homes in Detroit, MI, from coupling \ce{SrBr2}-based TES ofwith ASHP during 2018's peak demand week, under cost minimization and peak heating load shifting maximization.}
\label{fig:load_reduction}
\end{figure}

\subsection{Economic value of TES}
\begin{figure}[H]
\centering
     \centering
     \includegraphics[scale=0.41]{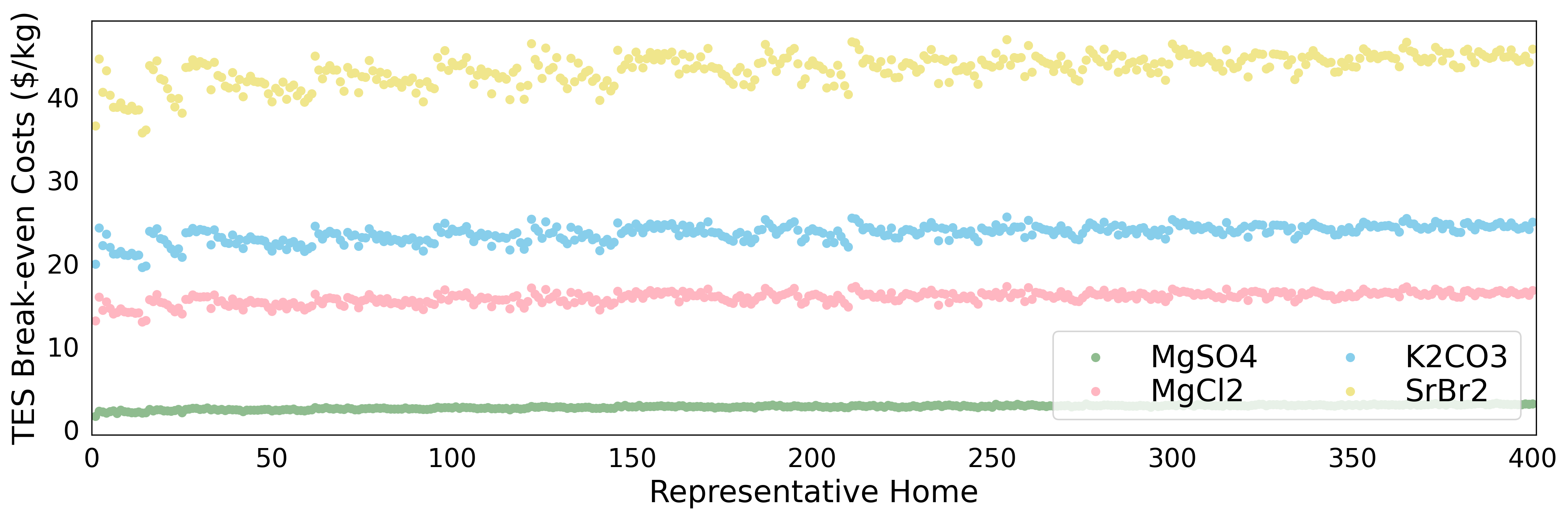}
        \caption{Break-even costs of \ce{BrSr2}-based TES.}
        \label{fig:break_even_kg_var}
\end{figure}

\subsection{Sensitivity Analysis}
\label{appendix:sa}

In our main results, we size TES variably based on each home's peak heating load. However, commercial TES systems would likely have a finite set of sizes. We approximate this commercialized future using two sizing methods: (1) incremental TES sizing (Figure \ref{fig:salt_mass}), where we size salt masses at each home's peak load and then round up to the nearest 25 kg, and (2) fixed TES sizing, where we size salt masses at 150 kg for all homes, which is roughly the average salt mass across 400 Detroit homes across salts when TES is sized variably. Under both sizing approaches, TES is oversized relative to peak heating loads, resulting in higher cost savings but lower TES break-even costs compared to variable sizing (i.e., our main) results. Specifically, under incremental and fixed sizings respectively, TES raises annual cost saving across household by up to 23\% and by 71\% but lowers beak-even cost by up to 9\% and 42\% compared to under variable (optimal) sizing.

We ran additional scenarios for the most promising salt in terms of break-even cost, \ce{SrBr2}, with a constant specific power output\citep{Zhang2020}. Compared to variable power rating based on Ragone plot, TES with constant specific power of 100 W/kg, which falls below the average power density (135 W/kg) in the Ragone plot, can increase annual cost saving by \$452 (4\%) across 400 representative Detroit homes. Greater cost savings are driven by discharging and charging at constant power rating at any given time. Average break-even cost across homes for \ce{SrBr2}-based TES with 100 W/kg constant power specific is \$16/kWh thermal, which is \$0.6 increase (3.4\%) compared to variable power rating TES.

Above, we examined a TES-ASHP system backed up by electric baseboard heat to take advantage of the resistance coil within the ASHP unit. However, ASHPs might instead use gas furnaces as backup heating systems. We run an additional scenario in which a \ce{SrBr2} TES-ASHP system is backed up by a gas furnace that is scheduled to be turned on when the ASHP COP falls below 2.7, or when the outdoor temperature falls below 25\degree F (-4\degree C). Across our 400 Detroit homes, on average, using gas furnace as backup lowers TES cost saving from 4.1\% to 3.1\%, or 1.0 percentage points. Annual city-wide TES cost saving with gas furnace back up is \$12.1 million, compared to \$16.4 million with electric backup, a 26\% decrease.

Finally, the above results assume TES round-trip efficiency of 98\% and capability of recirculating indoor air, resulting in only \ce{MgSO4}-based TES requiring humidification of incoming air during discharge process. We run additional scenarios with lower TES round-trip efficiency and higher humidification needs that result in higher parasitic loads between 44\% and 70\% across salts to test the robustness of our results under more restricted TES designs. While increasing non- humidification parasitic TES losses from 2\% to 5\% and 10\% reduces cost savings from TES by up to \$1.6 million (10\%) and \$3.6 million (22\%) in Detroit, a high humidification parasitic load can significantly reduce TES cost saving by up to 88\%. Specifically, under the highest parasitic load scenario, cost savings from \ce{MgSO4}-based, \ce{MgCl2}-based, \ce{K2CO3}-based TES decreases by \$0.7 million (58\%), \$22.1 million (88\%), \$13.3 (78\%) million, and \$ 12.5 million (76\%).

\begin{table}[H]
\tablefontsize
\centering
\begin{tabular}{lccc} 
\hline 
\hline 
 Scenario & Household Cost Savings (\$) & City-level Cost Savings (\$) & Break-even Costs (\$/kWh)  \\
  & (Change from Reference) & (Change from Reference) & (Change from Reference) \\
\hline
\multicolumn{4}{l}{\textit{Reference}:}\\
 \ce{SrBr2} & \$12 - \$157 & \$16.4M & \$12.7 - \$16.7 \\
\ce{MgSO4} & \$1 - \$14 & \$1.3M & \$1.2 - \$2.4 \\
\ce{MgCl2} & \$18 - \$241 & \$25.2M & \$2.5 - \$3.3 \\
 \ce{K2CO3} & \$13 - \$163 & \$17.1M & \$3.6 - \$4.8 \\
 \hline \\
\hline 
 \multicolumn{4}{l}{\textit{Examining value of TES with incremental sizes:}}\\
\ce{SrBr2} & \$15 - \$158 (18\% - 25\%) & \$20.0M (22\%) &\$12 - \$15 (-9\% - -11\%)\\
\ce{MgSO4} & \$2 - \$27 (62\% - 69\%) & \$2.2M (67\%) &\$0.9 - \$1.8 (-26\% - -21\%)\\
\ce{MgCl2} & \$24 - \$296 (23\% - 33\%) & \$32.8M (30\%) &\$1.7 - \$1.9 (-42\% - -31\%)\\
 \ce{K2CO3} & \$22 - \$262 (61\% - 67\%) & \$28.3M (66\%) &\$3.3 - \$4.4 (-9\% - -7\%)\\
\hline \\ \hline
\multicolumn{4}{l}{\textit{Examining value of TES with fixed sizes:}}\\
\ce{SrBr2} & \$21 - \$277 (70\% - 71\%) & \$28.0M (71\%) &\$11.0 - \$14.9 (-13\% - -11\%)\\
\hline \\ \hline
\multicolumn{4}{l}{\textit{Examining value of TES with different designs:}}\\
Constant Power Rating (\ce{SrBr2}) & \$13 - \$159 (1\% - 4\%) & \$17.1M (4\%) &\$13.0 - \$17.0 (2\% - 3\%)\\
Gas Furnace Backup (\ce{SrBr2}) & \$10 - \$119 (-19\% - -24\%) & \$12.1M (-26\%) &\$10.8 - \$13.3 (-15\% - -20\%)\\
\hline \\ \hline
\multicolumn{4}{l}{\textit{Higher Humidification Load:}}\\
\ce{SrBr2} & \$11 - \$82 (-51\% - -48\%) & \$8.4M (-49\%) &\$10.2 - \$11.6 (-39\% - -9\%)\\
\ce{MgSO4} & \$0.5 - \$8 (-47\% - -42\%) & \$0.7M (-45\%) &\$0.8 - \$1.4 (-40\% - -37\%)\\
\ce{MgCl2} & \$11 - \$123 (-49\% - -39\%) & \$11.6M (-54\%) &\$1.5 - \$2.1 (-39\% - -36\%)\\
\ce{k2CO3} & \$5 - \$96 (-64\% - -41\%) & \$9.1M (-47\%) &\$2.1 - \$3.1 (-41\% - -36\%)\\
\hline \\ \hline
\multicolumn{4}{l}{\textit{Highest Humidification Load:}}\\
\ce{SrBr2} & \$3 - \$42 (-75\% - -73\%) & \$3.9M (-76\%) &\$3.4 - \$4.8 (-73\% - -71\%)\\
\ce{MgSO4} & \$0.3 - \$3 (-77\% - -65\%) & \$0.6M (-58\%) &\$0.3 - \$0.7 (-78\% - -71\%)\\
\ce{MgCl2} & \$7 - \$77 (-68\% - -60\%) & \$3.1M (-88\%) &\$0.8 - \$1.2 (-69\% - -63\%)\\
\ce{K2CO3} & \$3 - \$11 (-93\% - -79\%) & \$3.8M (-78\%) &\$0.6 - \$1.1 (-82\% - -77\%)\\
\hline \\ \hline
\multicolumn{4}{l}{\textit{5\% Other Loss from Non-Humidification Parasitic Load:}}\\
\ce{SrBr2} & \$11 - \$140 (-11\% - -8\%) & \$14.8M (-10\%) &\$11.6 - \$13.5 (-18\% - -9\%)\\
\ce{MgSO4} & \$0.5 - \$13 (-11\% - -5\%) & \$1.2M (-8\%) &\$1.1 - \$2.0 (-15\% - -6\%)\\
\ce{MgCl2} & \$17 - \$224 (-7\% - -3\%) & \$23.9M (-5\%) &\$2.4 - \$2.9 (-11\% - -4\%)\\
\ce{K2CO3} & \$12 - \$145 (-11\% - -6\%) & \$15.6M (-9\%) &\$3.5 - \$4.2 (-12\% - -7\%)\\
\hline \\ \hline
 \multicolumn{4}{l}{\textit{10\% Other Loss from Non-Humidification Parasitic Load:}}\\
\ce{SrBr2} & \$10 - \$118 (-25\% - -20\%) & \$12.8M (-22\%) &\$10.7 - \$12.9 (-23\% - -16\%)\\
\ce{MgSO4} & \$0.7 - \$11 (-24\% - -19\%) & \$1.0M (-20\%) &\$1.0 - \$1.8 (-27\% - -19\%)\\
\ce{MgCl2} & \$16 - \$200 (-17\% - -13\%) & \$21.4M (-15\%) &\$2.2 - \$2.7 (-18\% - -11\%)\\
\ce{K2CO3} & \$11 - \$126 (-23\% - -18\%) & \$13.5M (-21\%) &\$3.0 - \$3.7 (-22\% - -18\%)\\
\bottomrule
\end{tabular} 
\caption{Results from sensitivity analysis of the study for each of the 400 representative homes in Detroit, MI}
\label{tab:sa_results}
\end{table}

%% file: main.bbl
\begin{thebibliography}{}

\bibitem[\protect\citeauthoryear{Aarts, de~Jong, Cotti, Donkers, Fischer, Adan, and Huinink}{Aarts et~al.}{2022}]{Aarts2022}
Aarts, J., S.~de~Jong, M.~Cotti, P.~Donkers, H.~Fischer, O.~Adan, and H.~Huinink (2022, 3).
\newblock Diffusion limited hydration kinetics of millimeter sized salt hydrate particles for thermochemical heat storage.
\newblock {\em Journal of Energy Storage\/}~{\em 47}.

\bibitem[\protect\citeauthoryear{APS}{APS}{2023}]{APS_2023}
APS (2023).
\newblock {Time of Day}.
\newblock https://www.aps.com/en/Residential/Service-Plans/Compare-Service-Plans/Time-of-Use-4pm-7pm-Weekdays.
\newblock Accessed: 2024-3-23.

\bibitem[\protect\citeauthoryear{Arteconi, Hewitt, and Polonara}{Arteconi et~al.}{2013}]{Arteconi_and_Polonara_2013}
Arteconi, A., N.~Hewitt, and F.~Polonara (2013).
\newblock {Domestic demand-side management (DSM): Role of heat pumps and thermal energy storage (TES) systems}.
\newblock {\em Applied Thermal Engineering\/}~{\em 51\/}(1), 155--165.

\bibitem[\protect\citeauthoryear{Cammarata, Verda, Sciacovelli, and Ding}{Cammarata et~al.}{2018}]{Cammarata2018}
Cammarata, A., V.~Verda, A.~Sciacovelli, and Y.~Ding (2018, 6).
\newblock Hybrid strontium bromide-natural graphite composites for low to medium temperature thermochemical energy storage: Formulation, fabrication and performance investigation.
\newblock {\em Energy Conversion and Management\/}~{\em 166}, 233--240.

\bibitem[\protect\citeauthoryear{Christen and Carlen}{Christen and Carlen}{2000}]{Christen_Carlen_2000}
Christen, T. and M.~W. Carlen (2000).
\newblock Theory of ragone plots.
\newblock {\em Journal of Power Sources\/}~{\em 91\/}(2), 210--216.

\bibitem[\protect\citeauthoryear{Clark, Gholamibozanjani, Woods, Kaur, Odukomaiya, Al-Hallaj, and Farid}{Clark et~al.}{2022}]{Clark_et_al_2022}
Clark, R.-J., G.~Gholamibozanjani, J.~Woods, S.~Kaur, A.~Odukomaiya, S.~Al-Hallaj, and M.~Farid (2022).
\newblock Experimental screening of salt hydrates for thermochemical energy storage for building heating application.
\newblock {\em Journal of Energy Storage\/}~{\em 51}, 104415.

\bibitem[\protect\citeauthoryear{ComED}{ComED}{2023}]{ComED_2023}
ComED (2023).
\newblock {Time of Day}.
\newblock https://www.comed.com/ways-to-save/for-your-home/manage-my-energy/time-of-day-pricing.
\newblock Accessed: 2024-3-23.

\bibitem[\protect\citeauthoryear{ConED}{ConED}{2023}]{ConEd_2023}
ConED (2023).
\newblock {Time of Day}.
\newblock {https://www.ouc.com/docs/rates---elctric-water-meter/ouc\_electric\_rates\_1023.pdf}.
\newblock Accessed: 2024-3-23.

\bibitem[\protect\citeauthoryear{De~Boer, Smeding, Zondag, and Krol}{De~Boer et~al.}{2014}]{DeBoer_et_al_2014}
De~Boer, R., S.~Smeding, H.~Zondag, and G.~Krol (2014, 05).
\newblock Development of a prototype system for seasonal solar heat storage using an open sorption process.

\bibitem[\protect\citeauthoryear{{de Gracia} and Cabeza}{{de Gracia} and Cabeza}{2017}]{Alvaro_and_Cabeza_2017}
{de Gracia}, A. and L.~F. Cabeza (2017).
\newblock Numerical simulation of a pcm packed bed system: A review.
\newblock {\em Renewable and Sustainable Energy Reviews\/}~{\em 69}, 1055--1063.

\bibitem[\protect\citeauthoryear{de~Jong}{de~Jong}{2022}]{de_Jong_thesis}
de~Jong, S. (2022).
\newblock Modelling salt hydration in packed beds for thermochemical heat storage purposes, master's thesis.

\bibitem[\protect\citeauthoryear{Deetjen, Walsh, and Vaishnav}{Deetjen et~al.}{2021}]{Deetjen_2021}
Deetjen, T.~A., L.~Walsh, and P.~Vaishnav (2021, jul).
\newblock {US residential heat pumps: the private economic potential and its emissions, health, and grid impacts}.
\newblock {\em Environmental Research Letters\/}~{\em 16\/}(8), 084024.

\bibitem[\protect\citeauthoryear{DOE}{DOE}{2017}]{DOE_2017}
DOE (2017).
\newblock {2017 Building America climate-specific guidance}.
\newblock \url{https://www.energy.gov/eere/buildings/climate-zones}, institution = {Office of Energy Efficiency and Renewable Energy}, note = {Accessed: 2023-02-10}.

\bibitem[\protect\citeauthoryear{DOE}{DOE}{2022}]{DOE_2022}
DOE (2022).
\newblock {Buildings Energy Efficiency Frontiers \& Innovation Technologies (BENEFIT) – 2022/23}.
\newblock Technical report, DOE.

\bibitem[\protect\citeauthoryear{Donkers, Sögütoglu, Huinink, Fischer, and Adan}{Donkers et~al.}{2017}]{Donkers_2017}
Donkers, P., L.~Sögütoglu, H.~Huinink, H.~Fischer, and O.~Adan (2017).
\newblock A review of salt hydrates for seasonal heat storage in domestic applications.
\newblock {\em Applied Energy\/}~{\em 199}, 45--68.

\bibitem[\protect\citeauthoryear{DTE}{DTE}{2023}]{DTE_2023}
DTE (2023).
\newblock {Time of Day}.
\newblock https://solutions.dteenergy.com/dte/en/Products/Time-of-Day-3-p-m---7-p-m-/p/TOD-3-7.
\newblock Accessed: 2023-3-27.

\bibitem[\protect\citeauthoryear{D’Ettorre, {De Rosa}, Conti, Testi, and Finn}{D’Ettorre et~al.}{2019}]{DEttorre_et_al_2019}
D’Ettorre, F., M.~{De Rosa}, P.~Conti, D.~Testi, and D.~Finn (2019).
\newblock Mapping the energy flexibility potential of single buildings equipped with optimally-controlled heat pump, gas boilers and thermal storage.
\newblock {\em Sustainable Cities and Society\/}~{\em 50}, 101689.

\bibitem[\protect\citeauthoryear{EIA}{EIA}{2022}]{EIA_2022}
EIA (2022).
\newblock {Annual Energy Outlook 2022, Table 4. Residential Sector Key Indicators and Consumption}.
\newblock Technical report, EIA.

\bibitem[\protect\citeauthoryear{Eversource}{Eversource}{2023}]{Eversource_2023}
Eversource (2023).
\newblock {Time of Day}.
\newblock https://www.eversource.com/content/docs/default-source/rates-tariffs/nh-summary-rates.pdf.
\newblock Accessed: 2024-3-23.

\bibitem[\protect\citeauthoryear{Farah, Liu, and Saman}{Farah et~al.}{2019}]{Farah_2019}
Farah, S., M.~Liu, and W.~Saman (2019).
\newblock Numerical investigation of phase change material thermal storage for space cooling.
\newblock {\em Applied Energy\/}~{\em 239}, 526--535.

\bibitem[\protect\citeauthoryear{Fischer, Toral, Lindberg, Wille-Haussmann, and Madani}{Fischer et~al.}{2014}]{Fischer_2014}
Fischer, D., T.~R. Toral, K.~Lindberg, B.~Wille-Haussmann, and H.~Madani (2014).
\newblock Investigation of thermal storage operation strategies with heat pumps in german multi family houses.
\newblock {\em Energy Procedia\/}~{\em 58}, 137--144.
\newblock Renewable Energy Research Conference, RERC 2014.

\bibitem[\protect\citeauthoryear{Fisher, Ding, and Sciacovelli}{Fisher et~al.}{2021}]{Fisher2021}
Fisher, R., Y.~Ding, and A.~Sciacovelli (2021, 6).
\newblock Hydration kinetics of k2co3, mgcl2 and vermiculite-based composites in view of low-temperature thermochemical energy storage.
\newblock {\em Journal of Energy Storage\/}~{\em 38}, 102561.

\bibitem[\protect\citeauthoryear{Fopah-Lele and Tamba}{Fopah-Lele and Tamba}{2017}]{Lele_and_Tamba_2017}
Fopah-Lele, A. and J.~G. Tamba (2017).
\newblock A review on the use of \ce{SrBr2}·\ce{6H2O} as a potential material for low temperature energy storage systems and building applications.
\newblock {\em Solar Energy Materials and Solar Cells\/}~{\em 164}, 175--187.

\bibitem[\protect\citeauthoryear{Gaeini, Shaik, and Rindt}{Gaeini et~al.}{2019}]{Gaeini2019}
Gaeini, M., S.~A. Shaik, and C.~Rindt (2019, 8).
\newblock Characterization of potassium carbonate salt hydrate for thermochemical energy storage in buildings.
\newblock {\em Energy and Buildings\/}~{\em 196}, 178--193.

\bibitem[\protect\citeauthoryear{GW}{GW}{2023}]{GW_2023}
GW (2023).
\newblock {Time of Day}.
\newblock https://www.georgiapower.com/content/dam/georgia-power/pdfs/residential-pdfs/tariffs/2024/tou-rd-9.pdf.
\newblock Accessed: 2024-3-23.

\bibitem[\protect\citeauthoryear{Hawwash, Hassan, and feky}{Hawwash et~al.}{2020}]{Hawwash2020}
Hawwash, A.~A., H.~Hassan, and K.~E. feky (2020, 3).
\newblock Impact of reactor design on the thermal energy storage of thermochemical materials.
\newblock {\em Applied Thermal Engineering\/}~{\em 168}, 114776.

\bibitem[\protect\citeauthoryear{IEA}{IEA}{2022a}]{IEA_2022c}
IEA (2022a).
\newblock { Heat Pumps}.
\newblock Technical report, IEA.

\bibitem[\protect\citeauthoryear{IEA}{IEA}{2022b}]{IEA_2022b}
IEA (2022b).
\newblock {Monthly Energy Review - September 2022}.
\newblock Technical report, IEA.

\bibitem[\protect\citeauthoryear{IEM}{IEM}{}]{IEM}
IEM.
\newblock {Iowa Environmental Mesonet}.
\newblock \url{https://mesonet.agron.iastate.edu/request/download.phtml}, note = {Accessed: 2022-12-29}.

\bibitem[\protect\citeauthoryear{James, Mahvi, and Woods}{James et~al.}{2022}]{James2022}
James, N., A.~Mahvi, and J.~Woods (2022, 12).
\newblock Optimizing phase change composite thermal energy storage using the thermal ragone framework.
\newblock {\em Journal of Energy Storage\/}~{\em 56}, 105875.

\bibitem[\protect\citeauthoryear{Kenisarin and Mahkamov}{Kenisarin and Mahkamov}{2016}]{Kenisarin_and_Mahkamov_2016}
Kenisarin, M. and K.~Mahkamov (2016).
\newblock Salt hydrates as latent heat storage materials:thermophysical properties and costs.
\newblock {\em Solar Energy Materials and Solar Cells\/}~{\em 145}, 255--286.

\bibitem[\protect\citeauthoryear{Kinzer, Ghosh, Crowley, Jatkar, and Chandran}{Kinzer et~al.}{2024}]{Kinzer_2023}
Kinzer, B., D.~Ghosh, D.~Crowley, A.~Jatkar, and R.~B. Chandran (2024).
\newblock {Modeling and Experimental Demonstrations Reveal Ragone Framework for Salt Hydrate Thermochemical Energy Storage.}
\newblock {\em ChemRxiv\/}.

\bibitem[\protect\citeauthoryear{Kumar, Ness, Chavez, Banerjee, Muley, and Stoia}{Kumar et~al.}{2020}]{Kumar_2020}
Kumar, N., R.~V. Ness, J.~Chavez, Reynaldo, D.~Banerjee, A.~Muley, and M.~Stoia (2020, 08).
\newblock {Experimental Analysis of Salt Hydrate Latent Heat Thermal Energy Storage System With Porous Aluminum Fabric and Salt Hydrate as Phase Change Material With Enhanced Stability and Supercooling}.
\newblock {\em Journal of Energy Resources Technology\/}~{\em 143\/}(4).
\newblock 042001.

\bibitem[\protect\citeauthoryear{LADWP}{LADWP}{2023}]{LADWP_2023}
LADWP (2023).
\newblock {Time of Day}.
\newblock https://www.ladwp.com/account/customer-service/understanding-your-rates/residential-electric-rates.
\newblock Accessed: 2024-3-23.

\bibitem[\protect\citeauthoryear{Langevin, Harris, and Reyna}{Langevin et~al.}{2019}]{Langevin_2019}
Langevin, J., C.~B. Harris, and J.~L. Reyna (2019).
\newblock Assessing the potential to reduce u.s. building co2 emissions 80% by 2050.
\newblock {\em Joule\/}~{\em 3\/}(10), 2403--2424.

\bibitem[\protect\citeauthoryear{Le, Huang, Shah, Wilson, Artain, Byrne, and Hewitt}{Le et~al.}{2019}]{Le_et_al_2019}
Le, K.~X., M.~J. Huang, N.~N. Shah, C.~Wilson, P.~M. Artain, R.~Byrne, and N.~J. Hewitt (2019).
\newblock Techno-economic assessment of cascade air-to-water heat pump retrofitted into residential buildings using experimentally validated simulations.
\newblock {\em Applied Energy\/}~{\em 250}, 633--652.

\bibitem[\protect\citeauthoryear{Li, Huang, Xu, Liu, and Wu}{Li et~al.}{2018}]{Li_2018}
Li, Y., G.~Huang, T.~Xu, X.~Liu, and H.~Wu (2018).
\newblock Optimal design of pcm thermal storage tank and its application for winter available open-air swimming pool.
\newblock {\em Applied Energy\/}~{\em 209}, 224--235.

\bibitem[\protect\citeauthoryear{Linnow, Niermann, Bonatz, Posern, and Steiger}{Linnow et~al.}{2014}]{Linnow2014}
Linnow, K., M.~Niermann, D.~Bonatz, K.~Posern, and M.~Steiger (2014).
\newblock Experimental studies of the mechanism and kinetics of hydration reactions.
\newblock {\em Energy Procedia\/}~{\em 48}, 394--404.

\bibitem[\protect\citeauthoryear{Love, Smith, Watson, Oikonomou, Summerfield, Gleeson, Biddulph, Chiu, Wingfield, Martin, Stone, and Lowe}{Love et~al.}{2017}]{Love_et_al_2017}
Love, J., A.~Z. Smith, S.~Watson, E.~Oikonomou, A.~Summerfield, C.~Gleeson, P.~Biddulph, L.~F. Chiu, J.~Wingfield, C.~Martin, A.~Stone, and R.~Lowe (2017).
\newblock The addition of heat pump electricity load profiles to gb electricity demand: Evidence from a heat pump field trial.
\newblock {\em Applied Energy\/}~{\em 204}, 332--342.

\bibitem[\protect\citeauthoryear{Lun and Tung}{Lun and Tung}{2020}]{Lun_and_Tung_2020}
Lun, Y. H.~V. and S.~L.~D. Tung (2020).
\newblock {\em Heat Pumps for Sustainable Heating and Cooling}.
\newblock Springer Cham.

\bibitem[\protect\citeauthoryear{Mahmoudi, Donkers, Walayat, Peters, and Shahi}{Mahmoudi et~al.}{2021}]{Mahmoudi2021}
Mahmoudi, A., P.~A. Donkers, K.~Walayat, B.~Peters, and M.~Shahi (2021, 12).
\newblock A thorough investigation of thermochemical heat storage system from particle to bed scale.
\newblock {\em Chemical Engineering Science\/}~{\em 246}, 116877.

\bibitem[\protect\citeauthoryear{Masip, Cazorla-Marín, Montagud-Montalvá, Marchante, Barceló, and Corberán}{Masip et~al.}{2019}]{Masip_2019}
Masip, X., A.~Cazorla-Marín, C.~Montagud-Montalvá, J.~Marchante, F.~Barceló, and J.~Corberán (2019).
\newblock Energy and techno-economic assessment of the effect of the coupling between an air source heat pump and the storage tank for sanitary hot water production.
\newblock {\em Applied Thermal Engineering\/}~{\em 159}, 113853.

\bibitem[\protect\citeauthoryear{Michel, Mazet, Mauran, Stitou, and Xu}{Michel et~al.}{2012}]{Michel2012}
Michel, B., N.~Mazet, S.~Mauran, D.~Stitou, and J.~Xu (2012, 11).
\newblock Thermochemical process for seasonal storage of solar energy: Characterization and modeling of a high density reactive bed.
\newblock {\em Energy\/}~{\em 47}, 553--563.

\bibitem[\protect\citeauthoryear{Moreno, Solé, Castell, and Cabeza}{Moreno et~al.}{2014}]{Moreno_2014}
Moreno, P., C.~Solé, A.~Castell, and L.~F. Cabeza (2014).
\newblock The use of phase change materials in domestic heat pump and air-conditioning systems for short term storage: A review.
\newblock {\em Renewable and Sustainable Energy Reviews\/}~{\em 39}, 1--13.

\bibitem[\protect\citeauthoryear{N’Tsoukpoe and Kuznik}{N’Tsoukpoe and Kuznik}{2021}]{NTsoukpoe_Kuznik_2021}
N’Tsoukpoe, K.~E. and F.~Kuznik (2021).
\newblock A reality check on long-term thermochemical heat storage for household applications.
\newblock {\em Renewable and Sustainable Energy Reviews\/}~{\em 139}, 110683.

\bibitem[\protect\citeauthoryear{Odukomaiya, Woods, James, Kaur, Gluesenkamp, Kumar, Mumme, Jackson, and Prasher}{Odukomaiya et~al.}{2021}]{Odukomaiya_2021}
Odukomaiya, A., J.~Woods, N.~James, S.~Kaur, K.~R. Gluesenkamp, N.~Kumar, S.~Mumme, R.~Jackson, and R.~Prasher (2021).
\newblock Addressing energy storage needs at lower cost via on-site thermal energy storage in buildings.
\newblock {\em Energy Environ. Sci.\/}~{\em 14}, 5315--5329.

\bibitem[\protect\citeauthoryear{Osterman and Stritih}{Osterman and Stritih}{2021}]{Osterman_and_Stritih_2021}
Osterman, E. and U.~Stritih (2021).
\newblock Review on compression heat pump systems with thermal energy storage for heating and cooling of buildings.
\newblock {\em Journal of Energy Storage\/}~{\em 39}, 102569.

\bibitem[\protect\citeauthoryear{OUC}{OUC}{2023}]{OUC_2023}
OUC (2023).
\newblock {Time of Day}.
\newblock https://www.coned.com/en/accounts-billing/your-bill/time-of-use.
\newblock Accessed: 2024-3-23.

\bibitem[\protect\citeauthoryear{RCW}{RCW}{}]{RCW}
RCW.
\newblock {An Interactive Visualization of the Housing Characteristic Dependencies in ResStock}.
\newblock \url{https://htmlpreview.github.io/?https://github.com/NREL/OpenStudio-BuildStock/blob/master/project_national/util/dependency_wheel/dep_wheel.html}.
\newblock Accessed: 2022-12-29.

\bibitem[\protect\citeauthoryear{Renaldi, Kiprakis, and Friedrich}{Renaldi et~al.}{2017}]{Renaldi_2017}
Renaldi, R., A.~Kiprakis, and D.~Friedrich (2017).
\newblock An optimisation framework for thermal energy storage integration in a residential heat pump heating system.
\newblock {\em Applied Energy\/}~{\em 186}, 520--529.

\bibitem[\protect\citeauthoryear{SCL}{SCL}{2023}]{SCL_2023}
SCL (2023).
\newblock {Time of Day}.
\newblock \url{https://powerlines.seattle.gov/wp-content/uploads/sites/17/2022/09/2023-and-2024-Retail-Rate-Summary_final.pdf}, note = {Accessed: 2024-3-23}.

\bibitem[\protect\citeauthoryear{Teamah and Lightstone}{Teamah and Lightstone}{2019}]{Teamah_and_Lightstone_2019}
Teamah, H. and M.~Lightstone (2019).
\newblock Numerical study of the electrical load shift capability of a ground source heat pump system with phase change thermal storage.
\newblock {\em Energy and Buildings\/}~{\em 199}, 235--246.

\bibitem[\protect\citeauthoryear{Teichert, Link, Schneider, Wolff, and Lienkamp}{Teichert et~al.}{2023}]{Teichert_et_al_2023}
Teichert, O., S.~Link, J.~Schneider, S.~Wolff, and M.~Lienkamp (2023).
\newblock Techno-economic cell selection for battery-electric long-haul trucks.
\newblock {\em eTransportation\/}~{\em 16}, 100225.

\bibitem[\protect\citeauthoryear{Vaishnav and Fatimah}{Vaishnav and Fatimah}{2020}]{Vaishnav_et_al_2020}
Vaishnav, P. and A.~M. Fatimah (2020).
\newblock The environmental consequences of electrifying space heating.
\newblock {\em Environmental Science \& Technology\/}~{\em 54\/}(16), 9814--9823.
\newblock PMID: 32648744.

\bibitem[\protect\citeauthoryear{Vykhodtsev, Jang, Wang, Rosehart, and Zareipour}{Vykhodtsev et~al.}{2022}]{Vykhodtsev_2022}
Vykhodtsev, A.~V., D.~Jang, Q.~Wang, W.~Rosehart, and H.~Zareipour (2022).
\newblock A review of modelling approaches to characterize lithium-ion battery energy storage systems in techno-economic analyses of power systems.
\newblock {\em Renewable and Sustainable Energy Reviews\/}~{\em 166}, 112584.

\bibitem[\protect\citeauthoryear{Waite and Modi}{Waite and Modi}{2020}]{Waite_and_Modi_2022}
Waite, M. and V.~Modi (2020).
\newblock Electricity load implications of space heating decarbonization pathways.
\newblock {\em Joule\/}~{\em 4}, 376--394.

\bibitem[\protect\citeauthoryear{Wang, Yu, and Ling}{Wang et~al.}{2019}]{Wang_2019b}
Wang, Y., K.~Yu, and X.~Ling (2019).
\newblock Experimental and modeling study on thermal performance of hydrated salt latent heat thermal energy storage system.
\newblock {\em Energy Conversion and Management\/}~{\em 198}, 111796.

\bibitem[\protect\citeauthoryear{Wang, Yu, and Ling}{Wang et~al.}{2020}]{Wang_2020}
Wang, Y., K.~Yu, and X.~Ling (2020).
\newblock Experimental study on thermal performance of a mobilized thermal energy storage system: A case study of hydrated salt latent heat storage.
\newblock {\em Energy and Buildings\/}~{\em 210}, 109744.

\bibitem[\protect\citeauthoryear{Wang, Wang, Ma, Lin, and Ren}{Wang et~al.}{2019}]{Wang_2019}
Wang, Z., F.~Wang, Z.~Ma, W.~Lin, and H.~Ren (2019).
\newblock Investigation on the feasibility and performance of transcritical co2 heat pump integrated with thermal energy storage for space heating.
\newblock {\em Renewable Energy\/}~{\em 134}, 496--508.

\bibitem[\protect\citeauthoryear{Wilson}{Wilson}{2017}]{Wilson_2017}
Wilson, E.~J. (2017, 9).
\newblock Resstock - targeting energy and cost savings for u.s. homes.
\newblock {\em NREL\/}.

\bibitem[\protect\citeauthoryear{Woods, Mahvi, Goyal, Kozubal, Odukomaiya, and Jackson}{Woods et~al.}{2021}]{Woods_et_al_2021}
Woods, J., A.~Mahvi, A.~Goyal, E.~Kozubal, A.~Odukomaiya, and R.~Jackson (2021).
\newblock Rate capability and ragone plots for phase change thermal energy storage.
\newblock {\em Nature Energy\/}~{\em 6}, 295–302.

\bibitem[\protect\citeauthoryear{Wu, Wang, Li, Xu, Yang, and Yang}{Wu et~al.}{2020}]{Wu_et_al_2020}
Wu, P., Z.~Wang, X.~Li, Z.~Xu, Y.~Yang, and Q.~Yang (2020).
\newblock Energy-saving analysis of air source heat pump integrated with a water storage tank for heating applications.
\newblock {\em Building and Environment\/}~{\em 180}, 107029.

\bibitem[\protect\citeauthoryear{Wu and Skye}{Wu and Skye}{2021}]{Wu_and_Skye_2021}
Wu, W. and H.~M. Skye (2021).
\newblock Residential net-zero energy buildings: Review and perspective.
\newblock {\em Renewable and Sustainable Energy Reviews\/}~{\em 142}, 110859.

\bibitem[\protect\citeauthoryear{Xcel}{Xcel}{2023a}]{Xcel_CO_2023}
Xcel (2023a).
\newblock {Time of Day}.
\newblock https://co.my.xcelenergy.com/s/billing-payment/residential-rates/time-of-use-pricing.
\newblock Accessed: 2024-3-23.

\bibitem[\protect\citeauthoryear{Xcel}{Xcel}{2023b}]{Xcel_TX_2023}
Xcel (2023b).
\newblock {Time of Day}.
\newblock https://www.xcelenergy.com/staticfiles/xe-responsive/Programs%20and%20Rebates/Business/TX-Time-of-use-rate-FAQ.pdf.
\newblock Accessed: 2024-3-23.

\bibitem[\protect\citeauthoryear{Xcel}{Xcel}{2023c}]{Xcel_MN_2023}
Xcel (2023c).
\newblock {Time of Day}.
\newblock https://www.xcelenergy.com/staticfiles/xe-responsive/Company/Rates%20&%20Regulations/Me_Section_5.pdf.
\newblock Accessed: 2024-3-23.

\bibitem[\protect\citeauthoryear{Yan, Li, Xu, Chao, Wang, Aristov, Gordeeva, Dutta, and Murthy}{Yan et~al.}{2021}]{Yan2021}
Yan, T., T.~Li, J.~Xu, J.~Chao, R.~Wang, Y.~I. Aristov, L.~G. Gordeeva, P.~Dutta, and S.~S. Murthy (2021, 5).
\newblock Ultrahigh-energy-density sorption thermal battery enabled by graphene aerogel-based composite sorbents for thermal energy harvesting from air.
\newblock {\em ACS Energy Letters\/}~{\em 6}, 1795--1802.

\bibitem[\protect\citeauthoryear{Yazawa, Shamberger, and Fisher}{Yazawa et~al.}{2019}]{Yazawa2019}
Yazawa, K., P.~J. Shamberger, and T.~S. Fisher (2019).
\newblock Ragone relations for thermal energy storage technologies.
\newblock {\em Frontiers in Mechanical Engineering\/}~{\em 5}.

\bibitem[\protect\citeauthoryear{Ye, Liu, Wang, Liu, Lv, and Yang}{Ye et~al.}{2022}]{Ye_et_al_2022}
Ye, Z., H.~Liu, W.~Wang, H.~Liu, J.~Lv, and F.~Yang (2022).
\newblock Reaction/sorption kinetics of salt hydrates for thermal energy storage.
\newblock {\em Journal of Energy Storage\/}~{\em 56}, 106122.

\bibitem[\protect\citeauthoryear{Yu, Li, Zhang, and Zhao}{Yu et~al.}{2021}]{Yu_2021}
Yu, M., S.~Li, X.~Zhang, and Y.~Zhao (2021).
\newblock Techno-economic analysis of air source heat pump combined with latent thermal energy storage applied for space heating in china.
\newblock {\em Applied Thermal Engineering\/}~{\em 185}, 116434.

\bibitem[\protect\citeauthoryear{Zhang, Dong, Wang, and Feng}{Zhang et~al.}{2020}]{Zhang2020}
Zhang, Y., H.~Dong, R.~Wang, and P.~Feng (2020, 5).
\newblock Air humidity assisted sorption thermal battery governed by reaction wave model.
\newblock {\em Energy Storage Materials\/}~{\em 27}, 9--16.

\bibitem[\protect\citeauthoryear{Zhao and Wang}{Zhao and Wang}{2019}]{Zhao_and_Wang_2019}
Zhao, B. and R.~Wang (2019).
\newblock Perspectives for short-term thermal energy storage using salt hydrates for building heating.
\newblock {\em Energy\/}~{\em 189}, 116139.

\bibitem[\protect\citeauthoryear{Zhu, Hu, Lei, Jiang, and Lei}{Zhu et~al.}{2015}]{Zhu_et_al_2015}
Zhu, N., P.~Hu, Y.~Lei, Z.~Jiang, and F.~Lei (2015).
\newblock Numerical study on ground source heat pump integrated with phase change material cooling storage system in office building.
\newblock {\em Applied Thermal Engineering\/}~{\em 87}, 615--623.

\bibitem[\protect\citeauthoryear{Ángel Á.~Pardiñas, Alonso, Diz, Kvalsvik, and Fernández-Seara}{Ángel Á.~Pardiñas et~al.}{2017}]{Pardinas_et_al_2017}
Ángel Á.~Pardiñas, M.~J. Alonso, R.~Diz, K.~H. Kvalsvik, and J.~Fernández-Seara (2017).
\newblock State-of-the-art for the use of phase-change materials in tanks coupled with heat pumps.
\newblock {\em Energy and Buildings\/}~{\em 140}, 28--41.

\end{thebibliography}
